\documentstyle[12pt]{article}

\newlength{\abstractwidth}
\renewcommand{\thefootnote}{\fnsymbol{footnote}}
\renewcommand{\thanks}[1]{\footnote{#1}} 
\newcommand{\starttext}{ \setcounter{footnote}{0}
\renewcommand{\thefootnote}{\arabic{footnote}}}

\newcommand{\be}{\begin{equation}}
\newcommand{\bea}{\begin{eqnarray}}
\newcommand{\eea}{\end{eqnarray}} \newcommand{\ee}{\end{equation}}
 
\renewcommand{\>}{\rangle}
\def\ba{\begin{eqnarray}}
\def\ea{\end{eqnarray}}

\def\B{{\cal B}}

\def\E{{\cal E}}
\def\K{{\cal K}}
\def\cM{{\cal M}}

\def\r{\rho}

\def\ra{\rightarrow}

\def\o{\omega}
\def\Re{{\rm Re}}

\def\det{{\rm det}}

\def\log{\,{\rm log}\,}

\def\o{\omega}

\def\a{\alpha}
\def\al{\alpha}
\def\b{\beta}
\def\g{\gamma}
\def\d{\delta}
\def\e{\varepsilon}

\def\m{\mu}
\def\n{\nu}
\def\o{\omega}

\def\r{\rho}

\def\t{\theta}

\def\O{\Omega}

\def\na{\nabla}
\def\p{\partial}

\def\ti{\tilde}
\def\vp{\varphi}
\def\vps{\underline{\varphi}}

\def\Z{{\bf Z}}

\def\R{{\bf R}}
\def\C{{\bf C}}

\def\i{\infty}
\def\I{\int}
\def\s{\sum}

\def\ddb{{\partial\bar\partial}}
\def\ddbar{{\partial\bar\partial}}

\def\phi{\varphi}

\def\sub{\subseteq}
\def\ra{\rightarrow}

\def\B{{\cal B}}
\def\F{{\cal F}}

\def\cO{{\cal O}}
\def\cM{{\cal M}}

\def\cC{{\cal C}}

\def\cO{{\cal O}}

\def\na{{\nabla}}

\def\K{{K\"ahler\ }}

 \def\v{\vskip .1in}

\def\[{{\bf [}}
\def\]{{\bf ]}}

\def\pl{\partial}

\def\cT{{\cal T}}
\def\PSH{{ PSH}}
\def\Ricci{{ Ricci}}

\begin{document}
\starttext \baselineskip=15pt \setcounter{footnote}{0}
\newtheorem{theorem}{Theorem}
\newtheorem{lemma}{Lemma}
\newtheorem{definition}{Definition}
\newtheorem{proposition}{Proposition}
\newtheorem{corollary}{Corollary}

\begin{center}
{\Large \bf
COMPLEX MONGE-AMP\`ERE EQUATIONS
\footnote{Contribution to the proceedings of the Journal
of Differential Geometry Conference in honor of Professor C.C. Hsiung, Lehigh University, May 2010.
Work supported in part by
National Science Foundation grants DMS-07-57372,
DMS-09-05873, and DMS-08-47524.}}
\bigskip\bigskip

{\large \bf D.H. Phong*, Jian Song$^\dagger$, and Jacob Sturm$^\ddagger$} \\

\bigskip
{}
\smallskip
$^*$ Department of Mathematics\\
Columbia University, New York, NY 10027\\

\medskip
$\dagger$ Department of Mathematics\\
Rutgers University, Piscataway, NJ 08854\\

\medskip
$\ddagger$ Department of Mathematics\\
Rutgers University, Newark, NJ 07102\\

\end{center}

\medskip

\begin{abstract}

This is a survey of some of the recent developments in the theory of complex Monge-Amp\`ere equations. The topics discussed include refinements and simplifications of classical a priori estimates, methods from pluripotential theory, variational methods for big cohomology classes, semiclassical constructions of solutions of homogeneous equations, and envelopes.

\end{abstract}

\newpage

\tableofcontents

\newpage

\section{Introduction}
\setcounter{equation}{0}

Monge-Amp\`ere equations are second-order partial differential equations whose leading term is the determinant of the Hessian of a real unknown function $\vp$.
As such, they are arguably the most basic of fully non-linear equations. 
The Hessian is required to be positive or at least non-negative, so the equations are elliptic or degenerate elliptic.
Monge-Amp\`ere equations can be divided into real
or complex, depending on whether $\vp$ is defined on a real or complex manifold.  In the real case, the Hessian is $\nabla_j\nabla_k \vp$, so the positivity of the Hessian is a convexity condition.
In the complex case, the Hessian is $\p_j\p_{\bar k}\vp$, and its positivity is rather  a plurisubharmonicity condition. Unlike convex functions, plurisubharmonic functions can have singularities, and this accounts for many significant differences between the theories of real and complex Monge-Amp\`ere equations. In these lectures, we shall concentrate on the complex case.

\medskip

The foundations of an existence and regularity theory for complex Monge-Amp\`ere equations in the elliptic case, with smooth data, were laid by Yau \cite{Y78}
and Caffarelli, Kohn, Nirenberg, and Spruck \cite{CNS, CKNS}. In \cite{Y78},
a complete solution was given for the Calabi conjecture, which 
asserts the existence of a smooth solution
to the equation
\bea
\label{Calabiconjecture}
(\o_0+{i\over 2}\ddb\vp)^n=
e^{f(z)}\,\o_0^n, 
\eea
on a compact $n$-dimensional
K\"ahler manifold $(X,\o_0)$ without boundary,
where $f(z)$ is a given smooth function satisfying the necessary condition
$\int_X e^f\o_0^n=\int_X\o_0^n$. The solution was by the method of continuity,
and the key estimates for the $C^0$ norms of $\vp$, $\Delta\vp$, and
$\na_j\na_{\bar k}\na_l\vp$ were formulated and derived there. In \cite{CKNS},
a complete solution was given for the Dirichlet problem
\bea
{\rm det}(\p_j\p_{\bar k}\vp)=F(z,\vp)\ \ {\rm on}\ \ D,
\qquad
\vp=\vp_b\ \ {\rm on}\ \ \p D,
\eea
where $D$ is a smooth, bounded, strongly pseudoconvex domain in ${\bf C}^n$,
$F\in C^\infty(\bar D\times {\bf R})$, $F(z,\vp)>0$, $F_\vp(z,\vp)\geq 0$, and $\vp_b\in C^\infty(\p D)$. A crucial ingredient of the existence and regularity developed there is the $C^0$ boundary estimates for the second order derivatives and their modulus of continuity.

\medskip
In his paper \cite{Y78}, Yau also began an existence and regularity theory for singular complex Monge-Amp\`ere equations on K\"ahler manifolds. Here the term ``singular" should be interpreted in a broad sense. It encompasses situations where the right hand side may be degenerate or have singularities \cite{Y78}, or where the manifold $X$ may not be compact or have singularities
\cite{CY80, MY83, CY86, TY86}, or where the boundary condition may be infinite
\cite{CY80}. Such extensions were required by geometric applications, and many important results were obtained, of which the references we just gave are just a small sample (see e.g. \cite{TY90, TY91, K83, W08, LYZ}, and especially \cite{Y93,Y94, Y96} and references therein). 

\medskip
The last fifteen years or so have witnessed remarkable progresses in the theory of singular Monge-Amp\`ere equations. A particularly strong impetus was provided by related problems from the minimal model program in algebraic geometry
(\cite{EGZ, ST08, TZ, BEGZ}) and from the problem of finding metrics of constant scalar curvature in a given K\"ahler class (see \cite{Y93, T97, D02}
and \cite{PS08} for a survey).
The solutions in these problems are often inherently singular, and thus they must be understood in a generalized sense. The foundations of a theory of generalized solutions for the complex Monge-Amp\`ere equation - or pluripotential theory - had been laid out by Bedford and Taylor in \cite{BT76, BT82}. There they constructed Monge-Amp\`ere measures for bounded potentials and capacities, established monotonicity theorems for their convergence, and obtained generalized solutions of the Dirichlet problem for degenerate right hand sides by the Perron method. A key catalyst for several of the recent progresses is the theorem of Kolodziej \cite{K98}, based on pluripotential theory, which provided $C^0$ estimates for Monge-Amp\`ere equations with right hand sides in $L^p$ for any $p>1$. Other important ingredients have been the extensions of pluripotential theory to unbounded potentials (\cite{GZ, BBGZ, Ceg, B06, CG} and references therein), the Tian-Yau-Zelditch theorem \cite{Y93, T90a, Z, Cat, L} on approximations of smooth metrics by Fubini-Study metrics \cite{PS06, PS07, PS09b, SZ07, SZ10, RZ08, RZ10a, RZ10b}, and refinements and extensions \cite{Gb, B09b, PS09a, PS09c, GL, Gp, Ch, TW1, TW2, DK} of the classic estimates in \cite{Y78, CKNS}. 

\medskip
The main goal of this paper is to survey some of the recent progresses. There have been many of them, and the theory is still in full flux.
While definitive answers may not yet be available to many questions, we thought it would be useful to gather here in one place, for the convenience of students and newcomers to the field, some of what is known. It was not possible to be comprehensive, and our selection of material necessarily reflects our own limitations. At the same time, we hope that the survey would be useful to a broad audience of people with relatively little familiarity with complex Monge-Amp\`ere
equations, and we have provided reasonably complete derivations in places,
when the topics are of particular importance or the literature not easily accessible.  
Each of us has lectured on parts of this paper at our home institutions, and at various workshops. In particular, the first-named author spoke at the 2010 conference at Lehigh University in honor of Professor C.C. Hsiung, one of the founders of the Journal of Differential Geometry. We would like to contribute this paper to the volume in his honor.

\section{Some General Perspective}
\setcounter{equation}{0}

Let $(X,\o_0)$ be a compact K\"ahler manifold. We consider complex Monge-Amp\`ere equations of the form
\bea
\label{CY}
(\o_0+{i\over 2}\ddb \vp)^n=F(z,\vp)\,\o_0^n
\eea
where $F(z,\vp)$ is a non-negative function. The solution $\vp$ is required to be $\o_0$-plurisub\-
harmonic,
that is, $\vp\in \PSH(X,\o_0)$, with $$PSH(X,\o_0)=\{\vp: X \rightarrow [-\infty, \infty); \vp~{\rm is ~upper\ semicontinuous}, \o_\vp\equiv\o_0+{i\over 2}\ddb\vp\geq 0\}.$$
We shall consider both the case of $X$ compact without boundary,
and the case of $\bar X$ compact with smooth boundary $\p X$, in which case we also impose a Dirichlet condition $\vp=\vp_b$, where $\vp_b\in C^\infty(\p X)$ is a given function. 

\subsection{Geometric interpretation}

Equations of the form (\ref{CY}) are fundamentally geometric in nature. The form $\o_\vp$ can be viewed as a form in the same cohomology class as $\o_0$. It defines a regular K\"ahler metric when it is $>0$, or a K\"ahler metric with degeneracies when it does have zeroes. It is well-known that the Ricci curvature form $Ricci(\o_\vp)$ of a K\"ahler form $\o_\vp$ is given by
\bea
\label{Ricci}
Ricci(\o_\vp)=-{i\over 2}\ddb \log \o_\vp^n.
\eea
Thus the equation (\ref{CY}) is just an equation for a possibly degenerate metric $\o_\vp$ in the same K\"ahler class as $\o_0$, satisfying a given constraint on its volume form $\o_\vp^n$ or, equivalently upon differentiation, a given constraint on its Ricci curvature $Ricci(\o_\vp)$.

\medskip
The modern theory of complex Monge-Amp\`ere equations began with the following two fundamental theorems, due respectively to Yau \cite{Y78}
and to Yau \cite{Y78} and Aubin \cite{A}.

\begin{theorem}
\label{Yau}
Let $(X,\o_0)$ be a compact K\"ahler manifold without boundary, and let $F(z)=e^{f(z)}$,
where $f(z)$ is a smooth function
satisfying the condition
\bea
\int_X e^f\o_0^n=\int_X\o_0^n.
\eea
Then the equation (\ref{CY}) admits a smooth solution $\vp\in \PSH(X,\o_0)$, unique up to an additive constant.
\end{theorem}

\noindent

\begin{theorem}
\label{YauAubin}
Let $(X,\o_0)$ be a compact K\"ahler manifold without boundary, and let $F(z,\vp)=e^{f+\vp}$ 
where $f(z)$ is a smooth function.
Then the equation (\ref{CY})
admits a unique smooth solution $\vp\in \PSH(X,\o_0)$.
\end{theorem}

Geometrically, Theorem \ref{Yau} provides a solution of the Calabi conjecture, which asserts that, on a compact K\"ahler manifold $X$ with $c_1(X)=0$,
there is a unique metric $\o_\vp$ with $Ricci(\o_\vp)=0$
in any K\"ahler class $[\o_0]$.
Indeed, the formula (\ref{Ricci}) shows that the Ricci form of any K\"ahler metric must be in $c_1(X)$. The assumption that $c_1(X)=0$ implies that $Ricci(\o_0)={i\over 2}\ddb f$ for some smooth function $f(z)$. It is now readily verified, by taking 
$F(z)=e^{f(z)}$ in the equation (\ref{CY}) and taking $\frac{i}{2}\ddb$ of both sides, 
that the solution of (\ref{CY}) satisfies the condition 
\bea
Ricci(\o_\vp)=0.
\eea

Similarly, Theorem \ref{YauAubin} implies the existence of a K\"ahler-Einstein metric with negative curvature on any compact K\"ahler manifold $X$ with $c_1(X)<0$. In this case, since $c_1(X)<0$, we can choose a K\"ahler form $\o_0$ in the cohomology class $-c_1(X)$. But the Ricci curvature form $Ricci(\o_0)$ is still in $c_1(X)$, and thus there is a smooth function $f(z)$ with $Ricci(\o_0)+\o_0={i\over 2}\ddb f$. Taking $F(z,\vp)=e^{f+\vp}$ in
the equation (\ref{CY}) and taking again $i\ddb$ of both sides, we see that the solution of (\ref{CY}) satisfies now the K\"ahler-Einstein condition 
\bea
Ricci(\o_\vp)=-\o_{\vp}.
\eea

We note that the K\"ahler-Einstein problem for compact K\"ahler manifolds $X$ with $c_1(X)>0$ is still open
at this time, despite a lot of progress
\cite{TY87, Si, N, T90b, T97, D10, D11a, D11b, CDa, CDb}. A well-known conjecture of Yau \cite{Y93} asserts the equivalence between the existence of such a metric on $X$ and the stability of $X$ in geometric invariant theory. This can be reduced, just as above for the cases $c_1(X)=0$ and $c_1(X)<0$, to a complex Monge-Amp\`ere equation of the form (\ref{CY}), but with $F(z,\vp)=e^{f(z)-\vp}$. Thus the conjecture of Yau asserts the equivalence between the solvability of a complex Monge-Amp\`ere equation and a global, algebraic-geometric, condition. Clearly,  bringing the algebraic-geometric conditions into play in the solution of a non-linear partial differential equation is an important and challenging problem. The two major successes in this direction are the theorem of Donaldson-Uhlenbeck-Yau \cite{D87, UY}, on the equivalence between the existence of a Hermitian-Einstein metric on a holomorphic vector bundle $E\to (X,\o)$ and the Mumford-Takemoto stability of $E$, and the recent results of Donaldson \cite{D08} on the equivalence between the existence of metrics of constant scalar curvature on toric 2-folds and their K-stability. However, there are still many unanswered questions in this direction.

\subsection{The method of continuity}

The original proof of Theorems \ref{Yau} and \ref{YauAubin} is by the method of continuity, and this has remained a prime method for solving complex Monge-Amp\`ere equations to this day. In this method, the equation to be solved is deformed continuously to an equation which we know how to solve. For example, 
one introduces for Theorem \ref{Yau} the deformation, 
\bea
(\o_0+{i\over 2}\ddb \vp)^n={\int_X\o_0^n\over \int_X e^{tf(z)}\o_0^n}e^{tf(z)}\o_0^n,
\qquad 0\leq t\leq 1
\eea
and for Theorem \ref{YauAubin} the deformation, 
\bea
(\o_0+{i\over 2}\ddb \vp)^n=e^{tf(z)+\vp}\o_0^n,
\qquad 0\leq t\leq 1.
\eea
These equations admit trivially the smooth solution $\vp=0$ at $t=0$. It is not difficult to show, by the implicit function theorem, that the set of parameters $t$ for which the equation is solvable is open. So to show that this set is the full interval $[0,1]$ reduces to show that it is closed. This in turn reduces to the proof of a priori estimates for the solutions $\vp$, assuming that they already exist and are smooth.

\section{A Priori Estimates: $C^0$ Estimates}
\setcounter{equation}{0}

We begin by discussing
$C^0$ estimates for the most basic
complex Monge-Amp\`ere equation.
Let $(X,\o_0)$ be a compact K\"ahler manifold without boundary, and consider the equation
\bea
\label{CYF}
(\o_0+{i\over 2}\ddb \vp)^n=F(z)\,\o_0^n
\eea
for a smooth function $\vp$ satisfying the condition $\vp\in \PSH(X,\o_0)$, with $F(z)$ a smooth strictly positive function. Since the equation is invariant under shifts of $\vp$ by constants, we may assume that ${\rm sup}_X\vp=0$. It is well-known that all functions in $PSH(X,\o_0)$ satisfy an exponential integrability condition, and hence their $L^p$ norms are all uniformly bounded by constants depending only on the K\"ahler class $[\o_0]$ and on $p$, for any $1\leq p<\infty$ (see e.g. Appendix A). But the $L^\infty$, or $C^0$ estimate, is fundamentally different. In this section, we discuss several methods for obtaining $C^0$ estimates.

\subsection{Yau's original method}

Yau's original method was by Moser iteration. Set $\psi={\rm sup}_X\vp-\vp+1\geq 1$ and let $\a\geq 0$. 
Since $(F-1)\o_0^n=(\o_0+{i\over 2}\ddb\vp)^n-\o_0^n
={i\over 2}\ddb\vp\sum_{j=0}^{n-1}(\o_0+{i\over 2}\ddb\vp)^{n-1-j}\o_0^j$,
we find, after multiplying by $\psi^{\a+1}$ and integrating by parts.
\bea
\label{Moser1}
\int_X\psi^{\a+1}(F-1)\o_0^n
=
(\a+1)\sum_{j=0}^{n-1}\int_X \psi^\a \,i\p\psi\wedge \bar\p\psi \,(\o_0+{i\over 2}\ddb\vp)^{n-1-j}\o_0^j.
\eea
All the integrals on the right hand side are positive. Keeping only the contribution with $j=n-1$, we obtain
\bea
|\int_X\psi^{\a+1}(F-1)\o_0^n| &\geq&  (\a+1)
\int_X \psi^\a \, i\p\psi \wedge\bar\p\psi\, \o_0^{n-1}
\nonumber\\
&=&
{(\a+1)\over 2({\a\over 2}+1)^2}\int_X i\p(\psi^{{\a\over 2}+1})\wedge \bar\p(\psi^{{\a\over 2}+1})\wedge\o_0^{n-1}.
\eea
and hence, with $C_1$ depending only on $\|F\|_{L^\infty}$, and all norms 
and covariant derivatives with respect to the metric $\o_0$,
\bea
\label{alpha}
\|\nabla (\psi^{{\a\over 2}+1})\|^2
\leq C_1{n({\a\over 2}+1)^2\over\a+1}
\int_X \psi^{\a+1}\o_0^n.
\eea
On the other hand, the Sobolev inequality asserts that
\bea
\label{sobolev}
\|u\|_{L^{2n\over n-1}}^2\leq C_2(\|\nabla u\|_{L^2}^2+\|u\|_{L^2}^2)
\eea
with $C_2$ the Sobolev constant of $(X,\o_0)$. Applied to $u= \psi^{p\over 2}$, it can be expressed as
\bea
\|\psi\|_{L^{p\beta}}^p
\leq C_2(\|\nabla(\psi^{p\over 2})\|_{L^2}^2+\|\psi\|_{L^p}^p),
\eea
with $\beta={n\over n-1}>1$. 
Setting $p=\a+2$, and applying the inequality (\ref{alpha}), we find
\bea
\|\psi\|_{L^{p\beta}}
\leq (C_3 p)^{1\over p}\|\psi\|_{L^p},
\qquad p\geq 2,
\eea
with a constant $C_3$ depending only on $n$, $\|F\|_{L^\infty}$, and the Sobolev constant of $(X,\o_0)$. We can iterate $p\to p\beta\to\cdots\to p\beta^k$ and get
\bea
\log \|\psi\|_{L^\infty}
\leq
\sum_{k=0}^\infty {\log (C_3p\beta^k)\over p\beta^k}+\log \|\psi\|_{L^p}
=C_{4,p}+\log\|\psi\|_{L^p}.
\eea
An a priori bound for $\|\psi\|_{L^p}$ for any fixed finite $p$
can be obtained from the exponential estimate for plurisubharmonic functions in  Appendix A. Alternatively, if we apply Moser iteration instead to the function $\varphi$ normalized to have average $0$,
we can obtain an a priori bound for $\|\varphi\|_{L^2}$
from the analogue of (\ref{alpha}) by taking $\alpha=0$, and applying the Poincar\'e inequality to the left hand side. Either way gives

\begin{theorem}
\label{C0-1}
Let $\vp$ be a smooth solution
of the equation (\ref{CYF}) on a compact K\"ahler manifold $(X,\o_0)$ without boundary,
$F>0$, and $\vp\in \PSH(X,\o_0)$. Then $\|\psi\|_{L^\infty(X)}$ is bounded by a constant depending only on $n$,
an upper bound for $\|F\|_{L^\infty(X)}$, and the K\"ahler form $\o_0$. The dependence on the K\"ahler form $\o_0$ can be stated more precisely as a dependence on the Sobolev constant and the Poincar\'e constant of $\o_0$, or on the exponential bound for $\o_0$.
\end{theorem}

The Moser iteration method is now widely used in the study of Monge-Amp\`ere and other non-linear equations. An important variant has been introduced by Weinkove \cite{W}, where the Moser iteration is applied to $e^\vp$ instead of $\vp$.
Applications of this variant are in \cite{SW, ST06, TWY}. 

\subsection{Reduction to Alexandrov-Bakelman-Pucci estimates}

It was suggested early on by Cheng and Yau that the Alexandrov-Bakelman-Pucci estimate can be applied to the complex Monge-Amp\`ere equation. They did not publish their work, but a detailed account was subsequently provided by Bedford \cite{B} and Cegrell and Persson \cite{CP}. Using the Alexandrov-Bakelman-Pucci estimate, Blocki \cite{B11a}
gives the following proof of the $C^0$ estimate. This proof is of particular interest as it is almost a local argument. We follow closely Blocki's presentation.

\medskip

Let $D$ be any bounded domain in ${\bf C}^n$, $u\in C^2(\bar D)$,
$u_{\bar kj}\geq 0$, and $u=0$ on $\p D$. Then
\bea
\label{ABP}
\|u\|_{C^0}
\leq C\,\|{\rm det}\,u_{\bar kj}\|_{L^2}^{1\over n}
\eea
where $C=C(n,{\rm diam}\,D)$ depends only on the diameter of $D$ and the dimension $n$.
To see this, we apply the ABP estimate (from \cite{GT}, Lemma 9.2) to get
\bea
\|u\|_{C^0}\leq c_n({\rm diam}\,D)\,\big(\int_\Gamma {\rm det}\,D^2u\big)^{1\over 2n}
\eea
where $\Gamma$ is the contact set, defined by
\bea
\Gamma=\{z\in D; u(w)\geq u(z)+\<Du(z),w-z\>,\ {\rm for\ all} \ w\in D\}.
\eea
On the contact set $\Gamma$, the function $u$
satisfies $D^2u\geq 0$, and for such functions,
we have the following inequality between the determinants of the real and complex Hessians,
\bea
{\rm det}\,u_{\bar kj}
\geq 2^{-n}({\rm det}\,D^2u)^{1\over 2}.
\eea
This proves the estimate (\ref{ABP}).

\smallskip

Let now $z\in D$, $h>0$, and define the sublevel set $S(z,h)$ by
\bea
S(z,h)=\{w\in D; u(w)<u(z)+h\}.
\eea
If $S(z,h)\subset\subset D$, then applying the previous inequality to $S(z,h)$ instead of $D$ gives
\bea
\|u-u(z)-h\|_{C^0(S(z,h))}
&\leq&
C(n,{\rm diam}\,D)
\|{\rm det}\,u_{\bar kj}\|_{L^2(S(z,h))}^{1\over n}
\nonumber\\
&\leq&
C(n,{\rm diam}\,D)
|S(z,h)|^{1\over 2nq}\|{\rm det}\,u_{\bar kj}\|_{L^{2p}}^{1\over n}
\eea
for any $p>1$, ${1\over p}+{1\over q}=1$. In particular, we obtain the following lower bound for $|S(z,h)|$
\bea
h\leq C(n,{\rm diam}\,D)
|S(z,h)|^{1\over 2qn}\|{\rm det}\,u_{\bar kj}\|_{L^{2p}}^{1\over n}.
\eea
On the other hand, we have the following easy upper bound for $|S(z,h)|$,
\bea
|S(z,h)|(-u(z)-h)\leq\int_{S(z,h)}(-u)
\leq \|u\|_{L^1(D)}
\eea
If we choose $z$ to be the minimum point for $u$,
and eliminate $|S(z,h)|$ between the two inequalities, we obtain
a lower bound for $u$ in terms of $h$ and $h^{-1}$.

\medskip

This can be applied to the $C^0$ estimate for the Monge-Amp\`ere equation
$(\o_0+{i\over 2}\ddb\vp)^n=F(z)\,\o_0^n$ on a compact
K\"ahler manifold $(X,\o_0)$. Let $z$ be the minimum point for $\vp$ on $X$,
and let $K(w,\bar w)$ be a K\"ahler potential for $\o_0$ in a neighborhood of $z$.
By adding a negative constant and
shifting $K(w,\bar w)$ by the real part of a second order polynomial in $w$ if necessary,
we can assume that $K(w,\bar w)\leq 0$ in a ball $B(z,2r)$
around $z$, $K(w,\bar w)\geq K(0)+h$ for $r\leq |w|\leq 2r$, and $K(w,\bar w)$
attains its minimum in $B(z,2r)$ at $0$. The constant $h>0$
depends only on the K\"ahler form $\o_0$.
Then the function $u=K+\vp$
attains its minimum in $B(z,2r)$ at $z$, and the corresponding set $S(z,h)\subset
B(z,2r)\setminus B(z,r)$ has compact closure. By the preceding inequalities,
we obtain a lower bound for $u(z)$, depending only on $\o_0$ and the $L^{2p}$
norm of $F$, for any $p>1$. Thus

\begin{theorem}
\label{C0-2}
Let the setting be the same as in Theorem \ref{C0-1}. Then for any $p>1$,
$\|\vp\|_{L^\infty(X)}$ can be bounded by a constant depending only on $n$, an upper bound
for $\|F\|_{L^{2p}(X)}$, and the K\"ahler form $\o_0$.
\end{theorem}

\subsection{Methods of pluripotential theory}

A third method for $C^0$ estimates
was introduced by Kolodziej \cite{K98}. This method 
combines the
classic approach of De Giorgi with modern techniques of pluripotential theory. It produces $C^0$ bounds even when the right hand side $F$ is only in $L^p(X)$ for some $p>1$. As shown by Eyssidieux, Guedj, and Zeriahi \cite{EGZ, EGZ09} and 
Demailly and Pali \cite{DP}, it can also be extended to a family setting, where the background K\"ahler form $\o_t$ is allowed to degenerate to a closed form $\chi$ which is just non-negative. Other family versions of Kolodziej's $C^0$ estimates are in \cite{KT, DZ, TZ}. As we shall see later, such family versions are important for the study of 
singular K\"ahler-Einstein metrics and Monge-Amp\`ere equations on complex manifolds with singularities.

\medskip

Let $(X,\o_0)$ be a compact K\"ahler manifold. Let $\chi\geq 0$ be a
$C^\infty$ closed semi-positive $(1,1)$-form which is not identically $0$.
Set
\bea
\o_t=\chi+(1-t)\o_0,
\qquad t\in (0,1)
\eea
and let $[\o_t^n]=\I_X\o_t^n$.
Consider the equation
\bea
(\o_t+{i\over 2}\ddb \vp_t)^n=F_t\o_t^n
\eea
for some strictly positive function $F_t$ and
$\vp_t\in \PSH(X,\o_t)\cap L^\infty(X)$. Then we have the
following family version of the $C^0$ estimates of Kolodziej
\cite{K98}, due to Eyssidieux-Guedj-Zeriahi \cite{EGZ, EGZ09}
and Demailly-Pali \cite{DP}:

\begin{theorem}
\label{C0} Let $A>0$ and
and $p>1$. Assume that $\chi\leq A\o_0$ and ${1\over [\o_t^n]}{\o_t^n\over \o_0^n}\leq A$.
Assume also that the functions $F_t$ are in $L^p(X,\o_t^n)$  and that
\bea
{1\over [\o_t^n]}\int_X F_t^p\o_t^n \leq A^p<\infty
\eea
for all $t\in (0,1)$. Normalize $\vp_t$ so that
${\rm sup}_X\vp_t=0$. Then there exists a constant $C>0$, depending only on $n$,
$\o_0$ and $A$, so that
\bea
{\rm sup}_{t\in [0,1)}\|\vp_t\|_{L^\infty(X)}\leq C.
\eea
\end{theorem}

\medskip

\noindent
{\it Proof}. Recall the notion of capacity of a Borel set $E$ with respect to a K\"ahler form $\o$,
\bea
{\rm Cap}_\o(E)
=
{\rm sup}\{\int_E(\o+{i\over 2}\ddb u)^n;
\ u\in \PSH(X,\o),\ 0\leq u\leq 1 \}.
\eea
Set
\bea
\label{estimate1}
f_t(s)=({{\rm Cap}_{\o_t}(\vp_t<-s)\over [\o_t^n]})^{1\over n}.
\eea
It suffices to show that there exists $s_\infty<\infty$ independent of $t$
so that 
\bea
f_t(s)=0 \qquad {\rm for}\ s>s_\infty.
\eea
Since $f_t(s)^n\geq {1\over [\o_t^n]}\int_{\vp<-s}\o_t^n$, it would follow that $\vp_t\geq -s_\infty$
a.e. with respect to the measure $\o_t^n$, and since $\vp_t$ is upper semi-continuous, that $\vp_t\geq -s_\infty$
everywhere. The following classic lemma of De Giorgi provides sufficient conditions
for the existence of $s_\infty$:

\begin{lemma}
\label{degiorgilemma}
Let $f:{\bf R}^+\to {\bf R}^+$ satisfy the following conditions:

{\rm (a)} $f$ is right-continuous;

{\rm (b)} $f$ decreases to $0$;

{\rm (c)} There exist positive constants $\alpha, A_\alpha$
so that for all $s\geq 0$ and all $0\leq r\leq 1$, we have
\bea
\label{degiorgi}
r\,f(s+r)\leq A_\alpha f(s)^{1+\alpha}.
\eea
Then there exists $s_\infty$, depending only on $\a,A_\alpha$ and the smallest value $s_0$ for which we have
$f(s_0)^\al\leq (2A_\alpha)^{-1}$ so that $f(s)=0$ for $s>s_\infty$.
In fact, we can take $s_\infty=s_0+2A_\alpha(1-2^{-\a})^{-1}f(s_0)^\alpha$.
\end{lemma}

We shall show that the above functions $f_t(s)$ satisfy the conditions
of Lemma \ref{degiorgilemma}. The right-continuity (a) of the function $f_t(s)$ is a consequence
of the fact that, for any K\"ahler form $\o$, and any sequence of increasing sequence of Borel sets $E_j\subset E_{j+1}$, we have ${\rm Cap}_\o(\cup_{j=1}^\infty E_j)
={\rm lim}_{j\to\infty}{\rm Cap}_\o (E_j)$.  
Clearly $f_t(s)$ decreases as $s$ increases. In fact, it does so uniformly to $0$ in $t$
as is shown by the following lemma:

\begin{lemma}
\label{uniformdecrease}
There exists a constant $C$ depending only on $\o_0$ and an upper bound $A$ for $\chi$ so that
\bea
f_t(s)^n\leq C\,s^{-1}.
\eea
\end{lemma}

\noindent
{\it Proof of Lemma \ref{uniformdecrease}}. Let $u\in \PSH(X,\o_t)$. Then
\bea
\int_{\vp_t<-s}(\o_t+{i\over 2}\ddb u)^n
&\leq& {1\over s}\int_X(-\vp_t)(\o_t+{i\over 2}\ddb u)^n
\\
&=&
{1\over s}\int_X(-\vp_t)\o_t^n
+
{1\over s}\int_X(-\vp_t){i\over 2}\ddb u\sum_{j=0}^{n-1}\o_t^j(\o_t+{i\over 2}\ddb u)^{n-1-j}\nonumber
\eea
Writing $\o_t^n\leq A[\o_t^n]\o_0^n$,
and noting that $\PSH(X,\o_t)\subset \PSH(X,(A+1)\o_0)$,
we can bound the first integral on the right hand side by $C[\o_t^n]$,
in view of Theorem \ref{Hormander}
on exponential estimates for plurisubharmonic functions.
The other integrals can be re-expressed as
\bea
&&
\int_X \vp_t {i\over 2}\ddb u\, \o_t^j(\o_t+{i\over 2}\ddb u)^{n-1-j}
=-\int_X {i\over 2}\ddb\vp_t\, u\, \o_t^j(\o_t+{i\over 2}\ddb u)^{n-1-j}
\nonumber\\
\qquad
&&
=
-\int_X u\,(\o_t+{i\over 2}\ddb\vp_t) \o_t^j(\o_t+{i\over 2}\ddb u)^{n-1-j}
+
\int_X u\,\o_t^{j+1}(\o_t+{i\over 2}\ddb u)^{n-1-j}.
\nonumber
\eea
For $0\leq u\leq 1$, we can write
\bea
|\int_X u(\o_t+{i\over 2}\ddb\vp_t) \o_t^j(\o_t+{i\over 2}\ddb u)^{n-1-j}|
\leq \int_X\,(\o_t+{i\over 2}\ddb\vp_t) \o_t^j(\o_t+{i\over 2}\ddb u)^{n-1-j}
=[\o_t^n]\nonumber
\eea
and similarly for the other integral. Thus we obtain an upper bound $C[\o_t^n]$,
and taking the supremum in $u$ establishes the desired inequality. Q.E.D.

\medskip

It remains to establish the property
(c) for $f_t(s)$. For this, we need the following two propositions:

\begin{lemma}
\label{capacityvolume}
Let $\vp\in PSH(X,\o)\cap L^\infty(X)$. Then for all $s>0$, $0\leq r\leq 1$,
\bea
r^n
{\rm Cap}_\o(\vp<-s-r)
\leq
\int_{\vp<-s}(\o+{i\over 2}\ddb\vp)^n.
\eea
\end{lemma}

\begin{lemma}
\label{volumecapacity}
There exist constants $\delta, C>0$ so that for any open set $E\subset X$,
and any $t\in [0,1)$, we have
\bea
{1\over[\o_t^n]}\int_E\o_t^n
\leq C \,{\rm exp} \left[-\delta\left(\frac{[\o_t^n]}{{\rm Cap}_{\o_t}(E)}\right)^{1/n}
\right].
\eea
\end{lemma}

Assuming these two lemmas for the moment, we can readily establish the
inequality (c) in Lemma \ref{degiorgilemma}. For $\al>0$ we have
\bea
[r\,f_t(s+r)]^n
&=&r^n{{\rm Cap}_\o(\vp_t<-s-r)\over [\o_t^n]}
\leq {1\over [\o_t^n]}\int_{\vp_t<-s}(\o+{i\over 2}\ddb\vp_t)^n
\nonumber\\
&=&
{1\over [\o_t^n]}\int_{\vp_t<-s}F_t\o_t^n
\leq
({1\over [\o_t^n]}\int_{\vp_t<-s}F_t^p\o_t^n)^{1\over p}
({1\over [\o_t^n]}\int_{\vp_t<-s}\o_t^n)^{1\over q}
\nonumber\\
&\leq &
A\,{\rm exp}\left(-{\delta\over q} 
\left[{[\o_t^n]\over {\rm Cap}_\o(\vp_t<-s)}
\right]^{1\over n}
\right)
\leq A_\al\,f^{(1+\a)n}.
\eea

It remains to prove the two lemmas. Let
$u\in PSH(X,\o)$ with $0\leq u\leq 1$, and write
\bea
r^n\int_{\vp<-s-r}(\o+{i\over 2}\ddb u)^n&=&
\int_{\vp<-s-r}(r\o+{i\over 2}\ddb \,ru)^n
\nonumber\\
&\leq&
\int_{\vp<-s-r+ru}(\o+{i\over 2}\ddb(ru-s-r) )^n
\nonumber\\
&\leq&
\int_{\vp<-s-r+ru}(\o+{i\over 2}\ddb\vp)^n
\eea
where we have applied the comparison principle.
Since $-r+ru$ is negative, this last integral is bounded by the integral
over the larger region $\{\vp<-s\}$, and Lemma \ref{capacityvolume} is proved.
The next lemma requires some properties of global extremal functions
\cite{GZ, Ze}:

\begin{lemma}
\label{extremalfunction}
Let $E\subset X$ be an open set, and define its global extremal function $\it\psi_{E,\o}$ as the upper semi-continuous envelope of the following function
$\ti\psi_{E,\o}$,
\bea
\label{extremal}
\ti\psi_{E,\o}
=
{\rm sup}\{u\in \PSH(X,\o); u=0\ {\rm on}\ E\}.
\eea
Then 

{\rm (a)} $\psi_{E,\o}\in \PSH(X,\o)\cap L^\infty(X)$

{\rm (b)} $\psi_{E,\o}=0$ on $E$

{\rm (c)} $(\o+\frac{i}{2}\ddb\psi_{E,\o})^n=0$ on $X\setminus \bar E$.

\end{lemma}

We can now prove Lemma \ref{volumecapacity}. Let $E'\subset E$ be any relatively compact open subset.
Then
\bea
{1\over [\o_t^n]}
\int_{E'}\o_t^n
&=&
{1\over[\o_t^n]}e^{-\delta{\rm sup}_X \psi_{E',\o_t}}
\int_{E'}e^{-\delta(\psi_{E'\o_t}-{\rm sup}_X\psi_{E',\o_t})}\o_t^n
\nonumber\\
&\leq& e^{-\delta \,{\rm sup}_X\psi_{E',\o_t}}
\,
A\int_Xe^{-\delta(\psi_{E',\o_t}-{\rm sup}_X\psi_{E',\o_t})}\o_0^n
\eea
where $A$ is an upper bound for ${1\over [\o_t^n]}{\o_t^n\over\o_0^n}$.
Since $\chi\leq A\o_0$ also by assumption,
$PSH(X,\o_t)\subset \PSH(X,(A+1)\o_0)$,
and Theorem \ref{Hormander}
implies that the integral on the right
hand side is bounded by a constant independent of $t$ and $E'$. 
We can now complete the proof of Lemma \ref{volumecapacity}. First, observe that
if ${\rm sup}_X\psi_{E',\o_t}\leq 1$, then
\bea
[\o_t^n]=
\int_{\bar E'}(\o_t+\frac{i}{2}\ddb\psi_{E',\o_t})
\leq {\rm Cap}_{\o_t}(\bar E')
\leq {\rm Cap}_{\o_t}(E)
\leq{\rm Cap}_{\o_t}(X)=[\o_t^n].
\eea
Thus $\frac{[\o_t^n]}{{\rm Cap}_{\o_t}(E)}=1$, and a constant $C_\delta$ can clearly be chosen so that the desired inequality holds. Next, assume that
${\rm sup}_X\psi_{E',\o_t}^n> 1$. We can write 
\bea
({\rm sup}_X\psi_{E',\o_t})^{-n}
&=&
({\rm sup}_X\psi_{E',\o_t})^{-n}\frac{\int_X(\o_t+\frac{i}{2}\ddb\psi_{E',\o_t})^n}
{[\o_t^n]}
=
({\rm sup}_X\psi_{E',\o_t})^{-n}\frac{\int_{\bar E'}(\o_t+\frac{i}{2}\ddb\psi_{E',\o_t})^n}
{[\o_t^n]}
\nonumber\\
&\leq&
\frac{\int_{\bar E'}(\o_t+\frac{i}{2}\ddb (\frac{\psi_{E',\o_t}}{{\rm sup}_X\psi_{E',\o_t}}))^n}
{[\o_t^n]}.
\eea
This last term is bounded by $[\o_t^n]^{-1}{\rm Cap}_{\o_t}(\bar E')
\leq [\o_t^n]^{-1}{\rm Cap}_{\o_t}(E)$. Thus we obtain
\bea
{1\over [\o_t^n]}\int_{E'}\o_t^n
\leq
{\rm exp}(-\delta (\frac{[\o_t^n]}{{\rm Cap}_{\o_t}(E)})^{1\over n}).
\eea
Taking limits as $E'$ increases to $E$ establishes Lemma \ref{volumecapacity}. The proof of the theorem is complete.

\bigskip
We observe that a more straightforward adaption of Kolodziej's original argument can be applied to the special case of algebraic manifolds with the
background class being big and semi-ample [ZZh].

\medskip
Finally, we note that all three proofs of $C^0$ estimates can be extended to the equation (\ref{Calabiconjecture}) on Hermitian manifolds. For the Moser iteration method, this is carried out in \cite{Ch, TW1, TW2}. It is interesting that, in the K\"ahler case,
only one term in the right hand side of
(\ref{Moser1}) was needed, while the other terms 
are also needed in the Hermitian case. The extension of the pluripotential
method to the Hermitian case is in \cite{DK}, while the extension of the Alexandrov-Bakelman-Pucci method is in \cite{B11a}.

\section{Stability Estimates}
\setcounter{equation}{0}

In this section, we shall establish the stability and uniqueness of the continuous solutions of the complex Monge-Amp\`ere equations due to Kolodiej \cite{K03}. We shall closely follow the arguments in \cite{K05}.  Let $(X, \omega)$ be an $n$-dimensional  compact K\"ahler manifold and let 
\bea
\F_{p, A}= \{ F\in L^p(X)~|~ F \geq 0, ~\int_X F^p \omega^n \leq A, ~\int_X F \omega^n = \int_X \omega^n\}
\eea
for $p>1$ and $A>0$.  

\begin{theorem} 
\label{cmstability} For any two $F, G\in \F_{p, A}$,  let $\varphi$ and $\psi \in \PSH(X, \omega) \cap C(X)$ be solutions of the following Monge-Amp\`ere equations
\bea
(\omega+ {i\over 2}\ddbar \varphi)^n = F(z) \omega^n , ~~~~~(\omega+ {i\over 2}\ddbar \psi)^n = G(z) \omega^n
\eea
normalized by 
\bea
\sup_X  \varphi  = \sup_X  \psi =0.
\eea
Then for any $\e>0$, there exists $C>0$ depending only on $\e, p, A$ and $(X, \omega)$, so that
\begin{equation}
\|\varphi-\psi\|_{L^\infty(X)} \leq C \, \|F-G\|_{L^1(X)}^{\frac{1}{n+3+\e}}.
\end{equation}

\end{theorem}

Theorem \ref{cmstability} has been generalized in \cite{DZ} to nonnegative, big and smooth closed $(1,1)$-form $\omega$ if $\omega$ is chosen appropriately. The uniqueness of the solutions 
$\vp\in \PSH(X, \omega)\cap C^0(X)$, ${\rm sup}_X\vp=0$, to the Monge-Amp\`ere equation
(\ref{CYF}) for $F\in\F_{p,A}$
for some $p>1$ and $A>0$ follows immediately from the stability theorem. 

\medskip

As before, we write $\omega_\varphi = \omega+ {i\over 2}\ddbar \varphi$ for any $\varphi \in \PSH(X, \omega)$.  The following lemma is a generalization of Lemma \ref{capacityvolume} and a detailed proof can be found in \cite{K05}.

\begin{lemma} \label{caplessvol} Let $\varphi$ and $\psi\in \PSH(X, \omega)\cap C(X) $  with $0\leq \varphi \leq C$. Then for any $s>0$, 
\begin{equation}
{\rm Cap}_\o ( \{ \psi + 2s < \varphi \} ) \leq \left( \frac{C+1}{ s}\right) ^{n} \int_{\psi + s < \varphi}  \omega_\psi ^n.      
\end{equation} 
\end{lemma}

The following lemma is well-known for smooth plurisubharmonic functions. The general result can also be found in \cite{K05}.

\begin{lemma} \label{kcomparison} Let $\B$ be an Euclidean ball in $\C^n$ with the standard Euclidean volume form $\Omega$. For any $u\in PSH(\B)\cap C(\overline{\B})$ and any nonnegative $F\in L^1(\overline{\B})$ with 
\bea
(\frac{i}{2}\ddbar u)^n \geq F\, \Omega, ~~~~(\frac{i}{2}\ddbar v)^n \geq F\,\Omega, 
\eea
we have 
\bea
(\frac{i}{2}\ddbar u )^k \wedge (\frac{i}{2}\ddbar v)^{n-k} \geq F\,\Omega
\eea 
for all $k=0, ..., n$.

\end{lemma}

By Theorem \ref{C0}, there exists $a>0$ depending on $3^pA$, $p$ and $(X, \omega)$ such that for any $F\in \F_{p, 3^p A}$, the solution $u\in \PSH(X, \omega)\cap L^\infty(X)$ of $(\omega + {i\over 2}\ddbar u)^n= F\omega^n$ satisfies 
\bea
\sup_X u - \inf_X u \leq a.
\eea
Without loss of generality, we can assume that $\int_X F \o^n = \int_X G \o^n=\int_X \omega^n = 1$ and 
\bea
\int_{\psi < \varphi} (F+G)\o^n \leq 1.
\eea
Furthermore, we let $0<t_0< (q-1)/2$, where $q= (3/2)^{1/n}$ and we further assume $t_0^{n+3+\e} < 1/3$. For any $0<t< t_0$, we define
for all nonnegative integers $k$, 
\bea
E_k = \{ \psi < \varphi - k a t\}
\eea 
for a fixed $\e>0$. 
From now on we assume that for $0<t<t_0$, 
\bea
\|F-G\|_{L^1(X)} = t^{n+3+ \e}.
\eea 
We have immediately 
\bea
\int_{E_0}  G\omega^n = \frac{1}{2} \int_{E_0} \left( (F+G) + (G-F)\right) \omega^n \leq \frac{1}{2}(1+\frac{1}{3}) = \frac{2}{3}.
\eea
Now we define a new function $H$ such that $H= 3G/2$ on $E_0$ and $H= c_0$ on the complement of $E_0$ so that $\int_X H \omega^n =1$. Obviously, $c_0>0$ and $H\in \F_{p, (3/2)^pA}$. Then there exists a unique $\rho \in \PSH(X, \omega) \cap C(X)$ such that 
\bea
\omega_\rho^n = H\omega^n, ~~~~ \sup_X \rho = 0.
\eea 
Furthermore, 
\bea
- a \leq \rho\leq 0 .
\eea
We now define
\bea
E = \{ \psi < (1-t) \varphi + t\rho - at\}, ~~~~ S = \{ F< (1-t^2) G\}.
\eea
The following lemma can be easily verified.
\begin{lemma}
\bea
E_2 \subset E \subset E_0.
\eea
\end{lemma}

\begin{lemma}On $E_0\setminus S$, for $k=0, ..., n$, we have
\bea
\omega_\varphi^k \wedge \omega_\rho^{n-k} \geq q^{n-k} (1-t^2)^{k/n} G \omega^n.
\eea
\end{lemma}

\noindent{\it Proof.}  On $E_0\setminus S$, we have
\bea
\left( (1-t^2)^{-1/n} \omega_\varphi\right) ^n \geq G \omega^n, ~~~~ (q^{-1} \omega_\rho)^n = G\omega^n.
\eea
The lemma  then follows from Lemma \ref{kcomparison}.  Q.E.D.

\begin{lemma}
Let $B = \int_{E_2} G \omega^n.$  Then
\begin{equation} 
B \leq  \frac{3}{q-1} t^{n +\e} .
\end{equation}
\end{lemma}

\noindent {\it Proof. } On $E_0\setminus S$, we have 
\begin{eqnarray}
\label{qcomp}
\omega^n_{(1-t)\varphi+ t\rho }& =&\sum_{k=0}^n   \left(\begin{array}{c}
    n \\ 
    k \\ 
  \end{array}\right) (1-t)^k t^{n-k} \omega_\varphi ^k \wedge \omega_\rho^{n-k} \nonumber\\
&\geq&\sum_{k=0}^n   \left(\begin{array}{c}
    n \\ 
    k \\ 
  \end{array}\right) q^{n-k}(1-t^2)^{k/n}(1-t)^k t^{n-k} G \omega^n .
\end{eqnarray}
The right hand side can in turn be estimated by
\bea
\label{gcomp2}
&&
\sum_{n=0}^\n \pmatrix{n\cr k\cr} q^{n-k}(1-t^2){k\over n}
(1-t)^kt^{n-k}G\o^n
=(  qt + (1-t) (1-t^2)^{1/n} )^nG \\
&&
\quad
\geq ( (1-t)(1-t^2) + qt)^n G\omega^n
\geq (1+ (q-1)t - t^2)^n G \omega^n \geq (1+(q-1)t/2) G\omega^n,
\nonumber
\eea
where we make use of the additional assumption that $t< t_0 < (q-1)/2.$
On the other hand, since  $\int_S F\omega^n \leq (1-t^2) \int_S G \omega^n$ by the definition of $S$, we have  
\bea
t^2 \int_S G\omega^n \leq \int_S (G-F)\omega^n \leq t^{n+3+\e}
\eea 
and so
\begin{equation} \label{qineq} 
\int_S G\omega^n \leq t^{n+1+\e}.
\end{equation}

\smallskip

The above inequality implies that 
\bea\label{eminuss}
\int_{E\setminus S} G \omega^n \leq   \frac{2}{q-1}t^{n+\e}.
\eea
Thus by (\ref{gcomp2}) and (\ref{qineq}), 
\bea
B\leq \int_{E}G\omega^n \leq \int_{E\setminus S} G\omega^n + \int_S G\omega^n\leq   \frac{3}{q-1} t^{n +\e} .
\eea
The lemma is then proved.  Q.E.D.

\begin{lemma}

\begin{equation}
{\rm Cap}_\omega (E_4) \leq \frac{3}{q-1} (2a) ^{-n} (a+1)^n t^\epsilon.
\end{equation}

\end{lemma}

\noindent {\it Proof.}   By Lemma \ref{caplessvol}, we have 
\begin{eqnarray*}
{\rm Cap}_\omega(E_4)  \leq  \frac{(a+1)^n}{(2at)^n} \int_{E_2} G \omega^n 
=  \frac{(a+1)^n}{(2at)^n} B 
\leq\frac{ (a+1)^n}{(2a)^{n}} \frac{3}{q-1}  t^\e.
\end{eqnarray*}

The following lemma is used in Kolodziej's original proof of the $L^\infty$ estimates. We refer the readers to the detailed proof in \cite{K05}. 

\begin{lemma}\label{prec0} Let $\varphi, \psi \in \PSH(X, \omega)\cap C(X)$ with $0\leq \varphi \leq C$. Let 
\bea
U(s) = \{ \psi - s< \varphi \}, ~~~~ \alpha(s) = {\rm Cap}_\omega (U(s)).
\eea
Assume that 

{\rm (1)} $\{ \psi - S < \varphi\}\neq \emptyset$ for some $S$,

{\rm (2)} For any Borel set $K$, 
\bea
\int_K \omega_\psi^n \leq f ({\rm Cap}_\omega(K)),
\eea
where 
$f(x) = \frac{x}{h(x^{-1/n})}$
and $h(x): \R^+ \rightarrow (0, \infty)$ is a continuous strictly increasing  function satisfying 
$\int_1^\infty \frac{1}{ t h^{1/n}(t) } dt < \infty$.

Then for any $D<1$, we have
\begin{equation}
D\leq \kappa (\alpha(S+D)),
\end{equation}
where 
\bea
\kappa(s) = c(n)(1+C) \left( \int_{s^{-1/n}}^\infty \frac{dx}{x h^{1/n}(x) }+ \frac{1}{h^{1/n}(s^{-1/n}) }\right) 
\eea
for some constant $c(n)$ depending only on $n$.

\end{lemma}

\noindent {\it Proof of Theorem \ref{cmstability}.}      By Lemma \ref{volumecapacity},  for any $\delta>0$ and open set $K$, there exists $C_\delta>0$, 
\bea
\int_K \omega^n \leq C_1 e^{ - (C_2 {\rm Cap}_\omega (K))^{-1} } \leq  C_\delta \left( {\rm Cap}_\omega (K) \right)^{1/\delta}.
\eea
Then $0\leq \varphi+a\leq a$.  We can easily check that we can choose $h(x) = x^{1/\delta}$  and there exists $C'_\delta>0$ such that 
\bea \kappa(s) = C'_\delta s^{1/(\delta n^2)}.
\eea
Now we can prove the theorem by contradiction. Suppose that 
\bea
\{ \psi < \varphi - (4a+1) t\} = \{\psi+a < \varphi+a - (4a+1) t\} \neq \emptyset . 
\eea
Then by applying Lemma \ref{prec0} with $\psi+a$, $\varphi+a$, $S= - (4a+1)t$ and $D= t$,
\bea
t \leq \kappa ( {\rm Cap}_\omega (E_4)) \leq \kappa\left(  \frac{3}{q-1} (2a)^{-n} (a+1)^n  t^\epsilon \right)= C'_\delta \left(  \frac{3}{q-1} (2a)^{-n} (a+1)^n  t^\epsilon \right) ^{1/(\delta n^2)}  .  \nonumber
\eea 
This is a contradiction if we  choose $\delta>0$ sufficiently small and then $t>0$ sufficiently small. Therefore 
$\{ \psi < \varphi - (4a+1) t\} = \emptyset$ and so  
\bea
\sup_X (\varphi - \psi) = \sup_X (\psi  -\varphi ) \leq (4a+1) t =  (4a+1) \left( \|F-G\|_{L^1(X)}\right)^{1/(n+3+\e)} 
\eea 
if we choose $t_0$ sufficiently small.  The theorem is proved. Q.E.D.

\section{A Priori Estimates: $C^1$ Estimates}
\setcounter{equation}{0}

In Yau's original solution of the Calabi conjecture \cite{Y78}, the $C^2$ estimates were shown to follow directly from the $C^0$ estimates. The $C^1$ estimates follow from the $C^0$ and $C^2$ estimates by general linear elliptic theory. However, for more general Monge-Amp\`ere equations where the right hand side may be an expression $F(z,\vp)$ depending on the unknown $\vp$ as well as for the Dirichlet problem, the $C^1$ estimates cannot be bypassed. In this section, we describe the sharpest $C^1$ estimates available at this time. They are due to \cite{PS09a, PS09c}, and they exploit a key differential inequality discovered by Blocki \cite{B09a}.

\medskip
Let $(X,\omega_0)$ be a compact K\"ahler manifold with smooth boundary $\partial X$ (which may be empty)
and complex dimension $n$.
We consider the Monge-Amp\`ere equation on $\bar X$
\bea
\label{MA}
(\o_0+{i\over 2}\ddb \vp)^n=F(z,\vp)\,\o_0^n.
\eea
Here $F(z,\vp)$ is a  $C^2$ function on $\bar X\times\R$
which is assumed to be strictly positive on the set $\bar X\times [\inf\,\vp,\infty)$.
The gradient estimates allow $\vp$ to
be singular along a subset $Z\subset X$, possibly empty, which does not intersect $\p X$. All covariant derivatives and curvatures listed below
are with respect to the metric $\o_0$. Then \cite{PS09c}

\begin{theorem}
\label{C1}
Let $(X,\o_0)$ be a compact K\"ahler manifold, with smooth boundary $\p X$ (possibly empty).
Assume that $\vp\in C^4(\bar X\setminus Z)$ is a solution
of the equation (\ref{MA}) on $\bar X\setminus Z$.
If $Z$ is not empty, 
assume further that $Z$ does not intersect $\p X$, and
that there exists a constant $B>0$ so that
\bea
\label{divisor}
&&
\vp(z)\to +\infty\ {\rm as}\ z\to\ Z,\nonumber\\
&&
\log |\nabla \vp(z)|^2-B\,\vp(z)
\to -\infty \ {\rm as}\ z\to\ Z.
\eea
Then we have the a priori estimate
\bea
\label{C1a}
|\nabla \vp(z)|^2 \leq C_1\,{\rm exp}(A_1\, \vp(z)),
\qquad z\in {\bar X}\setminus Z,
\eea
where $C_1$ and $A_1$ are constants that depend only on upper bounds for
${\rm inf}_X \vp$,
${\rm sup}_{X\times [{\rm inf}\,\vp,\infty)} F$,
${\rm sup}_{X\times [{\rm inf}\,\vp,\infty)}\left(|\nabla F^{1\over n}|+|\pl_\vp F^{1\over n}|\right)$,
${\rm sup}_{\pl X}|\vp|$,
${\rm sup}_{\pl X}|\nabla \vp|$, and the following constant,
\bea
\label{Lambda}
\Lambda=-{\rm inf}_X{\rm inf}_{M>0}{M^j{}_kR^k{}_j{}^p{}_q(M^{-1})^q{}_p
\over
{\rm Tr}\, M\,{\rm Tr}\,M^{-1}},
\eea
where $M= (M^q{}_p)$ runs over all self-adjoint and positive definite endomorphisms.
\end{theorem}

\medskip

When there is no boundary, and the function $F(z,\vp)$ is a function
$F(z)$ of $z$ alone, the equation
(\ref{MA}) is unchanged under shifts of $\vp$ by an additive constant.
Thus the infimum of $\vp(z)$ can be normalized to be $0$ by replacing $\vp(z)\to\vp(z)-{\rm inf}_X\vp$, so
we obtain the estimate
\bea
\label{C1elementary}
|\nabla\vp(z)|^2\leq C_1\,{\rm exp}(A_1(\vp(z)-{\rm inf}_X\,\vp)) \quad z\in \bar X
\eea
where the constant $C_1$ does not depend on ${\rm inf}\,\vp$, but depends only on the other quantities listed above.
We shall see that the Laplacian $\Delta\vp$ satisfies the same pointwise estimate.

\medskip
Not surprisingly, the constants
${\rm sup}_{M\times [{\rm inf}\,\vp,\i]} F$ and
${\rm sup}_{M\times [{\rm inf}\,\vp,\i]}|\nabla F^{1\over n}|+|\pl_\vp F^{1\over n}|$
in (\ref{C1a}) can be replaced by
${\rm sup}_{M\times [{\rm inf}\,\vp,{\rm sup}\,\vp]} F$
and
${\rm sup}_{M\times [{\rm inf}\,\vp,{\rm sup}\vp]}|\nabla F^{1\over n}|+|\pl_\vp F^{1\over n}|$
respectively. Thus, when $\|\vp\|_{C^0}$ is bounded,
we obtain gradient bounds for $\vp$ for completely general
smooth and strictly positive functions $F(z,\vp)$.
We have however stated them in the above form since we are particularly
interested in the cases when there is no upper bound for ${\rm sup}\,\vp$.
This is crucial for certain applications \cite{PS09a, PS09b}.
If a dependence on $\|\vp\|_{C^0}$
is allowed, then there are many earlier direct approaches.
The first appears to be due to
Hanani \cite{Ha}. More recently, Blocki \cite{B09a}
gave a different proof, and our approach  builds directly on his.
The method of P. Guan \cite{Gp} can be extended to Hessian equations,
while the method of B. Guan-Q. Li \cite{GL} allows a general Hermitian metric $\omega$
as well as a more general right hand side $F(z)\chi^n$,
where $\chi$ is a K\"ahler form.

\medskip
The proof is an application of the maximum principle.
Let $g_{\bar kj}$ and $g_{\bar kj}'$ be the two metrics defined by the K\"ahler forms $\o_0$ and $\o_0+{i\over 2}\ddb\vp$. The covariant derivatives and Laplacians with respect to $g_{\bar kj}$ and $g_{\bar kj}'$ are denoted by $\nabla$, $\Delta$, 
and $\nabla'$, $\Delta'$ respectively. A subindex $g$ or $g'$ will denote 
the metric with respect to which a norm is taken.
It is convenient to introduce the endomorphisms
\bea
\label{relativeend}
h^j{}_k=g^{j\bar p}g_{\bar p k}', \qquad (h^{-1})^j{}_k=(g')^{j\bar p}g_{\bar pk}.
\eea
Their traces are ${\rm Tr}\,h=n+\Delta\vp$ and ${\rm Tr}\,h^{-1}=n+\Delta'\vp$.

\smallskip
As a preliminary, we calculate $\Delta' \log\,|\na \vp|_g^2$,
$|\na\vp|_g^2$ being the expression of interest, and $\Delta'$ 
being the natural Laplacian to use, as it arises from differentiating the Monge-Amp\`ere equation. We have
\bea
\Delta'\log\,|\na\vp|_g^2
={\Delta'|\na\vp|_g^2\over |\na\vp|_g^2}-{|\na|\na\vp|_g^2|_{g'}^2\over |\na\vp|_g^4}.
\eea
If we express $\Delta'$ on scalars as $\Delta'=(g')^{p\bar q}\na_p\na_{\bar q}$, then we can write
\bea
\Delta'|\na\vp|_g^2&=&
\Delta'(\na_m\vp)\na^m\vp+
\na_m\vp\Delta'(\na^m\vp)
+
|\na\na\vp|_{gg'}^2+|\bar\na\na\vp|_{gg'}^2.
\eea
However, making use of the Monge-Amp\`ere equation, we obtain
\bea
\Delta'(\na_m\vp)=(g')^{p\bar q}\na_p\na_{\bar q}\na_m\vp
=
(g')^{p\bar q}\na_m(\na_p\na_{\bar q}\vp)=
\p_m
\log {(\o')^n\over \o^n}=\p_m\log \,F,
\eea
while
\bea
\Delta'(\na^m\vp)=(g')^{p\bar q}\na^m\na_{\bar q}\na_p\vp
+
(g')^{p\bar q}R_{\bar qp}{}^m{}_\ell \na^\ell\vp
=
\na^m\log\,F
+
(h^{-1})^p{}_rR^r{}_p{}^m{}_\ell\na^\ell\vp.
\eea
Thus 
\bea
\Delta'\log \,|\na\vp|_g^2
\geq
{2 {\rm Re}\na_m\log F \na^m \vp
\over |\na\vp|_g^2}
-\Lambda \,{\rm Tr}\,h^{-1}
+
{|\na\na\vp|_{gg'}^2+|\bar\na\na\vp|_{gg'}^2\over
|\na\vp|_g^2}-{|\na|\na\vp|_g^2|_{g'}^2\over |\na\vp|_g^4}\nonumber
\\
\eea
The first term on the right is easily bounded: first write,
\bea
{2 {\rm Re}\na_m\log F \na^m \vp
\over |\na\vp|_g^2}
\geq -2|\na\log\, F|_g{1\over |\na\vp|_g}
=-2 nF^{-{1\over n}}|\na F^{1\over n}|_g{1\over |\na\vp|_g},
\eea
and note that
\bea
|\na F^{1\over n}|_g\leq {\rm sup}_{X\times [0,\infty)}|\p_zF(z,\vp)^{1\over n}|
+
|\na\vp|_g{\rm sup}_{X\times[0,\infty)}|\p_\vp F(z,\vp)^{1\over n}|_g
\equiv
F_1''+|\na\vp|_gF_1',
\nonumber
\eea
while, using the Monge-Amp\`ere equation and the arithmetic-geometric mean inequality,
\bea
nF^{-{1\over n}}\leq {\rm Tr}\,h^{-1}.
\eea
Thus we find
\bea
\Delta'\log \,|\na\vp|_g^2
\geq
-(\Lambda+2F_1'+2{F_1''\over |\na\vp|_g})\,{\rm Tr}\,h^{-1}
+
{|\na\na\vp|_{gg'}^2+|\bar\na\na\vp|_{gg'}^2\over
|\na\vp|_g^2}-{|\na|\na\vp|_g^2|_{g'}^2\over |\na\vp|_g^4}.\nonumber
\\
\eea
The only troublesome term is the negative last term to the right. The key to handling it is a partial cancellation with the two squares preceding it.
This cancellation is rather general, and we formalize in the following lemma:

\begin{lemma}\label{CL2}
Let $X$ be a K\"ahler manifold and $g_{\bar kj}, g_{\bar kj}'$ a pair of K\"ahler metrics on $M$ (not necessarily in the same K\"ahler class). 
Let $\vp\in C^\i(X)$ and define
\bea
S=\langle \na\na \vp,\na \vp\rangle_g,\qquad
T=\langle \na \vp,\bar\na\na \vp\rangle_g, 
\eea
Then we have
\be 
\label{gradientlemma1}
{|\na\na \vp|^2_{gg'}+|\na\bar\na \vp|^2_{gg'}
\over
|\na\vp|_g^2}
\ \geq \ 
{|\na |\na \vp|_g^2|_{g'}^2
\over|\na\vp|_g^4}
-2\Re \langle{\na|\na\vp|_g^2\over |\na\vp|_g^4}
,T\rangle_{g'}+2{|T|_{g'}^2\over |\nabla\vp|_g^4}.
\ee
\end{lemma} 

\medskip 
\noindent
{\it Proof of Lemma \ref{CL2}.} First, we observe that 
for all tensors $A_{pi}$ and $B_{j}$ on $X$,
\be
\label{CS} 
|\langle A,B\rangle_g|_{g'}\ = \ |A_{pi}g^{i\bar j}\overline{B_{j}}|_{g'}\ \leq \ |A|_{gg'}|B|_g.
\ee 
Now $\na|\na\vp|_g^2=S+T$, and applying  (\ref{CS}) to $S$ and $T$ gives:
\bea
|\na \vp|^2_g\cdot(|\na\na \vp|^2_{gg'}+|\na\bar\na \vp|^2_{gg'})
&\geq& |S|_{g'}^2+|T|_{g'}^2
= |\na|\na\vp|_g^2-T|_{g'}^2+|T|_{g'}^2\ 
\nonumber\\
&=&
|\na |\na\vp|_g^2|^2_{g'}-2\Re\langle\na|\na\vp|_g^2,T\rangle_{g'}+2|T|^2_{g'}\
\eea
This proves the inequality (\ref{gradientlemma1}). 

\bigskip

Returning to the problem of $C^1$ estimates,
we can now formulate and prove an important inequality due to Blocki at interior critical points of an expression of the form
\bea
\log\,|\na\vp|_g^2-\gamma(\vp)
\eea
where $\gamma$ is an arbitrary function of a real variable. 
We apply Lemma \ref{CL2} as follows: on the right side
of (\ref{gradientlemma1}), we drop the third term $2|T|^2_{g'}/|\nabla\vp|_g^4$.
In the second term, the tensor T simplifies
upon replacing $\bar\nabla\nabla\phi$ by
$g'-g$, so that $T$ becomes
$T_j = (\na_i\phi)g^{i\bar k}g'_{\bar k j} - \na_j\phi$. We obtain
\bea
\ {|\na\na \vp|_{gg'}^2+|\na\bar\na \vp|^2_{gg'}\over |\na\vp|_g^2}
- {|\na |\na\vp|_g^2|_{g'}^2\over |\na\vp|_g^4} \ 
&\geq&
 \ 
2\Re\langle {\na|\na\vp|_g^2\over |\na\vp|_g^2}, {\na \vp\over |\na\vp|_g^2}\rangle_{g'}\ - \
2\Re\langle {\na|\na\vp|_g^2\over |\na\vp|_g^2}, {\na \vp\over |\na\vp|_g^2}\rangle_g
\nonumber\\
&=&
2 \gamma'(\vp){|\na\vp|_{g'}^2\over |\na\vp|_{g}^2}-
2\gamma'(\vp)
\eea
In the last line, we made use of the fact that
$\na\log |\na\vp|_g^2=\gamma'(\vp)\na \vp$ at an interior critical point of 
the function $\log\,|\na\vp|_g^2-\gamma(\vp)$.

On the other hand,
\bea
-\Delta'\gamma(\vp)
=-\gamma'(\vp)\Delta'\vp-\gamma''(\vp)|\na\vp|_{g'}^2
=\gamma'(\vp)\,{\rm Tr}\,h^{-1}-n\gamma'(\vp)
-
\gamma''(\vp)|\na\vp|_{g'}^2.
\eea
Combining this with the preceding inequality,
we obtain Blocki's inequality \cite{B09a},
\bea
\Delta'(\log|\na\vp|_g^2-\gamma(\vp))
&\geq&
[\gamma'(\vp)-\Lambda-2F_1'-2{F_1''\over |\na\vp|_g}]{\rm Tr}\,h^{-1}
\nonumber\\
&&
\quad
-(n+2)\gamma'(\vp)-\gamma''(\vp)|\na\vp|_{g'}^2
+2 \gamma'(\vp){|\na\vp|_{g'}^2\over |\na\vp|_g^2}.
\eea
The key to the desired estimate is the following choice of $\gamma(\vp)$ \cite{PS09a, PS09c}
\bea
\gamma(\vp)=A\vp-{1\over \vp+C_1}
\eea
where $C_1$ is chosen to be $C_1=-{\rm inf}_X\vp+ 1$, and $A$ is a large positive constant. Then
\bea
A\vp-1\leq \gamma(\vp)\leq A\vp,
\quad
A\leq\gamma'(\vp)\leq A+1,
\quad
\gamma''(\vp)=-{2\over (\vp+C_1)^3}<0
\eea
and we obtain
\bea
\Delta'(\log|\na\vp|_g^2-\gamma(\vp))
\geq
[A-\Lambda-2F_1'-2{F_1''\over |\na\vp|_g}]{\rm Tr}\,h^{-1}
+{2\over (\vp+C_1)^3}|\na\vp|_{g'}^2
-C_2.
\eea

It suffices to show that, at an interior maximum point $p$, the function
$\log |\na\vp|_g^2-\gamma(\vp)$ is bounded by an admissible constant.
We can assume that $|\na\vp(p)|_g^2\geq 1$, otherwise the statement follows trivially from the fact that $\gamma(\vp)\geq A\vp-1$, and $\vp$ is bounded from below.
Choose $A=\Lambda +2F_1'+2F_1''+1$. Then the preceding inequality simplifies further to
\bea
\Delta'(\log|\na\vp|_g^2-\gamma(\vp))
\geq
{\rm Tr}\,h^{-1}+
{2\over (\vp+C_1)^3}|\na\vp|_{g'}^2-C_2.
\eea
At an interior minimum point $p$,
the left hand side is non-positive.
This implies that ${\rm Tr}\, h^{-1}(p)$ is bounded above,
and hence the eigenvalues of $h(p)$ are bounded below by a priori constants.
In view of the Monge-Amp\`ere equation, they are then bounded above and below by a priori constants, since these constants are allowed to depend on ${\rm sup}_XF$. This implies that $|\na\vp|_{g'}^2\geq C_3 |\na\vp|_g^2$, and we obtain
\bea
|\na\vp|_g^2
\leq C_4(\vp+C_1)^3.
\eea
But we can assume that $\log|\na\vp(p)|_g^2-\gamma(\vp(p))\geq 0$, otherwise there is nothing to prove. Thus $\gamma(\vp(p))\leq \log|\na\vp(p)|_g^2$, and hence
\bea
A\vp(p)\leq \gamma(\vp(p))+1
\leq \log|\na\vp(p)|_g^2+1.
\eea
Substituting this in the previous inequality, we find
\bea
|\na\vp(p)|_g^2
\leq C_4(\log|\na\vp(p)|_g^2+C_5)^3.
\eea
This implies that $|\na\vp(p)|_g^2$ is bounded by an a priori constant. The proof of the $C^1$ estimates is complete.

\bigskip
We note that Lemma \ref{CL2} has other uses.
For example, the inequality (\ref{gradientlemma1}) also implies, by completing the square,
\bea
\label{gradientlemma2}
{|\na\na \vp|^2_{gg'}+|\na\bar\na \vp|^2_{gg'}
\over
|\na\vp|_g^2}
\ \geq \ 
{1\over 2}{|\na|\na\vp|_g^2|_{g'}^2\over|\nabla\vp|_g^4}.
\eea
This inequality can be used to simplify several estimates in the K\"ahler-Ricci flow, including the one on the gradient of the Ricci potential.

\section{A Priori Estimates: $C^2$ Estimates}
\setcounter{equation}{0}

The $C^2$ estimates for the complex Monge-Amp\`ere equation (\ref{MA}) are due to Yau \cite{Y78} and Aubin \cite{A}. The precise statement is,

\begin{theorem}
\label{C2} 
Let $\vp$ be a $C^4$ solution of the equation
(\ref{CY}) on a compact K\"ahler manifold $X$,
with smooth boundary $\p X$ (possibly empty).
Then
\bea
0\leq n+\Delta\vp(z)
\leq C\,{\rm exp}\,(A_2\,\,(\vp(z)-{\rm inf}_X\vp))
\eea
where the constant $C$ depends only on an upper bound for $F$, for
${\rm sup}_{X\times[\inf\vp,\infty)}|(\log F)_\vp(z,\vp)|$, 
for the scalar curvature $R$, and a lower bound for
$\Delta_z \log\,F$, for $(\log F)_{\vp\vp}(z,\vp)|\na\vp|^2$,
and for the lower bound
$\Lambda$ 
introduced in (\ref{Lambda}) for the bisectional curvature of $g_{\bar kj}$.

When $\p X$ is not empty, the constant also depends on the boundary value $\vp_b$ of $\vp$,
and on $\|\Delta\vp\|_{C^0(\p X)}$.

The conclusion still holds for $z\in X\setminus Z$, if 
the equation (\ref{MA}) holds on $X\setminus Z$,
$Z$ is a subset of $X$ not intersecting $\p X$, and $\vp(z)\to +\infty$ as $z\to Z$.

\end{theorem}
\smallskip

The derivation of the $C^2$ estimates is particularly transparent if we use the formalism of the relative endomorphisms $h$ of (\ref{relativeend}), as in \cite{PSS} and \cite{PS09a}, which we follow here. As in the proof of the $C^1$ estimates, we would like to estimate ${\rm Tr}\, h$ by the maximum principle. As a preliminary, we calculate
\bea
\Delta' {\rm Tr}\,h=(g')^{p\bar q}\p_{\bar q}\p_p\,{\rm Tr}\,h
&=&
(g')^{p\bar q}{\rm Tr}(\na_{\bar q}'((\na_p' h h^{-1})h)
\nonumber\\
&=&
(g')^{p\bar q}{\rm Tr}(\na_{\bar q}'(\na_p' h h^{-1})h)
+
(g')^{p\bar q}{\rm Tr}(\na_p' h h^{-1}\na_{\bar q}'h).
\eea
But $\na_{\bar q}'(\na_p h h^{-1})=-Rm_{\bar q p}'+Rm_{\bar qp}$, as a special case of the general formula comparing the curvatures of two Hermitian metrics on the same holomorphic vector bundle. Here the full curvature tensors $Rm_{\bar qp}$ and $Rm_{\bar qp}'$ are viewed as endomorphisms on the holomorphic tangent bundle. Thus
\bea
(g')^{p\bar q}{\rm Tr}(\na_{\bar q}'(\na_p' h h^{-1})h)
&=&
-(g')^{p\bar q}R'_{\bar qp}{}^j{}_k h^k{}_j
+
(g')^{p\bar q}R_{\bar qp}{}^j{}_k h^k{}_j
\nonumber\\
&=&
-R'_{\bar mk}g^{k\bar m}+(h^{-1})^p{}_mR^m{}_p{}^j{}_k h^k{}_j.
\eea
But the Ricci curvature $R'_{\bar mk}$ can be obtained from the Monge-Amp\`ere equation
\bea
R'_{\bar mk}=R_{\bar mk}-\p_k\p_{\bar m}\log\,F(z,\vp).
\eea
Thus we obtain
\bea
\label{trh}
\Delta'{\rm Tr}\,h=
-R+\Delta \log\,F(z,\vp)+(h^{-1})^p{}_mR^m{}_p{}^j{}_k h^k{}_j
+
(g')^{p\bar q}{\rm Tr}(\na_p' h h^{-1}\na_{\bar q}'h)
\eea
and hence
\bea
\Delta'\log {\rm Tr}\,h
&=&
{-R+\Delta \log\,F(z,\vp)+(h^{-1})^p{}_mR^m{}_p{}^j{}_k h^k{}_j\over {\rm Tr}\, h}
\nonumber\\
&&
+
\bigg\{
{(g')^{p\bar q}{\rm Tr}(\na_p' h h^{-1}\na_{\bar q}'h)
\over 
{\rm Tr} \,h}
-
{|\na'{\rm Tr}\,h|^2\over ({\rm Tr}\, h)^2}
\bigg\}.
\eea
A fundamental inequality due to Yau and Aubin is that the expression between brackets is non-negative, as a consequence of the Cauchy-Schwarz inequality. Also
\bea
\Delta (\log\,F(z,\vp))=(\Delta_z\log F)(z,\vp)
+(\log F)_\vp({\rm Tr}\,h-n)+(\log F)_{\vp\vp}|\na\vp|^2.
\eea
Thus 
\bea
\Delta'\log {\rm Tr}\,h
\geq -C_1( {\rm Tr}\,h)^{-1}- \Lambda {\rm Tr}\,h^{-1}-C_2
\geq -(C_1+\Lambda){\rm Tr}\,h^{-1}-C_2
\eea
since $({\rm Tr}\,h)^{-1}\leq {\rm Tr}\,h^{-1}$. Here
$C_1, C_2$ depend only on an upper bound for the scalar curvature $R$,
a lower bound for $(\Delta_z \log\,F)(z,\vp)$, a lower bound for $(\log F)_{\vp\vp}|\na\vp|^2$,
and an upper bound for $|(\log F)_\vp|$. We can write now
\bea
\Delta'(\log {\rm Tr}\,h-A_2 \vp)
\geq A_2({\rm Tr}\,h^{-1}-n)-(C_1+ \Lambda) {\rm Tr}\,h^{-1}-C_2
\geq {1\over 2}A_2 \,{\rm Tr}\,h^{-1}-C_3
\eea
for $A_2\geq 2(C_1+\Lambda)$ and $C_3=nA_2$. At a maximum point $z_0$ for
$\log {\rm Tr}\,h-A_2 \vp$, the eigenvalues of $h^{-1}$ are then bounded from above by absolute constants. Equivalently, the eigenvalues $\lambda_i$ of $h$ are bounded from below by absolute constants, and hence, in view of the Monge-Amp\`ere equation
$\prod_{i=1}^n\lambda_i=F$, they are also bounded from above by constants
depending also on ${\rm sup}_XF$. Thus for any $z\in X$,
\bea
\log {\rm Tr}\,h(z)
\leq
\log {\rm Tr}\,h(z_0)+A_2 (\vp(z)-\vp(z_0)).
\eea
This establishes the desired $C^2$ estimates. We note that the $C^0$ bounds for $\Delta \vp$ imply similar $C^0$ bounds for $\p_j\p_{\bar k}\vp$ by plurisubharmonicity, but not for
$\na_j\na_k\vp$.

\section{A Priori Estimates: the Calabi identity}
\setcounter{equation}{0}

To obtain estimates for derivatives of order higher than $2$, we need the
equation to be non-degenerate. Thus we allow constants to depend now on a lower bound for $F$. In particular, the $C^2$ 
estimates imply that the metrics $g_{\bar kj}$ and $g_{\bar kj}'$ are equivalent, 
up to such constants. We restrict ourselves to the equation
(\ref{CYF}), although the arguments can be extended to certain classes of more general $F(z,\vp)$, for example $F(z,\vp)=e^{f(z)\pm\vp}$.

\medskip

In Yau's solution of the Calabi conjecture \cite{Y78}, uniform bounds for the third order derivatives $\na_j\na_{\bar k}\na_m\vp$ were derived from a generalization to the complex case of an identity due to Calabi \cite{Ca2}. We present here a simplified proof of this identity which appeared in \cite{PSS}, and which is again based on the formalism of the relative endomorphism $h^j{}_k=g^{j\bar p}g_{\bar p k}'$. Norms and lowering and raising of indices are with respect to $g_{\bar kj}'$. Covariant derivatives with respect to $g_{\bar kj}$ and $g_{\bar kj}'$ are denoted by $\na$ and $\na'$ respectively.

\medskip
Define 
$S=|\na\bar\na\na\vp|^2$
as in \cite{Y78}.
In terms of $h^\al{}_\b$, 
we have $(g')^{\alpha\bar k}\nabla_j\vp_{\bar k\beta}
=(\nabla_j'h\,h^{-1})^\alpha{}_\beta$, and thus \cite{PSS} 
\bea
S=
|\nabla' h\,h^{-1}|^2.
\eea
The point is that the Laplacian of $S$ can now be evaluated
directly in terms of metrics and curvatures, instead of K\"ahler potentials. We readily find
\bea
\Delta' S&=&
(g')^{m\bar\g}
(\Delta'(\nabla_m'h\,h^{-1})^\b{}_l\overline{(\nabla_\g' h\,h^{-1})_{\bar\beta}{}^{\bar \ell}}
+(\nabla_m'h\,h^{-1})^\b{}_\ell\overline{\bar\Delta'(\nabla_\g' h\,h^{-1})_{\bar\beta}
{}^{\bar\ell}} )
\nonumber\\
&&
+
|\bar\nabla'(\nabla' h\,h^{-1})|^2+|\nabla'(\nabla' h\,h^{-1})|^2
\eea
where 
$|\bar\nabla'(\nabla' h\,h^{-1})|^2\equiv
(g')^{q\bar p}
\nabla_{\bar p}'(\nabla_j'h\,h^{-1})^\al{}_\b
\overline{\nabla_{\bar q}'(\nabla_m'h\,h^{-1})_{\bar\al}{}^{\bar\beta}}$, 
and $\Delta'=(g')^{q\bar p}\na_q'\na_{\bar p}'$,
$\bar\Delta'=(g')^{q\bar p}\na_{\bar p}'\na_q'$.
Commuting the $\nabla_q'$ and the $\nabla_{\bar p}'$ derivatives
gives,
\bea
(\bar\Delta'(\nabla'_j h\,h^{-1}))^\g{}_\al
&=&
(\Delta'(\nabla_j' h\,h^{-1}))^\g{}_\al
-
(R')^\g{}_\mu (\nabla_\g' h\,h^{-1})^\mu{}_\al
+
(R')^\mu{}_\al(\nabla_j' h\,h^{-1})^\g{}_\mu
\nonumber\\
&&
+
(R')^\mu{}_j(\nabla_\mu' h\,h^{-1})^\g{}_\al
\eea
while, in view of the Bianchi identity,
\bea
\Delta' (\nabla_j' h\,h^{-1})^l{}_m
&=&(\nabla')^{\bar p}\pl_{\bar p}(\nabla_j'h\,h^{-1})
=
-(\nabla')^{\bar p}R'_{\bar pj}{}^l{}_m
+
(\nabla')^{\bar p}R_{\bar pj}{}^l{}_m
\nonumber\\
&=&
-\nabla'_j(R')^l{}_m
+
(\nabla')^{\bar p} R_{\bar p j}{}^l{}_m.
\nonumber
\eea
with $R_{\bar pj}{}^l{}_m=-\p_{\bar p}(g^{l\bar q}\p_j g_{\bar qm})$.
Thus we obtain the exact formula
\bea
\label{explicit1}
\Delta' S
&=&|\bar\nabla'(\nabla' h\,h^{-1})|^2+|\nabla'(\nabla' h\,h^{-1})|^2
\nonumber\\
&&
- ((\nabla')^{\bar\g}R'_{\bar\beta\al}
\overline{(\nabla_\g' h\,h^{-1})^{\b\bar\al}}
+
((\nabla')^{\bar\g} h\,h^{-1})^{\b\bar\al}
\overline{\nabla'_\g R'_{\bar\b\ell}}\ )
\nonumber\\
&&
+(\nabla'_m h\,h^{-1})^\b{}_l
((R')^{l\bar\rho}
\overline{((\nabla')^{\bar m} h\,h^{-1})_{\bar\b\rho}}
-
R'_{\bar\rho\b}
\overline{((\nabla')^{\bar m}h\,h^{-1})^{\rho\bar l}}
+
(R')^{m\bar\rho}
\overline{(\nabla_\rho' h\,h^{-1})_{\bar\beta}{}^{\bar l}}
\ )
\nonumber\\
&&
+
(\nabla')^{\bar p}R_{\bar pm}{}^\b{}_l
\overline{(\nabla')^{\bar m} h\,h^{-1})_{\bar\b}{}^{\bar l}}
+
((\nabla')^{\bar \g} h\,h^{-1})_{\bar\mu}{}^{\bar\alpha}
\overline{(\nabla')^{\bar p}R_{\bar p \g}{}^\mu{}_\al}.
\eea
Since $R_{\bar pm}'=R_{\bar pm}-\p_m\p_{\bar p}\log\,F$, it can be viewed as known.
As already noted, the metrics $g_{\bar kj}$ and $g_{\bar kj}'$ are uniformly equivalent. Since the connection $\na'h h^{-1}$ is of order $O(S^{1\over 2})$, 
the above identity implies
\bea
\Delta' S\geq -C_1 S-C_2
\eea
where the constants $C_i$ depend on an upper bound for $\Delta\vp$, 
a lower bound for $F$, the $C^3$
norm of $F$, and the $C^1$ norms of the curvature $R_{\bar kj}{}^l{}_m$ of $g_{\bar kj}$. Using the expression (\ref{trh}) for $\Delta' {\rm Tr}\,h$,
and the fact that
\bea
(g')^{p\bar q}{\rm Tr}(\na_p'h h^{-1}\na_{\bar q}'h)
=
(g')^{p\bar q} (g')^{j\bar k}g^{m\bar r}\na_p\na_{\bar m}\na_j\vp\overline{\na_q\na_{\bar m}\na_k\vp}\geq C_3 S
\eea
since the metrics $g_{\bar kj}$ and $g_{\bar kj}'$ are equivalent,
we readily see that
\bea
\Delta'(S+ A\,{\rm Tr}\,h)
\geq C_4 S-C_5
\eea
for $A$ sufficiently large. We can now apply the maximum principle 
and obtain the following:

\begin{theorem}
\label{C3}
Let $\vp$ be a $C^5$ solution of the equation (\ref{MA})
on a compact K\"ahler manifold $(X,\o_0)$ with smooth boundary $\p X$
(possibly empty). Then that
$S$ is uniformly bounded by constants depending only on the $C^0$ norms
of $\p_j\p_{\bar k}\log\, F$ and $\na_j\na_{\bar k}\na_m\log F$, a lower bound for $F$, the $C^1$ norm of $R_{\bar kj}{}^l{}_m$,
and, when $\p X$ is not empty, on $\|S\|_{C^0(\p X)}$.
\end{theorem}

\medskip
We conclude this section by noting that most of the a priori estimates discussed here have counterparts for parabolic Monge-Amp\`ere equations. They were instrumental in Cao's proof of the all-time existence for the K\"ahler-Ricci flow
on manifolds of definite Chern classes \cite{Cao}. They also apply in many situations to manifolds of general type \cite{Ts, EGZ, ST09}. For the modified K\"ahler-Ricci flow, the $C^3$ estimates and the Calabi identity require an additional argument \cite{PSSW2}, as well as a full use of the square terms
in (\ref{explicit1}) which were dropped in the proof of Theorem \ref{C3}. Extensions to flows on Hermitian manifolds can be found in \cite{Gm, ZZ}.

\section{Boundary Regularity}
\setcounter{equation}{0}

In this section, we discuss a priori estimates for
the Dirichlet problem for the complex Monge-Amp\`ere equation on a K\"ahler  manifold $(X,\o_0)$ with smooth boundary $\p X$.

\subsection{$C^0$ estimates}

Let $\vp_b$ be a smooth function on $\p X$, and consider the Dirichlet problem
\bea
\label{Dirichlet0}
(\o_0+{i\over 2}\ddb \vp)^n=F(z,\vp,\na\vp)\o_0^n
\ \ {\rm on}\ \ X,
\qquad
\vp=\vp_b\ \ {\rm on}\ \ \p X,
\eea
where $n={\rm dim}\,X$, $F$ is a smooth strictly positive function, and $\vp\in \PSH(X,\o_0)\cap C^\infty(X)$.

The fact that $\vp\in \PSH(X,\o_0)$ implies that $n+\Delta\vp\geq 0$. If $h$ is the solution of the Dirichlet problem $\Delta h=-n$ on $X$, $h=\vp_b$ on $\p X$, then the comparison principle implies
\bea
\vp\leq h.
\eea
Thus, to obtain $C^0$ estimates, we need only a lower bound for $\vp$. As shown by Caffarelli, Kohn, Nirenberg, and Spruck \cite{CNS, CKNS}, this can be effectively obtained if we assume the existence of a smooth subsolution $\vps$ 
of the Dirichlet problem (\ref{Dirichlet1}), that is, a smooth function $\vps$ satisfying
\bea
\label{Dirichlet2}
(\o_0+{i\over 2}\ddb \vps)^n>F(z,\vps,\na\vps)\o_0^n
\ \ {\rm on}\ \ X,
\qquad
\vps=\vp_b\ \ {\rm on}\ \ \p X.
\eea
Indeed, in the method of continuity, the problem reduces to a priori estimates for the equation
\bea
\label{Dirichlet3}
(\o_0+{i\over 2}\ddb \vp)^n&=&
tF(z,\vp,\na\vp)\o_0^n+(1-t)(\o_0+{i\over 2}\ddb\vps)^n
\ \ {\rm on}\ \ X,
\nonumber\\
\vp&=&\vp_b\ \ {\rm on}\ \ \p X,
\eea
for $0\leq t\leq 1$. Let $\vp=\vps$ for $t=0$. We claim that, if a smooth solution exists in an interval $0\leq t<T$, then
\bea
\vps < \vp\ \ {\rm in} \ \ X,
\eea
for all $t<T$. To see this, note that the derivative in $t$ of ${(\o_0+{i\over 2}\ddb\vp)^n\over (\o_0+{i\over 2}\ddb\vps)^n}$ is strictly negative at $t=0$. Thus
$(\o_0+{i\over 2}\ddb \vp)^n<(\o_0+{i\over 2}\ddb\vps)^n$,
and $\vps<\vp$  for $t$ strictly positive and small, by the comparison principle. 
If there exists $t_0$, $0<t_0<T$, with $\vps(z_0)=\vp(z_0)$ for some $z_0\in X$,
let $t_0$ be the first such time. By continuity, $\vps(z)\leq \vp(z)$
for all $z\in X$ and $t=t_0$, so $z_0$ is a maximum of the function $\vps-\vp$ at $t_0$. In particular, at $t_0$ and $z_0$, we have $\na\vps=\na\vp$ and
\bea
(\o_0+{i\over 2}\ddb\vps)^n\leq (\o_0+{i\over 2}\ddb\vp)^n.
\eea
But the equation (\ref{Dirichlet3}) implies, again at $t_0$ and $z_0$,
\bea
(\o_0+{i\over 2}\ddb\vp)^n=
tF(z,\vps,\na\vps)+(1-t)(\o_0+{i\over 2}\ddb\vps)^n
<(\o_0+{i\over 2}\ddb\vps)^n,
\eea
which is a contradiction.

\subsection{$C^1$ boundary estimates}

The $C^1$ estimates at the boundary $\p X$ follow from
the bounds $\vps\leq u\leq h$, and the fact that all three functions have the same boundary values. When the right hand side $F(z,\vp,\na\vp)$
does not depend on $\na\vp$, the estimates established earlier in Section 4 show that the interior $C^1$ estimates can be reduced to the boundary $C^1$ estimates.

\subsection{$C^2$ boundary estimates of Caffarelli-Kohn-Nirenberg-Spruck and B. Guan}

The barrier constructions of Caffarelli, Kohn, Nirenberg, Spruck \cite{CKNS} and B. Guan \cite{Gb} provide $C^0(\p X)$ bounds for $\Delta \vp$, in terms of $C^0(X)$ bounds for $\vp$ and for $\na\vp$. The following slightly more precise formulation of their estimates can be found in \cite{PS09a}, under the simplifying assumption that the boundary $\p X$ is holomorphically flat \footnote{A hypersurface $\p X$ is holomorphically flat if, locally, there exist holomorphic coordinates $(z_1,\cdots, z_n)$ so that $\p X$ is given by ${\rm Re}\,z_n=0.$}:

\begin{theorem}
\label{C2b} 
Assume that $\p X$ is holomorphically flat, and that $\vp$ is a $C^3$ solution of the equation (\ref{Dirichlet0}), with $F(z)$ on the right-hand side. 
Then we have
\bea
{\rm sup}_{\p X}
(n+\Delta \vp)\leq C\,
{\rm sup}_{\p X}(1+|\na\vp|^2)\,{\rm sup}_X(1+|\na\vp|^2),
\eea
for a constant $C$ depending only on the boundary $\p X$, $\o_0$, 
and upper bounds for ${\rm sup}_X F$, and
${\rm sup}_X(\na \log\,F)$, and $\|\vp\|_{C^0(X)}$, $\|\na \vp\|_{C^0(X)}$.
\end{theorem}

By the interior estimates of Yau and Aubin in \S 6, the uniform bound for $\Delta\vp$ in the whole of $X$ can be reduced to its estimate on $\p X$. Thus 
the above bound implies that $\|\Delta\vp\|_{C^0(X)}$ is bounded
in terms of the constants indicated.
By plurisubharmonicity, it follows that  all the mixed partials $\|\p_j\p_{\bar k}\vp\|_{C^0(X)}$ are bounded as well.

\medskip
It is an interesting question whether bounds for the un-mixed partials
$\|\na_j\na_k\vp\|_{C^0}$ can be obtained as well without additional assumptions. 
Such bounds have been obtained by Blocki \cite{B09b} under the additional assumption that the background form $\o_0$ has non-negative bisectional 
curvature. 

\medskip
If we allow bounds to depend on a lower bound for $F$, then the equation
(\ref{Dirichlet2}) can be viewed as uniformly elliptic, since the eigenvalues of
the relative endomorphism $h^j{}_k=g^{j\bar p}g_{\bar p k}'$ are already known to be bounded from above, and using the lower bound for $F$, they are also bounded from below. The Monge-Amp\`ere equation is concave, so we can then apply to the following general theorem of the Evans-Krylov and Krylov theory, which we quote from Chen-Wu \cite{CW} (see also Gilbarg-Trudinger \cite{GT}
p. 482 and Q. Han \cite{H}). The statement is local, and can be formulated for domains with smooth boundary in ${\bf R}^n$:

\begin{theorem}
\label{EvansKrylov}
Assume that $\Omega\subset{\bf R}^n$ has smooth boundary, and the boundary data is smooth. Assume that $F(x,u,Du,D^2u)$ is smooth in all variables
$(x,u,p,A)$, uniformly elliptic and concave (or convex) in $D^2u$, and assume that $\|u\|_{C^{1,\g}(\bar\Omega)}$
is bounded for some $0<\g<1$. Then there are constants $0<\alpha<\g$ and $C$
so that, for any $0<\beta<\alpha$, we have
\bea
\|u\|_{C^{2,\beta}(\bar\Omega)}\leq C.
\eea
\end{theorem}

We note that, while both the local \cite{Ca} and the global \cite{TrWa}
$C^{2,\alpha}$ regularity is known for
real Monge-Amp\`ere equations when the right hand side $F$ is in $C^\alpha$,
the corresponding question is still not completely resolved in the complex case. For some recent progress on this issue,
see \cite{DZZ}, and particularly \cite{W2}, where
it is shown that the solution $\vp$ is of class $C^{2,\alpha}$ if the right hand side
$F(z)$ is strictly positive, $F^{1\over n}\in C^\alpha$, and $\Delta\vp$ is bounded.

\section{The Dirichlet Problem for the
Monge-Amp\`ere equation}
\setcounter{equation}{0}

The preceding a priori estimates imply the following classic existence theorem due to Caffarelli, Kohn, Nirenberg, and Spruck \cite{CKNS} and B. Guan \cite{Gb}:

\begin{theorem}
\label{CKNS-G}
Let $(\bar X,\o_0)$ be a compact K\"ahler manifold of dimension $n$, with smooth boundary $\p X$. Let $F(z,\vp)$ be a smooth, strictly positive function
of the variables $z$ and $\vp$, and let $\vp_b$ be a smooth function on $\p X$. Consider the Dirichlet problem
\bea
\label{Dirichlet1}
(\o_0+{i\over 2}\ddb \vp)^n=F(z,\vp)\,\o_0^n,
\qquad
\vp=\vp_b\ {\rm on}\ \p X.
\eea
If $F_\vp(z,\vp)\geq 0$ and the problem admits a smooth subsolution, that is, a smooth function $\vps$
satisfying
\bea
\label{sub1}
(\o_0+{i\over 2}\ddb \vps)^n>F(z,\vps)\,\o_0^n,
\qquad
\vps=\vp_b\ {\rm on}\ \p X,
\eea
then the Dirichlet problem (\ref{Dirichlet1}) admits a unique solution $\vp$, and  $\vp\in C^\infty(\bar X)$.
\end{theorem}

Indeed, the $C^0$ estimates of \S 8.1, the $C^2$ estimates of \S 6 and \S 8.2, the Evans-Krylov theory for higher derivatives of \S 8.2 can be applied to show that the equation (\ref{Dirichlet3}) admits a solution for $0\leq t\leq 1$. 

\medskip
Similar results for the real Monge-Amp\`ere equation can be found in \cite{GS} and \cite{Gb98}. An extension to Hermitian manifolds can be found in \cite{GL}.

\medskip
The a priori estimates show more than just the existence of a solution $\vp$ for
the equation (\ref{Dirichlet1}): the upper bound for $\Delta \vp$ does not depend on a lower bound for $F(z,\vp)$. This allows an immediate application to the existence of solutions to the Dirichlet problem for the completely degenerate, or homogeneous, complex Monge-Amp\`ere equation.
For this, we apply Theorem \ref{CKNS-G} to the Dirichlet problem
\bea
(\o_0+{i\over 2}\ddb \vp_\e)^n
=\e\,\o_0^n,
\qquad
\vp_\e=0\ {\rm on}\ \p X,
\eea
where $\e$ is a constant satisfying $0<\e<1$. The function $\vps_\e=0$ is a subsolution, and hence Theorem \ref{CKNS-G} implies the existence of a smooth solution $\vp_s$ with
$\Delta\vp_s$ bounded uniformly in $\e$. Thus a
subsequence of the functions $\vp_\e$ converges
in $C^{1,\alpha}$ to a $C^{1,\alpha}$ solution
of the equation (\ref{HCMA}) for all $0<\alpha<1$. We obtain in this manner 
the following theorem, whose present formulation is due to Blocki \cite{B09b} 
and which generalizes the theorem of Chen \cite{C00}
stated further below as Theorem \ref{Chen}:

\begin{theorem}
\label{HCMA1}
Let $(X,\o_0)$ be a compact K\"ahler manifold with smooth boundary $\p X$.
Then the Dirichlet problem
\bea
\label{HCMA}
(\o_0+{i\over 2}\ddb \vp)^n=0
\ {\rm on}\ \p X,
\qquad
\vp=0\ {\rm on}\ \p X
\eea
admits a unique solution, which is of class $C^{1,\alpha}(\bar X)$ for each
$0<\alpha<1$.
\end{theorem}

In some applications, as in the problem of geodesics in the space of K\"ahler potentials described below in Section \S 13, it is actually necessary to consider equations
of the form (\ref{HCMA}), but with the K\"ahler form
$\o_0$ replaced by a smooth background $(1,1)$-form $\o$ which is closed, non-negative, but not strictly positive. We discuss a specific situation where the existence and regularity of solutions can still be established by the a priori estimates that we described in sections \S 5, \S 6, \S 8.2, and \S 8.3 (in fact, some of the $C^1$ estimates given in \S 5 were designed for that purpose).

\medskip
Assume that $\o$ is a smooth, closed, and non-negative $(1,1)$-form, and that there exists an effective divisor $E$, not intersecting $\p X$, with the line bundle $O(E)$ admitting a metric $K$ satisfying
\bea
\label{omegaH}
\o_K\equiv \o+\delta{i\over 2}\ddb \log K>0.
\eea
for some strictly positive constant  $\delta$. Then we have the following theorem \cite{PS09a}:

\begin{theorem}
\label{HCMA2}
Let $X$ be a compact complex manifold with smooth boundary $\p X$. Assume that
$\o$ is a smooth non-negative $(1,1)$-form, $E$ is an effective divisor
not intersecting $\p X$, $K$ is a metric on $O(E)$, with $\o_K$ satisfying
the K\"ahler condition (\ref{omegaH}). Then the Dirichlet problem
\bea
\label{HCMAeq}
(\o+{i\over 2}\ddb\vp)^n=0
\ {\rm on}\ X,
\qquad
\vp=0\ {\rm on}\ \p X
\eea
admits a unique bounded solution. The solution is $C^\alpha(\bar X\setminus E)$
for any $0<\alpha<1$. If $\p X$ is holomorphically flat (in the sense that there exists holomorphic coordinates $z_i$ with $\p X=\{{\rm Re}\,z_n=0\}$ locally), then
the solution is $C^{1,\alpha}(\bar X\setminus E)$ for any $0<\alpha<1$.
\end{theorem}

We sketch the proof. Let 
\bea
\omega_s=(1-s)\omega_0+s\omega_K.
\eea
For $0<s<1$, $\omega_s$ is strictly positive definite. Consider the equation
\bea
\label{HCMAeq1}
(\o_s+{i\over 2}\ddb\vp_s)^n=F_s(z)\o_s^n\ \ {\rm on}\ \ X,
\qquad
\vp_s=0\ \ {\rm on}\ \ \p X,
\eea
for some smooth functions $F_s>0$ satisfying ${\rm sup}_XF_s<1$, to be specified more completely later. By Theorem \ref{CKNS-G}, this equation admits a smooth solution in $PSH(X,\o_s)$ for each $s>0$. Since the eigenvalues of $\o_s$ are bounded from above with respect to the K\"ahler form $\o_K$, the $\o_s$-plurisubharmonicity of $\vp_s$ implies that $\Delta_{\o_K}\vp_s\geq -C$
for a uniform constant $C$. The arguments for $C^0$ estimates in \S 8.1 imply that
the norms $\|\vp_s\|_{C^0(X)}$ are uniformly bounded in $s$.

To obtain $C^1$ estimates on compact subsets of $X\setminus E$, we choose $F_s$ as follows. First, define
\bea
\hat\o_s=\o_K+s{i\over 2}\ddb \log K^\delta.
\eea
Then $\hat\o_s$ is uniformly bounded from below for all $s$ sufficiently small. In particular, its curvature tensor is uniformly bounded together with all its derivatives.
On the other hand, since $\hat\o_s$ can also be expressed as
\bea
\hat\o_s=
\o_s+{i\over 2}\ddb \log K^{\delta}.
\eea
the equation (\ref{HCMAeq1}) can be rewritten as
\bea
(\hat\o_s+{i\over 2}\ddb (\vp_s-\delta\log \|\psi\|_{K}^2))^n
=\hat F_s\hat\o_s^n,
\eea
with $\psi$ a holomorphic section of $O(E)$,
$\|\psi\|_{K}^2=\psi\bar \psi K$,
and $\hat F_s\hat\o_s^n=F_s\o_s^n$. Choose $F_s$ to be constants tending so fast to $0$ that ${\rm limsup}_{s\to 0}\|F_s\|_{C^0(X)}=0$. The desired uniform bounds for $\na\vp_s$ on compact subsets of $X\setminus E$ follow from
the $C^1$ estimates of \S 5. With these estimates, it is then easy to show the
existence of a subsequence of $\vp_s$ converging in $C^\alpha$ on compact subsets of $X\setminus E$ to a solution of (\ref{HCMAeq}).

\section{Singular Monge-Amp\`ere equations}
\setcounter{equation}{0}

In the seminal paper \cite{Y78}, Yau not only solved the Calabi conjecture, but he also started the study of complex Monge-Amp\`ere equations in more general settings. These include settings when the right hand side may have zeroes or poles, or when the manifold $X$ is not compact and one looks for a complete K\"ahler-Einstein metric, or when $X$ is quasi-projective. We recall briefly some of these classical results below, before discussing some more recent developments.  
In these more recent developments, the underlying manifold may have singularities, and/or the background form $\o_0$ in the Monge-Amp\`ere equation may be degenerate.

\subsection{Classic works}

The classical literature on singular Monge-Amp\`ere equations and singular K\"ahler-Einstein metrics is particularly rich, as different equations are required by different geometric situations. We shall restrict ourselves to describing three results.

\medskip
First, we consider the case when the right hand side of the equation has zeroes and/or poles.
Let $(X,\o_0)$ be a compact K\"ahler manifold. Let $\{ L_i\}_{i=1}^I$ be a family of holomorphic line bundles over $X$. For each $i$, let $s_i$ be a holomorphic section of $L_i$ and $h_i$ a smooth hermitian metric on $L_i$. 
Let 
\bea
u_k = \sum_{i=1}^I a_{k,i} |s_i|_{h_i}^{2\alpha_{k,i}},
\ 
k=1, ..., K,
\qquad
v_l = \sum_{i=1}^I b_{l, i} |s_i|_{h_i}^{2\beta_{l,i}}
\ l=1, ..., L, 
\eea
where $a_{k, i}$ $b_{l, i}$, $\alpha_{k,i}$ and $\beta_{l, i}$ are nonnegative numbers, and consider the following Monge-Amp\`ere equation 
\begin{equation}
\label{yaueqn}
(\omega_0 + {i\over 2}\ddbar \varphi)^n = \frac{u_1 u_2... u_K}{ v_1 v_2 ... v_L} e^{f(z)} \omega_0^n,
\end{equation}
where $f=f(z)$ is a smooth function on $X$,. The following theorem is due to Yau \cite{Y78}:

\begin{theorem} 
\label{yausing} 
Assume the following two conditions:

\smallskip

{\rm (1)} $(u_1 u_2 ... u_K)(v_1 v_2 ... v_L)^{-1} e^f \in L^n (X)$ and $$\int_X \frac{u_1 u_2 ... u_K}{v_1 v_2 ... v_L} \,e^f\, \omega_0^n = \int_X \omega_0^n,$$

{\rm (2)} there exists $\epsilon>0$ such that 
$$ (v_1 v_2 ... v_L)^{-\epsilon} | \Delta \log (v_1 v_2... v_L)|^{(n-1)/n} \in L^1 (X\setminus D),$$
where $\Delta$ is the Laplacian with respect to $\omega_0$ and $D$ is the union of the zeros of $v_l$, $l=1,..., L$.

\smallskip

Then there exists a bounded $\omega_0$-psh function $\varphi$ solving the equation (\ref{yaueqn}). Furthermore, $\varphi$ is smooth outside the zeros of $u_k$ and $v_l$, for $k=1, ..., K$ and  $l=1, ..., L$, and $\varphi$ is unique up to a constant. 

\end{theorem}

In particular, if for each $k$ and $l$, $u_k = a_{i_k} |s_{i_k}|_{h_{i_k}}^{2 \alpha_{i_k}}$ and $v_l = b_{i_l} |s_{i_l}|_{h_{i_l}}^{2\beta_{i_l}}$ for some $1\leq i_k\leq K$ and $1\leq i_l \leq I$, the second assumption in Theorem \ref{yausing} holds automatically. The first assumption for Theorem \ref{yausing} is that the right hand side of the equation (\ref{yaueqn}) is in $L^n(X)$. Thus the theorem on $C^0$ estimates of Kolodziej can also be applied here, and we can obtain in this manner 
a new proof of Theorem \ref{yausing}.

\medskip
The next important geometric situation is that of open complex manifolds. There one is interested in complete K\"ahler-Einstein metrics of negative curvature. 
In \cite{CY80}, Cheng and Yau gave effective criteria for the existence of such metrics. In particular, they proved the existence of a complete K\"ahler-Einstein metric of negative scalar curvature on bounded, smooth, strictly pseudoconvex domains in $\C^n$.  This corresponds to solving Monge-Amp\`ere equations of the form (\ref{CY}), with the solution tending to $\infty$ at the boundary.
This also allowed Cheng and Yau to obtain essentially sharp boundary regularity results for the Dirichlet problem for the closely related equation $J(u)=1$ of Fefferman \cite{F76}. It was subsequently shown by Mok and Yau \cite{MY83} that any bounded domain of holomorphy admits a complete K\"ahler-Einstein metric.

\medskip
A third important class of non-compact manifolds is the class of quasi-projective manifolds. Let $M=\overline{M}\setminus D$ be a quasi-projective manifold, where $\overline{M}$ is a projective manifold and $D$ is a smooth ample divisor on $M$. The Calabi conjecture for quasi-projective manifolds asserts that,
 for any smooth real valued $(1,1)$-form $\eta \in c_1(K_{\overline{M}}^{-1} \otimes [D]^{-1})$, there exists a complete K\"ahler metric on $M$ with its Ricci curvature equal to $\eta|_M$. This was proved by Tian and Yau in \cite{TY86, TY90, TY91}.

\medskip
These works required at that time many new technical tools which remain useful to this day. They include the notion of bounded geometry, the Cheng-Yau H\"older spaces with weights, and particularly the observation repeatedly stressed in these works that the arguments are almost local in nature, and that the manifold can be allowed singularities, as long as the metric admits a non-singular resolution by a local holomorphic map.

\subsection{Monge-Amp\`ere equations on normal projective varieties}

The original theorems of Yau \cite{Y78} and Yau \cite{Y78} and Aubin \cite{A}
establish the existence of K\"ahler-Einstein metrics on a K\"ahler manifold $X$ when $K_X$ has zero or positive first Chern class.
We discuss now one of the new developments in the theory of complex Monge-Amp\`ere equations, namely an extension of these results to normal projective manifolds.
Normal projective manifolds are a very specific class of manifolds
with singularities.
For the convenience of the reader, we summarize here some of their basic definitions and properties.

\medskip
Let $X$ be a subvariety of ${\bf CP}^N$. 
A function on a neighborhood of a point $z_0\in X$ is holomorphic if it extends to a holomorphic function on a neighborhood of $z_0\in {\bf CP}^N$. Let $X_{sing}$ be the smallest subset of $X$ with $X\setminus X_{sing}$
a complex manifold. Then $X$ is said to be normal if for any $z_0\in X_{sing}$,
there is a neighborhood $U$ of $z_0$ so that any  bounded holomorphic function on 
$U\setminus X_{sing}$ extends to a holomorphic function on $U$.

\smallskip
A plurisubharmonic function on $U\subset X$ is by definition the restriction
to $X$ of a plurisubharmonic function in a neighborhood $\hat U$ of $U$ in ${\bf CP}^N$. By a theorem of Fornaess and Narasimhan \cite{FN},
if $X$ be a normal projective variety, and a function $\vp$
is plurisubharmonic on $U\setminus X_{sing}$ and is bounded, then $\vp$ is plurisubharmonic on $U$.

\smallskip
A line bundle $L$ on $X$ is an ample ${\bf Q}$-line bundle if 
$mL$ is the restriction to $X$ of $O(1)$ for some $m\in {\bf Z}^+$.
More generally, a line bundle $\tilde L\to\tilde X$
is an ample ${\bf Q}$-line bundle if there is an imbedding of $\tilde X$ into
projective space, with the pull-back of $O(1)$ equal to $m\tilde L$ for some $m\in{\bf Z}^+$.

\medskip

We can define now the notion of Monge-Amp\`ere measure on a normal projective variety $X$. Let ${\rm dim}\,X=n$, and let
$\pi: \tilde X \rightarrow X$ be a smooth resolution of singularities of $X$. Let $\tilde L = \pi^* L$ for any ample ${\bf Q}$-line bundle $L\rightarrow X$. By definition, $mL = O (1)$ for some $m \in {\bf Z}^+$ and so $$ m \tilde L = \pi^* O (1).$$
Let $m\omega$ be the restriction of the Fubini-Study metric on  ${\bf CP}^M$ to $X$ and let $\tilde \omega = \pi^* \omega$. For any bounded $\omega$-plurisubharmonic function $\varphi$ on $X$,  we let $\tilde \varphi = \pi^* \varphi$.
The measure $(\tilde\omega+{i\over 2} \ddbar \tilde \varphi)^n $ is a well-defined Monge-Amp\`ere measure on $\tilde X$. Since 
$\tilde \varphi$ is bounded and $\tilde\omega$-plurisubharmonic,
by the work of Bedford and Taylor \cite{BT76}, it puts no mass on
the exceptional locus $\pi^{-1} (X_{sing})$. Furthermore, 
\bea
\int_{X\setminus X_{sing}} (\omega+ {i\over 2}\ddbar \varphi)^n = \int_{\tilde X} (\tilde\omega + {i\over 2}\ddbar\tilde\varphi)^n= \int_{\tilde X} \tilde \omega^n<\infty.
\eea 
Therefore,  $(\tilde\omega + {i\over 2}\ddbar\tilde \varphi)^n$ can be pushed forth to a measure on $X$ and it coincides with the trivial extension of $(\omega+ {i\over 2} \ddbar \varphi)^n$ from $X\setminus X_{sing}$ to $X$. In particular, the Monge-Amp\`ere measure $(\omega+ {i\over 2}\ddbar \varphi)^n$ on $X$ does not depend on the resolution of singularities. 

\medskip

Let $\Omega$ be a smooth real valued semi-positive $(n, n)$-form on $X\setminus X_{sing}$ and let $\tilde \Omega= \pi^* \Omega$. We consider the following Monge-Amp\`ere equation on $X$
\begin{equation}
\label{masing}
(\omega + {i\over 2}\ddbar \varphi)^n = e^{\alpha \varphi + F} \Omega,
\end{equation}
where $\alpha = -1, 0 , 1$ and $F$ is a real valued function on $X$.
This Monge-Amp\`ere equation can be lifted to a Monge-Amp\`ere equation on $\tilde X$ if the solution $\varphi$ is bounded and $\omega$-plurisubharmonic, i.e., $\omega + {i\over 2}\ddbar \varphi \geq 0$. More precisely, let $\tilde \varphi = \pi^* \varphi$, and consider the Monge-Amp\`ere equation on $\tilde X$ 
\begin{equation}\label{mareso}
(\tilde \omega + {i\over 2}\ddbar \tilde \varphi)^n = e^{\alpha \tilde \varphi + \tilde F} \tilde \Omega,
\end{equation}
where $\tilde F = \pi^* F$.
If $\varphi$ is a bounded $\omega$-plurisubharmonic solution of  equation (\ref{masing}), then $\tilde \varphi=\pi^* \varphi$ is bounded and $\tilde\omega$-plurisubharmonic, and it solves  equation (\ref{mareso}).  Now let us assume that $\tilde \varphi$ is a bounded $\tilde\omega$-plurisubharmonic solution of the equation (\ref{mareso}). Any fibre of $\pi$ over a singular point of $X$ is connected by Zariski's connectedness theorem, and $\tilde\omega=0$ when restricted to the fibre.  Therefore $\tilde \varphi$ is constant along the fibre because it is plurisubharmonic and bounded.  Hence $\tilde \varphi$ descends to a function $\varphi$ on $X$. The function $\varphi$ is bounded and $\omega$-plurisubharmonic on $X\setminus X_{sing}$, and so it is $\omega$-plurisubharmonic function on $X$ because $X$ is normal. We have thus shown
 
\begin{lemma} 
The equation (\ref{masing}) admits a bounded $\omega$-plurisubharmonic solution $\varphi$ if and only if the equaiton (\ref{mareso}) admits a bounded $\tilde\omega$-plurisubharmonic solution $\tilde\varphi$.
\end{lemma}

Therefore, we can solve equation (\ref{mareso}) on a smooth manifold $\tilde X$ instead of solving equation (\ref{masing}) on a singular variety $X$. Furthermore the construction is resolution independent because given any two resolutions, we can move the measures to the same resolution and apply the uniqueness property of Monge-Amp\`ere equations there. So we obtain
the following lemma, which follows immediately from Theorem \ref{C0}:
 
\begin{lemma}
\label{masingsol} Let $\Theta$ be a smooth volume form on $\tilde X$. Then the equation (\ref{masing}) admits a bounded and $\omega$-plurisubharmonic solution for  $\alpha = 0 $,  if $\frac{e^{\tilde F} \tilde \Omega }{\Theta}\in L^p(\tilde X)$ for some $p>1$.
\end{lemma}

In fact, Lemma \ref{masingsol} also holds for $\alpha=1$ by \cite{EGZ}. We can apply it now to solving K\"ahler-Einstein equations on singular varieties. Recall some basic definitions for canonical models of general type:

\begin{definition}

{\rm (a)}
A projective variety $X$ is said to be a canonical model of general type if $X$ is a normal and the canonical divisor $K_X$ is an ample ${\bf Q}$-line bundle. 

{\rm (b)}
Let $X$ be a canonical model of general type. A form $\Omega$ is said to be a smooth volume form on $X$ if for  any point  $z\in X$, there exists an open  neighborhood $U$ of $z$ such that   $$\Omega = f_U (\eta\wedge \overline{\eta})^{\frac{1}{m}},$$
where $f_U$ is a smooth positive function on $U$ and $\eta$ is a local generator of $mK_X$ on $U$. In particular, any  smooth volume $\Omega$ induces a smooth hermitian metric $h= \Omega^{-1}$ on $K_X$.

{\rm (c)} $X$ is said to be a canonical model of general type with canonical singularities if for any resolution of singularities $\pi: \tilde X \rightarrow X$ and any smooth volume form $\Omega$ on $X$, 
\bea
\tilde \Omega = \pi^* \Omega
\eea
is a smooth real valued $(n,n)$-form on $\tilde X$. 

\end{definition}

We can now describe some recent results of Eyssidieux, Guedj, and Zeriahi \cite{EGZ} on the existence of K\"ahler-Einstein metrics of zero
or negative curvature on manifolds with canonical singularities. As above, let $\pi: \tilde X \rightarrow X $ be a resolution of singularities, let $m\omega$ be the restriction of the Fubini-Study metric of ${\bf CP}^M$ on $X$,
and let $\Omega$ be a smooth volume form $\Omega$ on $X$ such that 
\bea
{i\over 2}
\ddbar \log \Omega = \omega.
\eea
The following theorem on K\"ahler-Einstein metrics with negative curvature
was proved in \cite{EGZ}, using the $C^0$ estimates
of Theorem \ref{C0}:

\begin{theorem}\label{egz1} Let $X$ be a canonical model of general type with canonical singularities. Then there exists a unique bounded and $\omega$-plurisubharmonic function $\varphi$ solving the following Monge-Amp\`ere equation on $X$

\begin{equation}\label{98}
(\omega+ {i\over 2}\ddbar \varphi)^n = e^\varphi \Omega.
\end{equation}
In particular, on $X\setminus X_{sing}$,  $\omega_{KE}=\omega+ {i\over 2} \ddbar \varphi$ is smooth  and $$Ricci(\omega_{KE}) = - \omega_{KE}.$$

\end{theorem}

Next we discuss the case of zero curvature. Recall that
$X$ is said to be a Calabi-Yau variety if $X$ is a projective normal variety and $mK_X$ is a trivial line bundle on $X$ for some $m \in {\bf Z}^+$. 
Since $mK_X$ is a trivial line bundle on $X$, there exists a constant global section $\eta$ of $mK_X$. Let $\Omega_{CY} = (\eta\wedge \eta)^{\frac{1}{m}}$. Then $\Omega$ is a smooth volume form on $X$. Furthermore, ${i\over 2}\ddbar\log \Omega_{CY} = 0$.

\begin{definition} 
A Calabi-Yau variety $X$ is said to be a Calabi-Yau variety with canonical singularities if for any resolution of singularities $\pi: \tilde X \rightarrow X$, 
$$\tilde \Omega_{CY} = \pi^* \Omega_{CY}$$ is a smooth real valued $(n,n)$-form on $\tilde{X}$.
\end{definition}

Let $X$ be a Calabi-Yau variety with canonical singularities. We choose the smooth K\"ahler form $\omega_L \in c_1(L)$ induced from the Fubini-Study metric on ${\bf CP}^N$ in the same way as in the earlier discussion. Then we have the following theorem due to \cite{EGZ}
 
\begin{theorem} Let $X$ be a Calabi-Yau variety with canonical singularities.  Then for any ample ${\bf Q}$-line bundle, there exists a unique bounded and $\omega_L$-plurisubharmonic function $\varphi$ solving the following Monge-Amp\`ere equation on $X$
\begin{equation}
(\omega_L+ {i\over 2}\ddbar \varphi)^n = c_L \Omega_{CY},
\end{equation}
where $c_L \int_X \Omega_{CY}= \int_X \omega_L^n$. 
In particular, on $X\setminus X_{sing}$,  $\omega_{CY}=\omega_L+{ i\over 2} \ddbar \varphi$ is smooth  and 
\bea
Ricci(\omega_{CY}) =0.
\eea
\end{theorem}

\subsection{Positivity notions for cohomology classes}

Another extension of the theory is the existence of K\"ahler-Einstein metrics, which are then necessarily singular, on manifolds $X$ whose first Chern class $c_1(K_X)$ is neither zero nor positive definite.

\medskip
To discuss the classes which are allowed, we recall briefly the definitions of some 
basic cones in the space of cohomology classes.
They were introduced by Demailly \cite{D1}
and play an important
role in his differential geometric approach to positivity problems in algebraic geometry.
Let $X$ be a compact \K manifold
and $\al\in H^{1,1}(X,\R)$ be a cohomology class. Then
\bea
\al\ \in \ 
{ 
\{\t: \hbox{closed $(1,1)$ forms}\}
\over
\{\t: \hbox{exact $(1,1)$ forms}\}
}
=
{ 
\{T: \hbox{closed $(1,1)$ currents}\}
\over
\{T: \hbox{exact $(1,1)$ currents}\}
}
\eea
We say $\al$ is pseudo-effective (psef) if there is a
closed $(1,1)$ current $T\in\al$ such that $T\geq 0$. We
say $\al$ is big if there exists $T\in\al$ with 
$T\geq \e \o$ for some $\e>0$ and some \K form~$\o$.
\v
Let ${\rm PSEF}(X)$ be the set of psef classes and ${\rm BIG}(X)$
the set of big classes. Then we clearly have ${\rm BIG}(X)\sub {\rm PSEF}(X)$.
Moreover, ${\rm PSEF}(X)$ is a closed convex cone in the vector space
$H^{1,1}(X,\R)$ and ${\rm BIG}(X)$ is an open convex cone.
If $T$ is psef then $T+\e\o$ is big for all $\e>0$. This shows
that ${\rm BIG}(X)$ is precisely the interior of ${\rm PSEF}(X)$.
\v
Let ${\rm KAH}(X)$ be the set of \K classes in $H^{1,1}(X,\R)$. Thus
$\al\in {\rm KAH}(X)$ if and only if there exists a \K form $\o\in\al$.
Thus ${\rm KAH}(X)$ is an open cone and we 
clearly have ${\rm KAH}(X)\sub {\rm BIG}(X)$.
This inclusion  may be proper. Let ${\rm NEF}(X)$ be the closure of
${\rm KAH}(X)$. An element of $\al\in {\rm NEF}(X)$ is called a
nef class. In summary,
\bea
\matrix{
{\rm BIG}(X)& \sub & {\rm PSEF}(X)\cr
\cup & & \cup\cr
{\rm KAH}(X)&\sub & {\rm NEF}(X)\cr
}
\eea
The cones on the left are open and those on the right are their closures.
\v
Let
\bea
\cT(X,\al)=\{ T\in \al: T \
\hbox{a closed $(1,1)$ current}, \ T\geq 0\}
\eea
Then $\al$ is pseudo-effective (psef) iff $\cT(X,\al)\not=\emptyset$.
We endow $\cT(X,\al)$ with the weak topology, so that
$T_j\rightharpoonup T$ iff 
$\I_X T_j\wedge \eta\ra \I_X T\wedge \eta$ for all
smooth $(n-1,n-1)$ forms $\eta$. The space
$\cT(X,\al)$ is compact in the weak topology.

\v
Fix a smooth volume form $dV$ on $X$ and let $L^1(X)=L^1(X,dV)$.
For $\t\in\al$ smooth, let
$\PSH^1(X,\t)= \{\phi\in L^1(X): \t+{i\over 2}\ddb\phi\geq 0\}$
endowed with the $L^1(X)$ topology.
The map $\phi\mapsto \t+{i\over 2}\ddb\phi$ defines 
$\PSH^1(X,\t)/\R\ \ra \ \cT(X,[\t])$
a homeomorphism of compact topological spaces.
The map
$ \sup: {PSH^1}(X,\t)\ \ra \ \R$
is continuous (this is Hartogs' lemma). Thus we have a homeomorphism
\bea
\{\phi\in {PSH^1}(X,\t): \sup \phi = 0\} 
\ra \cT(X,[\t]).
\eea

Let $X$ be a compact \K manifold and $\al\in H^{1,1}(X,\R)$ a big class.
Fix $\t\in \al$, a closed smooth $(1,1)$ form. Then, by definition, there
exists $\varphi\in PSH^1(X,\t)$ such that $\t+{i\over 2}\ddb\varphi\geq \e\o$ for
some \K metric $\o$ and some $\e>0$. Demailly's theorem says that
we may choose $\phi$ such that $\phi$ has analytic singularities. This
means that locally on $X$, 
\bea
\varphi\ = \ c\,\log(\s_{j=1}^N |f_j|^2)\ + \ \psi
\eea
where $c>0$, $f_j$ are holomorphic and $\psi$ is smooth. In particular,
the set where $T$ is smooth is a Zariski open subset of $X$. 
\v

Let $X$ and $\t$ be as above. Thus $\t$ is big, but it general it will not be positive. We  define $V_\t\in \PSH(X,\t)$, the extremal function of $\t$
(the analogue of the ``convex hull") by 
\bea
V_\t(x)\ = \ \sup\{\phi(x): \phi\in \PSH^1(X,\t),\,\sup_X\phi\leq 0\,\}\ 
\eea
Thus $V_\t=0$ if $\t$ is a \K metric.

\smallskip
The extremal function $V_\t$ has a number of nice properties. To describe them, and because $V_\t$ is not bounded in general, we need to extend the definition of Monge-Amp\`ere measures. We shall use the definition that does not charge pluripolar sets, and which can be described as follows.

\smallskip

Let $T_1,...,T_p$ be closed positive $(1,1)$ currents.  
For any $z\in M$ there there exists an open set $U$
containing $z$, and pluri-subharmonic functions $u_1,...,u_j$, for which
$T_j= {i\over 2}\ddb u_j$. Let $U_k = \cap_j \{u_j>-k\}\sub U$. Then the
non-pluripolar product $\langle T_1\wedge\cdots\wedge T_p\rangle$
is the closed $(p,p)$-current defined by \cite{BEGZ}
\bea
T_1\wedge\cdots\wedge T_p|_{U}\ = \ \lim_{k\to\i}{\bf 1}_{U_k}\bigwedge_{j=1}^p
{i\over 2}\ddb\max(u_j,-k)
\eea
where ${\bf 1}_{U_k}$ denotes the characteristic function of $U_k$.
The non-pluripolar product coincides with the Bedford-Taylor definition of
$T_1\wedge\cdots\wedge T_p$ if the potentials are all bounded.
We still denote $T_1\wedge\cdots\wedge T$ by $T^n$.
We can now describe the properties of extremal functions:

\begin{theorem} Let $X$ be a compact \K manifold and $\t$ a big
$(1,1)$ form. Let $V_\t$ be the extremal function of $\t$. Then

{\rm (1)} $V_\t\in \PSH^1(X,\t)$

{\rm (2)}  $V_\t$ has minimal singularities: if $\phi \in \PSH^1(X,\t)$ then
$\phi\leq V_\t +C$ for some $C\geq 0$.

{\rm (3)} On the set where $\t$ is smooth,
$V_\t$ is continuous, and ${i\over 2}\ddb V_\t$ is locally
bounded.

{\rm (4)} $({i\over 2}\ddb V_\t)^n$ has $L^\i$ density with respect to $dV$. In particular,
if $\varphi\in\PSH^1(X,\t)$ then $\I_X|\varphi|\,({i\over 2}\ddb V_\t)^n<\i$.

{\rm (5)} $V_\t$ is maximal with respect to $\t$, that is
\bea
({i\over 2}\ddb V_\t)^n\ = \ {\bf {1}}_{\{V_t=0\}}\t^n
\eea
\end{theorem}

\subsection{Prescribing the Monge-Amp\`ere measure}

If $[\t]$ is a big class and $T_1,...,T_n\in \cT(X,[\t])$ then
$\I_X\langle T_1\wedge\cdots\wedge T_n \rangle \leq \I_X[\t]^n$. If we have
$T_1=\cdots T_n=T=\t+{i\over 2}\ddb\varphi$ and if equality holds, we say  $T$ has full mass. Define
\bea
&&
\cT^0(X,\t)\ = \ \{T\in T(X,\t) : T\ \ {\rm has\ full\ mass}\ \}
\nonumber\\
&&
\cT^1(X,\t)\ = \ \{T=\t+{i\over 2}\ddb\varphi\in T^1(X,\o): \I_X |\varphi|\, 
(\t+{i\over 2}\ddb\varphi)^n < \i\ \}
\eea

Let $\cM_X$ be the space of probability measures on $X$, let $\cM_X^0$ consist of
those measures which take no mass on pluri-polar sets, and let $\cM_X^1$ consist
of those measures of ``finite energy" (to be defined below). Then 
it has been shown by Guedj and Zeriahi \cite{GZ}
and by Berman, Boucksom, Guedj, and Zeriahi \cite{BBGZ} respectively
that the map
$T \mapsto T^n$ defines bijections:
\bea
\cT^0(X,\o) \ra \cM_X^0, \ \ \cT^1(X,\o) \ra \cM_X^1.
\eea

\subsection{Singular KE metrics on manifolds of general type}

Before stating the results, it is convenient to recast the K\"ahler-Einstein equation
for negative curvature in a slightly different form from usual. Let $X$ be a compact K\"ahler manifold, and let $K_X$ be its canonical bundle.
If $c_1(K_X)$ is a \K class, then the Aubin-Yau theorem
says that for any \K metric $\o$ there is a unique smooth 
$\psi\in \PSH(X,-\Ricci(\o))$ such that
\be
\label{KE}(-\Ricci(\o)+{i\over 2}\ddb\psi)^n\ =\ e^\psi\o^n
\ee
To put this in the usual form, let $\eta = -\Ricci(\o)+{i\over 2}\ddb\psi$. 
Then $\eta>0$ and 
\bea
\Ricci(\o)-\Ricci(\eta)={i\over 2}\ddb\psi = \Ricci(\o)+\eta
\eea
which implies $\Ricci(\eta)= -\eta$.

\medskip
Now assume that $c_1(K)$ is big and nef. Then Tsuji \cite{Ts}  proved that
there is a subvariety $Z\sub X$ and a smooth function $\psi\in C^\i(X\backslash Z)$
such that
$\eta = -\Ricci(\o)+{i\over 2}\ddb\psi>0$ and 
such that (\ref{KE}) holds on $X\backslash Z$. Thus $\Ricci(\eta)=-\eta$
on $X\backslash Z$.
\v
 Tian-Zhang \cite{TZ} proved that $\psi$ extends to a locally bounded
$\psi\in \PSH(X,-\Ricci(\o))$ 
satisfying $\I (-\Ricci(\o)+{i\over 2}\ddb\psi)^n=\I (-\Ricci(\o))^n$ (i.e., $\psi$
 has full MA measure) and that (\ref{KE}) holds on all of $X$.
\v
Now assume that $c_1(X)$ is big. Then 
[EGZ] showed  that there is a unique $\psi\in \PSH^1(-\Ricci(\o))$ 
of full Monge-Amp\`ere measure  such that (\ref{KE}) holds.
The [EGZ] proof uses the existence of  canonical models.
Tsuji  described an interesting approach to proving the existence of a singular
K\"ahler-Einstein metric without resorting to the existence of a canonical model in
\cite{Ts}.
Then Song-Tian \cite{ST09} gave an independent proof, via the K\"ahler-Ricci flow.
A  new proof was also given by [BEGZ] which used a
generalized comparison principle. More recently, a proof using variational methods
was given in [BBGZ].

\section{Variational Methods for Big Cohomology Classes}
\setcounter{equation}{0}

A basic property of the Monge-Amp\`ere determinant is that it can be interpreted as the variational derivative of a concave energy functional. In fact, if $\o$ is a smooth K\"ahler form on a compact complex manifold $X$, and we set
\be\label{E}
 E_\o(\varphi)\ = \ {1\over n+1}\s_{j=0}^n \I_X \varphi(\o+{i\over 2}\ddb\varphi)^j\wedge
\o^{n-j}
\ee
for $\o_\vp\equiv\o+{i\over 2}\ddb\varphi>0$. Then for smooth and small
variations $\delta\vp$, we have
\be\label{derivative} 
\delta E= \I_X \delta\vp\,\o_\varphi^n\ \ {\rm and}\ \ 
\delta^2E=-n\I_X d\delta(\vp)\wedge d^c\delta\vp\wedge \o_\varphi^{n-1}
\ee
This shows that ${\delta E\over\delta \vp}=\o_\vp^n$ and that $E$ is concave.

The functional $E(\vp)$ is actually equal to $E(\vp)=-(J(\vp)-\int_X\vp\o^n)$, where
\bea
J(\vp)=\sum_{j=0}^{n-1}{n-j\over n+1}\int_X i\p\vp\wedge\bar\p\vp\wedge \o_\vp^{n-1-j}\wedge \o^j.
\eea
This relation will play an important role below
The functional $E$ is sometimes denoted by $-(n+1)F^0$ in the literature,
where $F^0$ is the Aubin-Yau energy functional. 

\smallskip

The goal of this section is to describe some recent work of Berman, Boucksom, Guedj, and Zeriahi \cite{BBGZ} taking this variational viewpoint further. Their approach, which
works in the generality of   big cohomology classes,  allows one to use direct methods of the calculus of
variations to obtain solutions to a variety of Monge-Amp\`ere equations.

\subsection{Finite dimensional motivation}

We start with a finite dimensional model for the method introduced in \cite{BBGZ}.
Let $P\sub\R^n$ be a  convex domain and $E:P\ra\R$ be a strictly concave smooth function. Then $\na E: TP\ra \R$ where $TP=P\times\R^n$ is the tangent bundle of $P$. In other words,
$\na E: P\ \ra \ (\R^n)^*$
where  $(\R^n)^*$ is the dual of $\R^n$.
Since $E$ is strictly concave, $\na^2E(\vp):\R^n\ra(\R^n)^*$ is strictly
negative
definite for all $\phi\in P$ and $\na E$ is a diffeomorphism of $P$ onto a  convex domain $\cM\sub(\R^n)^*$.
\v
More generally, suppose $E$ is concave (but not strictly concave) and
that there is a flat direction, that is, an element $\vp_0\in P$ such that 
$E(\vp+t\vp_0)=E(\vp)+t$ for all $t\in\R$ and let $\cT^1=P/(\R\cdot \vp_0)$
which we view as a domain in the vector space $W=\R^n/(\R\cdot\vp_0)$. Then
$\na E:P\ra (\R^n)^*$ is invariant under $\R\cdot\vp_0$ so
\bea
\na E: \cT^1 \ra \cM^0
\eea
 where $\cM^0\sub (\R^n)^*$ is the $n-1$ dimensional affine space
\bea
\cM^0\ =\ \ \{\m\in\R^n: \langle\vp_0, \m\rangle=1\}\ = \ 
\hbox{the elements of $(\R^n)^*$ of {\it full mass}}
\eea
Now we want to impose the following condition: $E$ is strictly concave
in every direction other than $\R\cdot\vp_0$. There are several equivalent
ways of making this condition precise:

\begin{lemma} Let $E:P\ra\R$ be a concave function and
assume that $E(\vp+t\vp_0)=E(\phi)+t$ for all $t\in\R$.
Then the following conditions are equivalent

{\rm (1)} The map $t\mapsto E((1-t)\psi_0+t\psi_1)$ is strictly concave 
whenever 
$\psi_1-\psi_0\notin \R\cdot\phi_0$.

{\rm (2)} The negative definite map
$\na^2E(\vp): W\ra W^*$  is strictly negative definite if $\vp\in \cT$.

{\rm (3)} The map $\na E: \cT^1\ra\cM$
 is a diffeomorphism where $\cM\sub \cM^0$ is the image of $\na E$.
 
{\rm (4)} Fix $\psi_0\in \cT^1$, let $\m_0=\na E(\psi_0)$, $L_0(\vp)=\vp\cdot\m_0$ and define
 $J_0:\cT^1\ra [0,\i)$ by
 \be\label{jzero}
 J_{\psi_0}(\vp)\ \ =  \ E(\psi_0) + \na E(\vp_0)\cdot (\vp-\psi_0)-E(\vp)\ = \ 
L_0(\vp)-E(\vp)+ C
 \ee
 Then $J_0= J_{\psi_0}:\cT^1\ra [0,\i)$ is strictly convex.
 \end{lemma}

The proof of the lemma is easy and will be omitted. Henceforth,
we shall assume that $E$ satisfies any one of the equivalent conditions
enumerated in Lemma 1.
\v
The function $J_0$ is proper so
the sets $\{J_0\leq C\}\sub \cT^1$ are compact and exhaust $\cT^1$. In
the case $\cT^1$ is bounded, we can define $\cT$ to be the closure of $\cT^1$
and extend $J_0: \cT\ra [0,\i]$ as a continuous map between compact spaces.

   \v
For $\m\in\cM^0$ we let 
\bea
E^*(\m)=\sup_{\vp\in \cT} (E(\vp)-\vp\cdot\m)\ = \ 
\sup_{\vp\in\cT}F_\m(\vp)\ \in (-\i,\i]
\eea
We let $\cM^1\sub\cM^0$ be the elements with {\it finite energy}:
$$ \cM^1\ = \ \{\m\in\cM^0: E^*(\m)<\i\}
$$
If $\m\in\cM$ then, by definition, there exists $\varphi\in\cT^1$
such that $\na E(\vp)=\m$ and hence $\vp$ is a critical
point of $F_\m$. Since $F_\m$ is strictly concave we conclude
that $F_\m$ achieves its maximum at $\vp$ and thus
\bea
\cM\sub \cM^1\sub\cM^0
\eea
Conversely, if $\m\in\cM^1$ then $\sup F_\m<\i$ and $\m\in\cM$ if and only if
$F_\m$ {\it achieves} its sup (i.e. there exists $\vp\in\cT^1$ such that
$F_\m(\vp)=E^*(\m)$).
\v
The inclusion $\cM\sub\cM^1$ may be strict.
We give two simple examples:
\v
Example A. Suppose $\cT=\R$ and $E:\R\ra\R$ is smooth and concave,
and $E(x) = \log x-x$ for $x$ large. Then $\cM=\cM^1=(-1,1)$. So for 
example A, 
\be\label{23} \cM=\cM^1
\ee

\v
Example B. Suppose $\cT=\R$ and $E:\R\ra\R$ is smooth and convex,
and $E(x) = -x-{1\over x}$ for $x$ large. Then $\cM=(-1,1)$
but
$\cM^1= [-1,1]$. 
Thus in example B property (\ref{23}) fails.

\v
There is a simple criterion for guaranteeing that (\ref{23}) holds.  Suppose
$J:\cT\ra [0,\i)$ is a proper function and let
$\m\in\cM^1$. We say that $F_\m$ is coercive with
respect to $J$ if there exists $\e,C>0$,
 depending on $\m$, such that  
 \be\label{coerce}L_\m-E\ =\ -F_\m\geq \ \e J-C\ = \ \e(L_{\m_0}-E)-C'
 \ee
 We say $F$ is $J$-proper
 if $-F_\m\ra \i$ as $J\ra\i$. Clearly $J$-coercive implies $J$-proper.
Note that (\ref{coerce}) holds trivially if $\e=0$: If $\m_0=\na E(\psi_0)$
 then 
 $J_0(\phi)=F_{\m_0}(\psi_0)-F_{\m_0}(\phi) \geq 0$.

\begin{lemma}\label{lemmaproper}  
Let $\m\in\cM^1$ and assume that $F_\m$ is $J$-proper for some $J$.
Then  $\m\in\cM$. 
\end{lemma}
To see this,  we
fix $\vp_0\in \cT$ and let $A=J(\vp_0)+1$.
Then the set $\cT_A=\{J\leq A\}$ is compact so there exists $\hat \vp$ such
that $F_\m(\hat \vp)\geq F_\m(\vp)$ for all $\vp\in\cT_A$ and
hence for all $\vp\in \cT$. Thus $F_\m$ achieves its sup at $\hat \vp$
which implies $\na F_\m(\hat\vp)=0$, that is, $\na E(\hat\vp)=\m$.
This shows $\m\in\cM$.
\v
In order to apply Lemma \ref{lemmaproper} we need to find an
appropriate $J$. It turns out that if $F_\m$ is $J$-coercive for some
$J$, then it is coercive for the function $J_0$ constructed
in (\ref{jzero}):

\begin{lemma} Let $\m\in\cM^1$. Then the following are equivalent

{\rm (1)} $\m\in\cM$

{\rm (2)} $F_\m$ is $J_0$-coercive.

{\rm (3)} $F_\m$ is $J$-coercive for some exhaustion function $J$.

{\rm (4)} $F_\m$ is $J$-proper for some exhaustion function $J$.
\end{lemma}
\noindent {\it Proof.} We need only show that 1) implies 2). Let $\m\in\cM$.
We must show that for some $\e,C>0$ the following holds.
\be
-F_\m(\vp) \geq \e(F_{\m_0}(\psi_0)-F_{\m_0}(\vp))-C
\ee
which we rewrite as
\be
\vp\cdot\m_\e-E(\vp)=\vp\cdot {(\m-\e\m_0)\over 1-\e}-E(\vp)\geq {\e F_{\m_0}(\psi_0)-C\over 1-\e}
\ee
But $\cM$ is an open set so for $\e$ sufficiently small, $\m_\e\in\cM$
which implies $ \vp\cdot\m_\e-E(\phi)$ is bounded below.
\v
Now fix $\psi_0\in\cT$ and define $J_0$ as in (\ref{jzero}).
\begin{lemma} $\cM=\cM^1$ if and only if $F_\m$ is $J_0$-coercive
for all $\m\in\cM^1$.
\end{lemma}
Now we record a condition on $J_0$ that guarantees $F_\m$ 
is coercive for all $\m$. This key condition will hold in the
infinite dimensional setting and is used in the proof of
coerciveness.

\begin{lemma}\label{quadratic} Suppose that $J_0$ grows quadratically or, more
generally, suppose that
\be\label{pgrowth} |J_0(t\phi)|\ \leq \ C(t^p|J(\phi)|+1)\ 
\hbox{if $t\leq 1$}
\ee
for some $p>1$. Then $F_\m$ is $J_0$-coercive for all $\m\in\cM^1$.
\end{lemma}
Observe that we can rewrite the coercive condition (\ref{coerce})
as follows:
\be\label{coerce1} L_{\m_0}-L_\m\leq (1-\e)J_0+C_1
\ee
Now (\ref{pgrowth}) implies that for arbitrary $\phi\in\cT$, we have
$\phi_p={\phi\over |J(\phi)|^{1/p}}\in \{J_0\leq C_2\}$ which is a compact
set. This implies
\be
{|L_{\m_0}-L_\m|(\vp)\over |J_0(\phi)|^{1/p}}= |L_{\m_0}-L_\m|(\vp_p)
\leq C_3
\ee
which clearly implies (\ref{coerce1}).

\subsection{The infinite dimensional setting}

We describe the work of \cite{BBGZ}.
Let $\t$ be a big $(1,1)$ form on a compact \K manifold $X$.
Let $P=\PSH^1(X,\t)$ and $\cT=\cT(X,\t)$ and $\cM'$ the space of
positive metrics on $X$. Our goal is to define a concave
function $E: \PSH^1(X,\t)\ra\R$ with the property: $\na E: P\ra \cM'$
is the map $\vp\mapsto MA(\phi)$, where $MA(\vp)$ is the Monge-Amp\`ere
measure, defined by the non-pluripolar product.
Moreover, we will show that if $\vp_0=1$, the constant function,
then $\vp_0$ is a flat direction for $E$ and $\cT$ is a
complement to $\R\vp_0$ in $P$. Then, as in the finite dimensional
case, we obtain a map $\na E:\cT\ra\cM'$. Also, in analogy with the finite
dimensional case, we define  the function $F_\m(\vp)=
E(\vp)-\I_X\vp\,d\m$ for $\m\in\cM'$. Our goal is to prove (\ref{23})
and this will be done by establishing (\ref{coerce}) for a
suitably chosen $J$.

\v
The functional $E(\vp)$ has been defined in (\ref{E}) for $\t=\o$ a K\"ahler form,
and $\vp$ a smooth form in $\PSH^1(X,\t)$. To extend it to $\t$ a big form and
$\vp\in \PSH^1(X,\t)$, we proceed as follows:
Define $E(\vp)=\inf\,\{E(\psi): \psi\in \PSH^1(X,\o)\cap C^\i(X)$
and $\psi\geq\vp\}$. 

\v
Now let $\t$ be a big $(1,1)$ form (not necessarily K\"ahler). Let
$\vp\in \PSH^1(X,\t)$ and assume $\vp$ has minimal singularities
(i.e. that $\phi-V_\t$ is bounded). Define
\be
\label{EE}
E(\vp)\ = \ {1\over n+1}\s_{j=0}^n\I_X (\vp-V_\t)(\t+{i\over 2}\ddb\vp)^j
(\t+{i\over 2}\ddb V_\t)^{n-j}
\ee
Of course this coincides with (\ref{E}) in the case where $\t=\o$
is \K (since in that case, $V_\t=0$). If $\vp$ is arbitrary,
then again  extend using the monotonicity property of $E$
as before: $E(\vp)=\inf\{E(\psi): \psi\geq \vp$ and $\psi$ has
minimal singularities$\}$.
\v
We wish to implement the finite dimensional program in the infinite
dimensional setting.
Since $E$ may assume the value $-\i$, we must, at the outset, restrict
$E$ to the set ${\cal E}^1=\{\vp\in \PSH^1(X,\t): |E(\vp)|<\i\}$ so
that $E: {\cal E}^1\ra \R$ is concave. If $\vp_0=1$, then $\vp_0$
is a flat direction so, as in the finite dimensional program, we let $\cT^1={\cal E}^1/\R$
Again, we let $F_\m(\vp)=E(\vp)-\I_X(\vp-V_\t)\,d\m = E(\vp)-L_\m(\vp)$
and let $\cM^1\sub\cM^0$ be those measures with finite energy.

\subsection{Statement of theorems and sketch of proofs}

\begin{theorem}\label{bbgz1}(\cite{BBGZ})\label{fe} Let $\t$ be a big $(1,1)$ form on a compact complex manifold $X$.
If $\vp\in{\cal T}^1(X,\t)$ then $\na E(\vp)=MA(\vp)$ has full mass and
finite energy, that is, $MA(\vp)\in \cM^1$.
Conversely, if $\m\in \cM^1$ then there exists a unique $\vp\in\cT^1$
such that $MA(\vp)=\m$. Moreover, the solution $\phi$ satisfies the
following bound:

\be\label{enbound} {1\over n}E^*(\m)\leq J_0(\vp)\ \leq nE^*(\m)
\ee
where $J_0=J_{\psi_0}$ is defined in (\ref{jzero}).
\end{theorem}

We sketch the proof of the converse, using the finite dimensional program as our guide.
Recall that
\bea\label{jinf} 
J_\psi(\vp)&=&E(\psi)-E(\vp)+\I_X(\vp-\psi)MA(\psi)=\nonumber\\
&=& \sum_{j=0}^{n-1}{j+1\over n+1}\I_X \pl(\vp-\psi)\wedge\bar\pl (\vp-\psi)
\wedge\t_\psi^j\wedge\t_\vp^{n-1-j}\\ \nonumber
\eea

The first equality is a definition and the second follows via integration by parts.
In particular, we see
\be\label{jbound} n^{-1}J_\psi(\vp)\leq J_\vp(\psi)\leq nJ_\psi(\vp)
\ee
Now we can establish the bound (\ref{enbound}):  Assume $\vp\in {\cal E}^1$
and that $MA(\vp)=T^n=\m$, where $T=\t+{i\over 2}\ddb\vp$. Then (\ref{jinf})
and (\ref{jbound}) imply
$$ 
{1\over n}J_{0}(\vp)\leq J_\vp(V_0)=E(\vp)+\I_X(V_0-\vp)MA(\vp)=
E(\vp)-\I_X (\vp-V_0)\,d\m=F_\m(\vp)
$$
Since $MA(\vp)=\m$ one shows, in analogy with the finite dimensional picture, that $E^*(\m)=\sup_{\psi\in {\cal E}^1}F_\m(\psi)=F_\m(\vp)$. This proves (\ref{enbound}).

\v
Now let us fix $\m\in\cM^1$. We wish to prove the existence of $T\in\cT^1$
such that $T^n=\m$. The first step is to prove that $F_\m$ is $J_0$ coercive,
where $J_0$ is defined with respect to the potential $\psi_0=V_0$, as in
(\ref{jzero}). To do this, we wish to use Lemma \ref{quadratic}, which means
that we must prove that $E$ grows quadratically. But this follows easily
from the definitions (here we assume $V_\t=0$ for simplicity): If $\sup_X\vp=0$ then
\bea
 |E(t\vp)|={1\over n+1}\left|
\s_{j=0}^{n}\I (t\vp)(t(\o+{i\over 2}\ddb\vp)+ (1-t)\o)^j\o^{n-j}
\right| \ \leq c(n)[t^2E(\vp)+1]
\eea
For $C>0$ let ${\cal E}_C$ be the compact set ${\cal E}_C=\{\vp\in {\cal E}^1: E(\vp)\geq C\}$. 
Since $F_\m$ is coercive, there exists $C>0$ such that 
$\sup_{\phi\in{\cal E}^1} F_\m =\sup_{\phi\in{\cal E_C}} F_\m$. 
\v
The first main difficulty in implementing the finite dimensional program is proving that $F_\m$ is upper semi continuous on ${\cal E}_C$.
To do this, \cite{BBGZ} restrict first to the case where $\m\in \cC=\{\m\in\cM^1: \m\leq A\cdot Cap$ for some $A>0\}$ (here $Cap$ is the pluri-subharmonic capacity).
\v
Fix $\m\in C$. Since $E$ is easily seen to be usc, it suffices to show that $L_\m$
is continuous. Let $T: K \ra L^1(X,d\m)$ be the map $T(\vp)=\vp-V_\t$ (here
$K\sub{\cal E}^1$ is any compact convex subset). It's not hard to show  that $T(K)$ is closed and that it has a closed graph. If we could prove  $T(K)$ is compact,
then we would be done. To see this, let $\vp_j\ra \vp$ in $K$. If $T(K)$ is
compact then $T(\vp_j)\ra f$ (after passing to a subsequence). But the closed
graph property implies $f=\vp$. Thus $T(\vp_j)\ra T(\vp)$ so $T$ is continuous
which implies $L_\m$ is continuous.
\v
Instead of proving that $T(K)$ is compact, we prove something weaker, namely
that $T(K)\sub L^1(X,d\m)$ is ``convex combination compact". Recall that
if $B$ is a Banach space and if $T\sub B$ is closed and convex, then $T$ is
convex combination compact if for every sequence $\tau_1, \tau_2,...\in T$ there exists $\tau_1',\tau_2',...\in T$ such that $\tau_j'$ is a finite convex combination
of $\tau_j,\tau_{j+1},...$ and such that $\tau_j'$ converges. 
Observe that showing $T(K)$ is convex combination compact suffices for our purposes: Let $\vp_j\ra \vp\in K$.
Then one shows $L_\m(\vp_j)$ is bounded and hence 
$L_\m(\vp_j)=\I T(\vp_j) \ \ra \ \ell\in\R$ (after passing to a subsequence).
On the other hand, the convex combination compactness of $T$ implies there exist $\psi_1,\psi_2,...\in K$ 
such that $\psi_j$ is a convex combination of $\vp_j,\vp_{j+1}...$ and
$T(\psi_j)\ra f$ for some $f\in T(K)$. Since $\psi_j\ra \vp$ we see that $f=T(\vp)$
by the closed graph property. Thus
\be
\I T(\psi_j)\,d\m \ra \I f\,d\m \ = \ \I T(\vp)\, d\m
\ee
On the other hand, $\lim_j \I T(\psi_j)\,d\m=\lim_j \I T(\vp_j)\,d\m$ so
$\I T(\vp_j)\ra \I T(\vp)$. This shows $L_\m$ is continuous on $K$.

\v
To show that $T(K)$ is convex combination compact, it suffices to prove $T(K)$ is weakly compact (by
the Hahn-Banach Theorem). On the other hand, the Dunford-Pettis theorem says that
 to show that $T(K)$ is weakly
compact, it suffices to show that $T(K)$ is uniformly integrable, that is,
there exists $\e_k\downarrow 0$ such that
\be
\left|
\I_{T(\vp)<-k} T(\vp)\,d\m
\right|\ \leq \ \e_k \ \ \hbox{for all $\vp\in K$}
\ee
But $\m\in \cC$ implies
\be
\I_0^\i t\m\{t< V-\vp\}\,dt\ \leq \ A\I_0^\i t Cap\{t<V-\vp\}\,dt\ \leq \ C_1
\ee
Thus $T(K)$ is a bounded subset of $L^2(X,d\m)$ and hence, by H\"older's inequality,
$T(K)$ is uniformly integrable. 

\v
The second main difficulty is to prove the Euler-Lagrange equation is satisfied
at the point where $F_\m$ reaches its maximum. The problem is that the minimum
$\vp$ may only satisfy $\o+{i\over 2}\ddb\vp\geq 0$, and a variation
$\vp_t=\vp+\delta\vp$ may no longer satisfy this condition even if $\delta\vp$ is smooth and small. This
is overcome using a technique of \cite{BB} which we now describe. Define for each $u$ an upper semi-continuous function on $X$, its $\t$-psh envelope,
\be
P(u)\ = \ \sup\{\vp\in \PSH^1(X,\t): \vp\leq u \hbox{\ on $X\}$}
\ee
Then $P(u)\in \PSH^1(X,\t)$ so the function
$g(t)= E(P(\vp+t\delta\vp))-L_\m(\vp+t\delta\vp)$ is now 
well-defined for $t\in\R$, and has a maximum at $t=0$. It is shown in \cite{BBGZ} 
that we still have
$0=g'(0)= \I_X \delta\vp (MA(\vp)-d\m)$.
Thus the equation $MA(\vp)=\m$ has a solution
in the case where $\m\in\cC$.
\v
We now remove the assumption $\m\in \cal C$. 
Thus we let $\m$ be a non-pluripolar probability measure with finite energy. Our
goal is to prove that there exists $\vp\in \cT^1$ such that $MA(\vp)=\m$.
To do this, we make use of the following
lemma of Cegrell \cite{Ceg}:
\begin{lemma}\label{ceg} Let $\m$ be a non-pluripolar probability measure. Then there
exists $\n\in\cC$ and $f\in L^1(X,\n)$ such that $\m=f\n$.
\end{lemma}

Write $\m=f\n$ as in the lemma. For  $k>0$ choose $\e_k\geq 0$ so
that $\m_k=(1+\e_k)\min(f,k)\n$ is a probability measure. Then $\m_k\leq 2k\,Cap$
so $\m\in\cC$. Thus, by what has been proved thus far, $\m_k=T_k^n$
for some $T_k\in \cT^1$. Next we observe that $\m_k\leq 2\m$. Using the fact that $E^*(\m)<\i$, one
shows that $E^*(\m_k)\leq C<\i$ for some $C>0$.
The bound (\ref{enbound}) then implies that the
$T_k$ all lie in the compact set $\{J_0\leq nC\}$. Thus, after passing
to a subsequence, we conclude $T_k\ra T$ for some $T\in \cT^1$. Since
$T^n_k=\m_k$ we can take the limit as $k\ra\i$ to conclude $T^n=\m$
(this follows by Fatou's lemma in the case $n=1$, and a generalization
of Fatou's lemma, due to \cite{BEGZ}, in the case $n>1$).
\v
We remark that if $\m$ is a non-pluripolar probability measure (not necessarily
of finite energy), then one can still apply Lemma \ref{ceg} to conclude
$\m=f\n$ and one can still construct $\m_k=(1+\e_k)\min(f,k)\n$  as above.
Then Theorem \ref{bbgz1} implies the existence of $T_k\in \cT^1$ such
that $T_k^n=\m_k$. Since $\m$ is not assumed to have finite energy,
we cannot conclude that the $T_k$ lie in a compact subset of $\cT^1$. 
On the other hand, we have $T_k\sub \cT^1\sub \cT$ and $\cT$ is compact.
Thus, after passing to a subsequence, $T_k\ra T$ for some $T\in \cT$. 
Then \cite{BBGZ} show that one can again take the limit as $k\ra \i$ to
conclude $T^n=\m$. This gives a variational proof of the following theorem
of \cite{BEGZ}:

\begin{theorem}\label{bbgz2} Let $\m$ be a non-pluripolar probability measure on $X$.
There there exists $T\in \cT(X,\t)$ such that $T^n=\m$.
\end{theorem}

We close by describing one more application of the variational method in \cite{BBGZ}.

\begin{theorem}\label{bbgz3} Let $(X,\o)$ be a manifold of general type and let
$\t=-\Ricci(\o)$. Then there exists
$\psi\in \PSH^1(X,\t)$ such that 
\be\label{general}(\t+{i\over 2}\ddb\psi)^n=e^\psi\o^n
\ee
\end{theorem}
To prove this theorem, we consider the functional $F = E- L$ where
$L(\vp)={1\over 2}\log\I_X e^\vp\o^n$ and follow the same three steps
as 
in the proof of Theorem \ref{bbgz1}. The first step is to show that $F$ is 
$J_0$ coercive. The second step 
is to show that  $L$ is continuous (which implies that $F$ is upper semi-continuous).
And the third step is to prove that the critical points for $F$ are solutions
to (\ref{general}).
\v
Steps one and three proceed exactly as in the proof of Theorem \ref{bbgz1}.
Thus we restrict ourselves to step two, which is much easier than the corresponding
step in Theorem \ref{bbgz1}. Indeed, if $\vp_j\ra \vp$ is a convergent sequence
in $\PSH^1(X,\t)$ then after passing to a subsequence, $\vp_j\ra \vp$ almost
everywhere. On the other hand, Hartogs' lemma implies that $\vp_j$ is bounded
above. Thus $\I_X e^{\vp_j}\o^n\ra \I_X e^{\vp}\o^n$ and hence $L$ is continuous.

\v
\noindent
We observe that Theorem \ref{bbgz3} follows from Theorem \ref{egz1}. To see this, let
$\pi:X\ra X_{\rm can}$ be the canonical model of $X$ and choose $\vp_{\rm can}$ to
be the solution of (\ref{98})

\begin{equation}\label{98a}
(\omega+ {i\over 2}\ddbar \varphi_{\rm can})^n = e^{\varphi_{\rm can}} \Omega.
\end{equation}
where ${i\over 2}\ddb\log\O=\o$. Choose $\Theta$ such that ${i\over 2}\ddb\log\Theta=\t$.
Let $\ti\O=\pi^*\O$ and $\ti\vp_{\rm can}=\vp_{\rm can}\circ\pi$. Define $\psi$
by the equation $e^{\it\psi_{\rm can}}\ti\O=e^\psi\Theta$.
Applying $\pi^*$ to both sides of (\ref{98a}) we obtain (\ref{general}).

\medskip
The variational method establishes the existence of generalized solutions to the complex Monge-Amp\`ere equation. It is then important to determine when the generalized solution is actually smooth. One such result is the theorem of
Szekelyhidi-Tosatti \cite{ST} which asserts the smoothness of the generalized 
solution when it is known to be bounded and the right hand side is smooth.

\smallskip
It may also be noteworthy that the above solutions to the Monge-Amp\`ere equation  can be obtained as limits of the critical points of certain naturally-defined finite-dimensional analogues of the infinite-dimensional functionals. In the Fano case, the proof makes use of the Moser-Trudinger inequality proved in \cite{PSSW}. Finite-dimensional approximations also play a major role in the construction of solutions to the homogeneous Monge-Amp\`ere equation in \S 13.

\section{Uniqueness of Solutions} 
\setcounter{equation}{0}

It is a remarkable fact that 
the map
$\cT^0(X,\t)\ra \cM_X^0$, defined
by $T\mapsto T^n$, is bijective.
The surjectivity is part of Theorem \ref{bbgz2}. The injectivity was proved in \cite{BEGZ},
by adapting the proof of Dinew \cite{Di09}, who proved the injectivity in the K\"ahler case. We give here a slightly streamlined version of the proof of Dinew. 
\begin{theorem}\label{Dinew}
Let $(X,\o)$ be a compact K\"ahler manifold and 
${\cal E}(X,\o)\sub \PSH^1(X,\o)$ denote the potentials $\vp$ such that 
$(\o+{i\over 2}\ddb \vp)^n$ has full mass.  Let $\vp,\psi\in {\cal E}(X,\o)$. Assume
$(\o+{i\over 2}\ddb\vp)^n=(\o+{i\over 2}\ddb\psi)^n$. Then $\vp-\psi$ is constant.
\end{theorem}

To prove the theorem, we wish to make use of the comparison principle. If
we apply it  directly to $V=\{\vp<\psi\}$ we get
$\I_V\o_\psi^n\leq \I_V\o_\vp^n$, which
is not useful (in fact, the inequality is an equality since $\o_\phi^n=\o_\psi^n$).
Instead, we  shall apply the comparison principle to the set 
$V_\e(\t,\r)=\{(1-\e)\vp+\e\t< (1-\e)\psi+\e\r\}$  
where $\t,\r$ are potentials to be chosen later.
We obtain, for
$T$ a positive closed current, and $k\geq 1$,
\be\label{comp1} \I_{V_\e} T\wedge((1-\e)\o_\psi+\e\o_\r)^k \ \leq
\I_{V_\e} T\wedge((1-\e)\o_\vp+\e\o_\t)^k
\ee
Assume $\I_{V_\e} T\wedge \o_\psi^k=\I_{V_\e} T\wedge \o_\vp^k$. Then the
leading terms cancel so
 \be\label{comp} \I_{V_\e} T\wedge\o_\psi^{k-1}\wedge\o_\r\ \leq \ 
\I_{V_\e} T\wedge\o_\vp^{k-1}\wedge\o_\t\ + \ O(\e)
\ee
We shall also need the following generalization
of Lemma \ref{kcomparison} proved in
Dinew \cite{Di09a}:
\begin{theorem}\label{est1} Let $\vp,\psi\in {\cal E}(X,\o)$ and $\m$ a positive non-pluripolar
measure. Assume that $\o_\vp^n\geq fd\m$ and $\o_\psi^n\geq gd\m$ for some 
non-negative $f,g\in L^1(d\m)$. Then for $0\leq k\leq n$
\be 
\o_\vp^k\wedge \o_\psi^{n-k}\ \geq \ f^{k\over n}g^{n-k\over n}d\m
\ee
In particular, if $\o_\vp^n=\o_\psi^n$, then $\o_\vp^n=\o_{t\vp+(1-t)\psi}^n $ for
all $t\in [0,1]$.
\end{theorem}

We return to the proof of Theorem \ref{Dinew} and
follow  the argument in  \cite{Di09}. Let $\m=\o_\vp^n$
and  define $f:\R\ra [0,1]$ by $f(t)=\m(\{\vp<\psi+t\}$. Then $f$ is left continuous.
Moreover, since $f$ is increasing, there is a countable set $\Sigma$ such that $f$
is continuous on 
$\R\backslash\Sigma$.
\v
The key step in the proof is to show that ${\rm Image}(f)=\{0,1\}$: Assume
not. Then there exists $\alpha\in\R\backslash\Sigma$ such that $0<f(\al)<1$.
To see this, let $\beta\in\R$ satisfy $0<f(\b)<1$. Now choose an increasing sequence
$t_j\in\R\backslash\Sigma$ such that $t_j\ra \b$. Then we can take $\al=t_j$ for
any sufficiently large $j$. After replacing $\vp$ by $\vp-\al$ we may assume
$\al=0$. 
\v
Since $0<f(0)<1$ we may choose $0<q<1$ such that $1-q<f(0)<q$. With this choice of
$q$ we see that $\m(\{\vp<\psi\})<q$ and $\m(\{\psi<\vp\})<q$ and $\m(\{\vp=\psi\})=0$
(the last equality follows from the fact that $\al\in\R\backslash\Sigma$). 
To get a contradiction, we consider the probability measure $\hat \m = g\m$
where $g={1\over q}$ on $\{\vp<\psi\}$ and $g=c\m$ on $\{\vp\geq\psi\}$ for an
appropriately chosen $c>0$. The theorem of Guedj-Zeriahi \cite{GZ} implies there
exists $\r\in \PSH^1(X,\o)$ with $\sup\r=0$ and $\o_\r^n=\hat\m$.
Then, setting $a=({1\over q})^{1/n}$, and $\t=0$, and $k=n$, Theorem \ref{est1} implies 
\be\label{est2}
 \o_\psi^{n-1}\wedge\o_\r\ \geq\ a\,\o_\psi^n \hbox{\ \ on the set $V_\e(\t,\r)
 \sub\{\vp<\psi\}$}
\ee
Substituting in (\ref{comp}), and taking the limit, $\e\downarrow 0$, we get
$V_\e\uparrow \{\vp<\psi\}$ so
\be
a\I_{\vp<\psi} \o_\vp^n\leq \I_{\vp<\psi}\o_\vp^{n-1}\wedge\o
\ee
If instead we take $k=1$ and $T=\o_\psi^{n-1}$  we obtain
$$
a\I_{\vp<\psi} \o_\psi^n\leq \I_{\vp<\psi}\o_\psi^{n-1}\wedge\o
$$
Interchanging  $\vp$ and $\psi$ in the second estimate, we get
$a\I_{\vp\geq\psi} \o_\vp^n\leq \I_{\vp\geq\psi}\o_\vp^{n-1}\wedge\o$, where
we make use of the fact that $\m\{\vp=\psi\}=0$. Adding this
to the previous inequality we conclude that $a\leq 1$, a contradiction.
\v
The next step is to show 
\be\label{higher} \I_{\\vp<\psi}\o_\vp^j\o_\psi^{k-1}\o^{l+1}\ \leq \ 
\I_{\vp<\psi}\o_\vp^j\o_\psi^{k}\o^{l}\ = \ 0
\ee
for all $j,k,l$ such that $j+k+l=n$. To see this, we use induction on $l$. Let $\r=0$ and $\t=\vp$ and 
$T=\o_\vp^j\o^l$. Then applying (\ref{comp}) we obtain,  for every $\d>0$, 
$$ \I_{\vp<(1-\e)\psi-\d} \o_\vp^j\o_\psi^{k-1}\o^{l+1}\ \leq \ 
\I_{\vp<(1-\e)\psi-\d}\o_\vp^j\o_\psi^{k}\o^{l}+O(\e)
$$
Taking the limit as $\e\downarrow 0$ we obtain (\ref{higher}) but with
$\{\vp<\psi\}$ replaced  by $\{\vp\leq\psi-\d\}$. Now take the limit
$\d\downarrow 0$ to obtain (\ref{higher}).
\v
Taking $l=n$ we obtain
$\I_{\vp<\psi}\o^n=0$. Similarly $\I_{\vp>\psi}\o^n=0$. Since these sets are
plurifine open, they must be empty. This proves the theorem.

\section{Semiclassical Solutions of Monge-Amp\`ere Equations}
\setcounter{equation}{0}

In this section we discuss a particular method for solving Monge-Amp\`ere equations, namely by semiclassical limits of Bergman kernels. The equation accessible this way is the homogeneous complex Monge-Amp\`ere equation. It is a very specific equation, but one which is of great interest in the problem of finding
K\"ahler metrics of constant scalar curvature, and which is intimately linked with the notion of stability in geometric invariant theory
(see \cite{D99} and \cite{PS08} for a survey). We shall see how semiclassical limits lead to generalized solutions of this equation. This may be noteworthy in itself from the viewpoint of PDE theory, as generalized solutions of partial differential equations usually arise rather from either variational or Perron methods.

\subsection{Geodesics in the space of K\"ahler potentials}

We begin with some geometric motivation. Let $(X,\o_0)$ be a compact K\"ahler manifold without boundary of dimension $n$.
Then the space ${\cal K}$ of K\"ahler potentials
\bea
{\cal K}=\{\vp\in C^\infty(X);\ \o_\vp\equiv\o_0+{i\over 2}\ddb\vp>0\}
\eea
is formally an infinite-dimensional Riemannian manifold with tangent space $T_\vp({\cal K})
=\{\delta\vp\in C^\infty(X)\}$ and metric
\bea
\|\delta\vp\|^2=\int_X |\delta\vp|^2\o_\vp^n.
\eea
The geodesic equation for ${\cal K}$ is the Euler-Lagrange equation for the energy
functional
\bea
{\cal E}=\int_0^T\int_X \dot\vp^2 \o_\vp^n dt
\eea
for paths $[0,T)\ni t\to \vp(\cdot,t)\in {\cal K}$. Under a variation $\delta\vp$ of this path, we have
\bea
\delta {\cal E}
=\int_0^T\int_X (2\dot\vp\,\delta\dot\vp+\dot\vp^2\Delta_\vp\delta\vp)\o_\vp^n dt,
\eea
where $\Delta_\vp$ is the Laplacian with respect to the K\"ahler form $\o_\vp$.
For variations $\delta\vp$ fixing the end points, integrating by parts gives
\bea
\delta {\cal E}
=
-2\int_0^T\int_X\delta\vp(\ddot\vp-|\na\dot\vp|_{\o_\vp}^2)\o_\vp^n dt.
\eea
and thus the geodesic equation for paths $[0,T)\ni t\to \vp(\cdot,t)\in {\cal K}$ is
\bea
\ddot\vp-|\na\dot\vp|_{\o_\vp}^2=0.
\eea
A key observation due to Donaldson \cite{D99} and Semmes \cite{Se}
is that this geodesic equation is equivalent to a homogeneous complex Monge-Amp\`ere equation 
\bea
\label{geo}
(\pi^*\o_0+{i\over 2}\ddb \Phi)^{n+1}=0.
\eea
for the function
\bea
\Phi(z,w)\equiv \vp(z,\log\,|w|)
\eea
on $M\equiv X\times A$, with $A=\{w\in {\bf C}; e^{-T}<|w|< 1\}$, 
Here $\ddb\Phi$ is taken with respect to all $(n+1)$ variables $(z,w)$,
and $\pi^*\o_0$ is the pull-back of $\o_0$ to $M$. Note that $\pi^*\o_0$ is not strictly positive, viewed as a $(1,1)$-form on $M$.

\smallskip
We shall consider both {\it geodesic segments}, joining two points $\vp_0$ and $\vp_1$
of ${\cal K}$, and {\it geodesic rays}, extending from a point $\vp_0\in {\cal K}$ to infinity. There is no loss of generality in taking $\vp_0=0$.

\smallskip
The equation for geodesic segments is then a standard Dirichlet problem
on the manifold $M=X\times A$, $\p A=\{w\in {\bf C};|w|=1\ {\rm or}
\ |w|=e^{-T}\}$ with $T$ finite, and boundary value $\Phi_b$ defined by
\bea
\Phi_b(z,w)&=& 0 \ \ {\rm for}\ \ |w|=1\nonumber\\
\Phi_b(z,w) &=& \ \vp_1\ \ {\rm for}\ \ |w|=e^{-T}.
\eea

The equation for geodesic rays is more unusual: here $T=\infty$, so the annulus $A$ reduces to the punctured disk $D^\times=\{w\in {\bf C}; 0<|w|< 1\}$, and
$M=X\times D^\times$. A boundary value $\Phi_b$ is
assigned only on the component $|w|=1$ of the boundary of $M$,
\bea
\Phi_b(z,w) =0 \ \ {\rm for}\ \ |w|=1
\eea
but there is no condition near $w=0$. In practice, and with motivation from geometric invariant theory, we shall restrict to the case where $\o_0=c_1(L)$, where $L\to X$ is a positive line bundle, and consider the geodesic rays associated to a {\it test configuration} of $L\to X$. A test configuration (see the precise definition in \S 13.2 below) is a one-parameter subgroup (1PS) degeneration of the line bundle $L\to X$. It produces a limiting singular line bundle, or more precisely a polarized scheme, $L_0\to X_0$, which can be viewed as an implicit boundary value for the Dirichlet problem for the homogeneous complex Monge-Amp\`ere equation at $w=0$. The net result is that we shall associate a canonical geodesic ray to each test configuration, starting from an arbitrary point $\vp_0\in {\cal K}$. Thus a test configuration can also be viewed as providing a direction where the Cauchy problem admits a generalized solution for infinite time, and hence as a generalized vector field on the space ${\cal K}$ of K\"ahler metrics.

\smallskip

The motivation for test configurations and geodesic rays is the following.
It has been shown by Donaldson \cite{D99} and Mabuchi \cite{M87}
that ${\cal K}$ is a symmetric space with non-positive curvature. The geodesic rays in ${\cal K}$ are then just a generalization of the one-parameter subgroups of finite-dimensional symmetric spaces of negative curvature. In Donaldson's program for the problem of constant scalar curvature metrics \cite{D99}, a numerical invariant for geodesic rays can 
be defined as the limiting value as $t\to\infty$ of the rate of change of the Mabuchi K-energy along geodesic rays. An infinite-dimensional notion of stability in geometric invariant theory (GIT) can then be defined as the negativity of this numerical invariant of geodesic rays, and we obtain in this way an infinite-dimensional version of the Yau-Tian-Donaldson conjecture, which asserts the equivalence between the existence of a metric of constant scalar curvature in the K\"ahler class $c_1(L)$ and the K-stability of $L\to X$ in the sense of GIT.

\medskip
Very recently, Lempert and Vivas \cite{LV} have produced examples of K\"ahler manifolds $(X,\o_0)$, specifically tori with certain symmetries, where there are no
$C^3$ geodesic connecting two given potentials $\vp_1,\vp_2\in {\cal K}$. A $C^\infty$ geodesic would correspond to a function $\Phi(z,w)$ which is $C^\infty$
in both $z$ and $w$, a solution of the completely degenerate Monge-Amp\`ere equation, and which is strictly plurisubharmonic with respect to $\o_0$ for each
$w\in A$. On the other hand, generalized solutions of (\ref{geo}) in the sense of pluripotential theory, where $\Phi(z,w)$ is only known to be plurisubharmonic in all variables $(z,w)$, and where $\Phi$ is only of class $C^{1,\alpha}$, have been known for some time. Their existence can be deduced from the general theory of boundary value problems of \S 9 \cite{C00, B09b, PS09a}, as well as from explicit semi-classical constructions \cite{PS06, PS07, SZ07, SZ10, RZ08}. 
Other constructions are due to
Arezzo-Tian \cite{AT} (by the Cauchy-Kowalevska method), and to
Chen \cite{C06}, Chen-Sun \cite{CS}, Chen-Tang \cite{CT}, Ross-Witt Nystrom \cite{RW} and other authors under various types of assumptions.
A partial regularity theory for the homogeneous complex Monge-Amp\`ere equation has been proposed by Chen and Tian \cite{CT}.
We discuss some of these developments in this section.

\subsection{Geodesics from a priori estimates}

We show how generalized geodesics can be obtained from the existence theorems for the Dirichlet problem of Section \S 8. In the present context, $(X,\o_0)$ is a given compact K\"ahler manifold without boundary, and the manifolds 
with boundary of the theorems in Section \S 8 are now given by $M=X\times A$,
with $A$ an annulus, or $M=X\times D^\times$, with $D^\times$ a punctured disk.

\smallskip
First, we consider the case of geodesic segments, linking $\vp_0=0$ to $\vp_1\in {\cal K}$. Here the set-up of the geodesic equation is exactly the same as for Theorem \ref{HCMA1}, with the difference that the form $\pi^*\o_0$ in
(\ref{geo}) is not strictly positive. However, it is easy to bring ourselves back to the case of a strictly positive form, by constructing a smooth function
$\underline{\Phi}$ satisfying
\bea
\pi^*\o_0+{i\over 2}\ddb\underline{\Phi}>0
\ {\rm on}\ M,
\qquad
\underline{\Phi}=\Phi_b\ {\rm on}\ \p M.
\eea
Setting then $\Omega_0=\pi^*\o_0+{i\over 2}\ddb\underline{ \Phi}$ and
$\Phi=\Psi+\underline{\Phi}$, the equation
$(\pi^*\o_0+{i\over 2}\ddb\Phi)^{n+1}=0$ is equivalent to the equation
$(\Omega_0+{i\over 2}\ddb\Psi)^{n+1}=0$,
with $\Omega_0$ now a K\"ahler form so that Theorem \ref{HCMA1} applies at once. The function $\underline{\Phi}$ is obtained by the following elementary construction (see e.g. Lemma 14 of \cite{PS08}): $\underline{\Phi}
=t \vp_1+f(w)$, $t=\log |w|$, with $f(w)$ the solution of the Dirichlet problem $\Delta f(w)=C$ on $A$, $f=0$ on $\p A$, and $C$ is a large positive constant.
Thus we have established the following theorem, first proved by Chen \cite{C00}:

\begin{theorem}
\label{Chen}
Let $\vp_0$ and $\vp_1$ be two points in ${\cal K}$. Then there exists a unique generalized geodesic of class $C^{1,\alpha}$, for any $0<\alpha<1$, joining $\vp_0$ and $\vp_1$.
\end{theorem}

We turn next to geodesic rays. The notion of test configuration alluded to in \S 13.1
can be defined precisely as follows (see Donaldson \cite{D02}):

\medskip
\begin{definition}
\label{testconfiguration}
Let $L\to X$ be a positive line bundle over a compact complex manifold $X$.
A test configuration ${\cal T}$ for $L\to X$ consists of

{\rm (1)} a scheme ${\cal X}$ with a ${\bf C}^\times$ action $\rho$;

{\rm (2)} a ${\bf C}^\times$ equivariant line bundle ${\cal L}\to {\cal X}$, ample on all fibers;

{\rm (3)} and a flat ${\bf C}^\times$ equivariant map $\pi:{\cal X}\to {\bf C}$, where ${\bf C}^\times$ acts on ${\bf C}$ by multiplication, with the property that $(\pi^{-1}(1),
{\cal L}_{\vert_{\pi^{-1}(1)}})$ is isomorphic to $(X,L^r)$ for some $r>0$.
\end{definition}

We shall also denote ${\cal T}$ by
\bea
{\cal T}=\bigg(\rho:{\bf C}^\times \to {\rm Aut}({\cal L}\to {\cal X}\to {\bf C})\bigg)
\eea
where ${\rm Aut}({\cal L}\to {\cal X}\to {\bf C})$ is the space of automorphisms of the fibrations ${\cal L}\to {\cal X}\to {\bf C}$. A test configuration ${\cal L}$ is said to be trivial if ${\cal L}=L\times {\bf C}$, with the action $\rho(\tau)(\ell,w)=(\ell, \tau w)$, for $(\ell,w)\in L\times {\bf C}$, $\tau\in {\bf C}^\times$. 

It is convenient to introduce the notation $X_w=\pi^{-1}(w)$,
$L_w={\cal L}_{\vert_{\pi^{-1}(w)}}$, and to view ${\cal X}$ as ${\cal X}=\cup_{w\in{\bf C}}X_w$. A typical example of a test configuration would be ${\cal X}=\cup_\tau \sigma_\tau(X)$ where $X$ is a submanifold of ${\bf CP}^N$, and $\sigma_\tau=e^{\tau B}$, is a one-parameter subgroup of $GL(N+1)$.

Note that all fibers $(X_w,L_w)$ are biholomorphic to $(X_1,L_1)$ for $w\not=0$.
However, the ``central fiber" $(X_0,L_0)$ will usually have singularities. It can be viewed as the limit in the sense of schemes of $(X_w,L_w)$ as $w\to 0$, and it is invariant under the action of $\rho$.

\smallskip
For the construction of geodesic rays, we need the following geometric properties of a test configuration. Let
\bea
\label{resolution}
p:\tilde{\cal X}\to {\cal X}\to{\bf C}
\eea
a resolution of singularities, which can be chosen to be equivariant, in the sense
that the homomorphism $\rho$ lifts to a homomorphism $\tilde\rho:
{\bf C}^\times\to {\rm Aut}(p_*{\cal L}\to {\cal X}\to {\bf C})$, and that all diagrams commute. Let ${\cal X}_D=\pi^{-1}(D)$,
${\cal X}_D^\times=\pi^{-1}(D^\times)$,
and ${\cal L}_D\to {\cal X}_D=\pi^{-1}(D)$ and
${\cal L}_D^\times\to {\cal X}_D^\times$ be the fibrations
above $D$ and $D^\times$, with similar definitions for 
$\tilde{\cal X}_D$, $\tilde{\cal X}_D^\times$, and $\tilde{\cal L}_D$ and 
$\tilde{\cal L}_D^\times$. 
Let the K\"ahler form $\o_0$ on $X$ be the curvature of a metric $h_0$ on $L$.
Then an element $\vp$ of the space of K\"ahler potentials ${\cal K}$ can be identified with a metric $h=h_0 e^{-\vp}$ on $L$. We can now state the properties of test configurations which we need: 

\smallskip

$\bullet$ There is a biholomorphism 
\bea
M=X\times D^\times &\to& \tilde{\cal X}^\times\nonumber\\
L\times D^\times &\to& \tilde{\cal L}^\times
\eea
defined by $(\zeta,w)\to \rho(w)(\zeta)\in L_w$, for $\zeta\in L=L_1$.
If we extend the metric $h_0$ on $L$ trivially as a metric on $L\times D^\times$,
then it can be identified through the above biholomorphism
with a metric $H_0$ on $\tilde{\cal L}^\times$. 
The curvature of $H_0$ is (the push-forth of) $\Omega_0$. If $H$ is any other metric on $\tilde{\cal L}^\times$ and $\Omega$ its curvature, then $H=H_0e^{-\Psi}$, and the
equation $(\Omega_0+{i\over 2}\ddb\Phi)^m=0$ on $M$ is equivalent to
the equation
\bea
\label{HCMA2}
(\Omega+{i\over 2}\ddb \tilde\Phi)^{n+1}=0
\ \ {\rm on}\ \ \tilde{\cal X}_D^\times
\eea
with $\tilde\Phi=\Phi-\Psi$. 

$\bullet$ Recall that $h_0$ is any fixed metric on $L$ with positive curvature $\o_0$. It is shown in \cite{PS07a}, \cite{PS09b} that there is a metric $H$
on $\tilde{\cal L}$ (in particular over the whole of $\tilde{\cal X}_D$, including the central fiber), 
which restricts to $h_0$ on $X_1$,
and which has curvature $\Omega\geq 0$ over $\tilde{\cal X}$,
and $\Omega> 0$ on $\tilde{\cal X}_D^\times$. 

$\bullet$ Furthermore, there exists an effective divisor $E$ supported only in
the central fiber of $\tilde{\cal X}$ and a smooth metric $K$ on $O(E)$ so that
\bea
\Omega_K\equiv\Omega
+\e{i\over 2}\ddb \log K
\eea
is smooth and strictly positive definite everywhere on $\tilde{\cal X}_D$,
for all small, strictly positive $\e$.

\medskip
We return now to the problem of constructing generalized geodesic rays. Given a test configuration, the above constructions show that a generalized geodesic ray
is a solution of the equation (\ref{HCMA2}) on $\tilde{\cal X}_D^\times$.
The above geometric properties of test configurations mean precisely that the hypotheses of Theorem \ref{HCMA2} are satisfied with $M=\tilde{\cal X}_D$. Note that in this case, the boundary of the manifold is clearly holomorphically flat. Thus we obtain the following theorem \cite{PS09a}:

\begin{theorem}
\label{PS09}
Let $L\to X$ be a positive line bundle over a compact complex manifold $X$.
Let ${\cal T}$ be any test configuration of $L\to X$, and let $p$ be an equivariant resolution
as in (\ref{resolution}).
Then for any metric $h_0$ on $L$ with positive curvature,
there is a $C^{1,\alpha}$ generalized geodesic, with bounded potential, starting from $h_0$.
More precisely, there is a bounded solution $\tilde\Phi$ of the equation
$(\Omega+{i\over 2}\ddb \tilde\Phi)^{n+1}=0$ on $\tilde{\cal X}_D$,
which is of class $C^{1,\alpha}$ on $\tilde{\cal X}_D^\times$.
\end{theorem}

\subsection{Algebraic approximations: the Tian-Yau-Zelditch theorem}

The conjecture of Yau on K\"ahler-Einstein metrics \cite{Y93}, and the Yau-Tian-Donaldson conjecture on the existence of metrics of constant scalar curvature in a K\"ahler class $c_1(L)$
if $L$ is K-stable \cite{D02}, are conjectures relating the solvability of a geometric partial differential geometric to a global algebraic condition. One strategy, advocated by Yau
and implemented particularly successfully by Donaldson \cite{D01, D10} in the proof of the necessity of stability, has been to approximate transcendental objects such as metrics by 
their algebraic counterparts, such as Bergman kernels and Fubini-Study metrics. A precise
example is the Tian-Yau-Zelditch theorem \cite{Y93, T90a, Z} (also proved independently by Catlin \cite{Cat} and refined by Lu \cite{L}), which can be stated as follows.

\medskip
Let $L\to X$ be a positive line bundle over a compact complex manifold $X$.
If $\underline{s}=\{s_\alpha\}_0^{N_k}$
is a basis for the space $H^0(X,L^k)$ of holomorphic sections of $L^k$, ${\rm dim}\,
H^0(X,L^k)=1+N_k$, then the Kodaira map $\iota_{\underline s}$ is defined by
\bea
\iota_{\underline s}:\ X\ni z \to
[s_0(z):\cdots:s_{N_k}(z)]\in {\bf CP}^{N_k}.
\eea
The Kodaira imbedding theorem says that $\iota_{\underline s}$ is an imbedding for $k$ sufficiently large. The hyperplane bundle $O(1)$ over ${\bf CP}^{N_k}$ 
pulls back to $L^k$. The Fubini-Study metric $h_{FS}={1\over \sum_{\alpha=0}^{N_k}|s_\alpha|^2}$ on $O(1)$ and $\o_{FS}=-{i\over 2}\ddb \log \,h_{FS}$ on ${\bf CP}^{N_k}$ pull back 
then to the metrics
\bea
\iota_{\underline s}^*(h_{FS})={1\over
\sum_{\alpha=0}^{N_k}|s_\alpha(z)|^2},
\qquad
\iota_{\underline s}^*(\o_{FS})=-{i\over 2}\ddb \log \sum_{\alpha=0}^{N_k}|s_\alpha(z)|^2.
\eea
on $L^k$ and $X$ respectively. Note that $h(k)\equiv(\iota_{\underline s}^*(h_{FS}))^{1\over k}$ is then a metric on $L$, and its curvature $\o(k)\equiv
{1\over k}\iota_{\underline s}^*(\o_{FS})$ is a metric on $X$,
which can be viewed as algebraic objects.

Let now $h$ be a metric on $L$ with positive curvature $\o=-{i\over 2}\ddb\log h>0$.
The Tian-Yau-Zelditch theorem asserts that the metric $h$ on $L$ and $\o$ on $X$ can be approximated asymptotically in $k$ by the metrics $h(k)$ and $\o(k)$,
if the basis $\underline{s}=\{s_\alpha\}_0^{N_k}$ used to construct the Kodaira imbedding
is an orthonormal basis of $H^0(X,L^k)$ with respect to the $L^2$ metric, $\|s\|^2\equiv\int_X |s|^2 h^k \o^n$. To see how this comes about, we write
\bea
\label{FS}
\log {h(k)\over h}
=
-{1\over k}\log\,\rho_k(z),
\qquad
\o-\o(k)
=
-{i\over 2k}\ddb \log \,\rho_k(z)
\eea
where $\rho_k(z)$ is the Bergman kernel (or density of states, since it integrates to ${\rm dim}\,H^0(X,L^k)$),
\bea
\rho_k(z)\equiv \sum_{\alpha=0}^{N_k}|s_\alpha(z)|^2 h^k(z).
\eea
The desired statement follows then from an asymptotic expansion for $\rho_k(z)$ \cite{Z}
\bea
\rho_k(z)=k^n(1+\sum_{p=1}^\infty A_p(z)k^{-p}),
\eea 
This expansion is itself a consequence of an asymptotic expansion 
obtained by Boutet de Monvel and Sj\"ostrand \cite{BS}
for the Szeg\"o kernel on strongly pseudo-convex domains, refining an
earlier expansion along the diagonal obtained by Fefferman \cite{F}.
The coefficient $A_1(z)$ has been shown by Lu \cite{L} to be given by ${1\over 2}R(z)$, where $R$ is the scalar curvature of $\o$. This turns out to provide a key link with the problem of constant scalar curvature metrics \cite{D01, D10}. Other generalizations and applications of asymptotic expansions and of the Tian-Yau-Zelditch theorem can be found in \cite{MM, Wg}.

\subsection{Semi-classical constructions}

Let ${\cal K}_k$ denote the space  $\{(\iota_{\underline s}(h_{FS}))^{1\over k}\}$
of pull-backs of the Fubini-Study metrics, under Kodaira imbeddings defined by an arbitrary
basis $\underline s$ of $H^0(X,L^k)$. Since $h_{FS}$ is invariant under $SU(N_k+1)$, we have
\bea
{\cal K}_k
=
SL(N_k+1)/SU(N_k+1).
\eea
Note that the right hand side is a symmetric space with negative curvature.
A suggestive consequence of the Tian-Yau-Zelditch theorem is that, in a pointwise sense,
we have
\bea
\label{approximation}
{\cal K}={\rm lim}_{k\to\infty}{\cal K}_k
=
{\rm lim}_{k\to\infty}SL(N_k+1)/SU(N_k+1).
\eea
It is a natural question whether this pointwise approximation can be extended to the approximation of more geometric properties, for example of extended geometric objects.
A prime example is whether geodesics in ${\cal K}$ can be approximated by geodesics in ${\cal K}_k$. This translates precisely into whether the solutions of the homogeneous complex Monge-Amp\`ere equations can be approximated by one-parameter subgroups of Bergman kernels \cite{PS06}. We shall see below that the answer is affirmative, see
\cite{PS06, PS07, PS09b}.
Some refinements of these approximations and their rate of convergence can be found
in Berndtsson \cite{Be1, Be2},
and in \cite{SZ07, SZ10} in the case of toric varieties. For toric varieties, a similar approximation has been extended to harmonic maps by Rubinstein and Zelditch \cite{RZ08}. 

\medskip
We provide now some details. Fix a metric $h_0$ on $L$ with positive curvature $\o_0$.
Let $\Phi(z,w)$ be a solution of the Dirichlet problem (\ref{geo})
with boundary value $\vp_0$ when $|w|=1$. If we view it as either
a geodesic segment or a geodesic ray in ${\cal K}$ emanating from the corresponding potential $\vp_0$, then this geodesic segment or ray should be the limit
of a sequence of one-parameter subgroups in ${\cal K}_k$, as $k\to\infty$. Let $B_k
\in GL(N_k+1)$ be the infinitesimal generator
of the one-parameter subgroup in ${\cal K}_k$, so that the subgroup is given by $w^{B_k}$, $w\in {\bf C}^\times$. We can assume that $B_k$ is diagonal,
with eigenvalues $\lambda_\alpha^{(k)}$, $0\leq\alpha\leq N_k$.
Let $\underline{s}=\{s_\alpha\}_0^{N_k}$ be a basis of $H^0(X,L^k)$
which is orthonormal with respect to the $L^2$ norm defined by $h_0$ and the volume form $\o_0^n$. The subgroup $w^{B_k}$ acts on the basis $\{s_\alpha\}_0^{N_k}$ to produce the basis
$w^{B_k}\cdot\underline{s}\equiv
\{w^{\lambda_\alpha^{(k)}}s_\alpha\}_0^{N_k}$. The corresponding pull-backs of the Fubini-Study metrics can be written explicitly as
\bea
\label{Phik}
\Phi_k(z,w)
=
{1\over k}\log\sum_{\alpha=0}^{N_k}|w|^{2\lambda_\alpha^{(N_k)}}|s_\alpha(z)|^2h_0(z)^k
\eea
Thus the main problem is to choose the appropriate generators $B_k$ and to show that
the $\Phi_k(z,w)$ converge, in a suitable sense, to a solution $\Phi(z,w)$ of the homogeneous complex Monge-Amp\`ere equation. Note that, by the Tian-Yau-Zelditch theorem, the functions $\Phi_k(z,w)$ converge to the correct boundary value when $|w|=1$.

\smallskip
We address the problem of choosing the generators $B_k$. 

\smallskip
Consider first the case of a generalized geodesic segment joining two points $\vp_0$ and $\vp_1$ in ${\cal K}$.
There is in this case a natural choice of infinitesimal generator $B_k$: $B_k$ is just the matrix of change of bases, from a basis $\{s_\alpha^{(0)}\}_0^{N_k}$ orthonormal with respect to the $L^2$ metric defined by $h_0$, $\o_0^n$, to a basis $\{s_\beta^{(1)}\}_0^{N_k}$ orthonormal with respect to the $L^2$ metric defined by $h_1$, $\o_1^n$.

Next, consider the case of a test configuration. 
The group action $\rho(w)$, $w\in {\bf C}^\times$ preserves the central fiber $(X_0,L_0)$.
Thus it induces a one-parameter subgroup $w^{B_k}$
on the space of holomorphic sections,
\bea
w^{B_k}: \ H^0(X_0,L_0^k)
\to H^0(X_0,L_0^k).
\eea
The generators $B_k$ are the generators that we are looking for.

\smallskip
Once the generators $B_k$, and hence their eigenvalues $\lambda_\alpha^{(k)}$ have been chosen, we need a criterion for when the expressions (\ref{Phik}) have the desired convergence properties. This is provided by the following lemma
\cite{PS08}:
 
\begin{lemma}
Fix $h_0\in{\cal K}$
as before. Consider the general Ansatz
\begin{equation}
\Phi_k(z,w)
=
{1\over k}
\log\,\sum_{\alpha=0}^{N_k}|w|^{2\lambda_\alpha^{(k)}}|^2|s_\alpha^{(k)}(z)|_{h_0^k}^2
\nonumber
\end{equation}
where $\{s_\alpha^{(k)}(z)\}$ is an orthonormal basis for $H^0(X,L^k)$ with respect to
the $L^2$ norm defined by $(h_0^k,\omega_0^n)$, and $\lambda_\alpha^{(k)}$
are real numbers, $0\leq\alpha\leq N_k$ for each $k$.
Then if

\smallskip

{\rm (1)} There exists a constant $C$ independent of both $\alpha$ and $k$ so that
\begin{equation}
|\lambda_\al^{(k)}|\ \leq \ C\,k
\nonumber
\end{equation}

{\rm (2)} There exists a constant $C_T$ independent of $k$ so that
\begin{equation}
\int\int_{X\times \{e^{-T}<|w|\leq 1\}}(\pi^*\o_0+{i\over 2}\ddb\Phi_k)^{n+1}
\leq \ C\,k^{-1}.
\nonumber
\end{equation}
then
$\Phi(z,w)={\rm lim}_{k\to\infty}[{\rm sup}_{\ell\geq k}\Phi_k(z,w)]^*$
is continuous at $|w|=1$, and satisfies in the sense of pluripotential theory
\begin{equation}
(\pi^*\o_0+{i\over 2}\ddb\Phi)^{n+1}=0
\ {\rm on}\ X\times\{e^{-T}<|w|<1\},
\qquad \Phi(z,w)=0\ {\rm for}\ |w|=1.
\nonumber
\end{equation}

\end{lemma}

It can be shown that, in both choices of generators $B_k$ for geodesic segments and for geodesic rays, the condition on the uniform growth of the eigenvalues $\lambda_\alpha^{(k)}$ is satisfied. The only non-trivial step remaining is to verify the condition on the decay of the masses of the Monge-Amp\`ere measures of $\Phi_k(z,w)$. The key observation here is that these masses are essentially cohomological, and given by 
\begin{equation}
\label{MAmass}
\int\int_{X\times \{e^{-T}<|w|<1\}}
(\pi^*\o_0+{i\over 2}\ddb\Phi)^{n+1}
=
\dot E(\vp(\cdot,0))-\dot E(\vp(\cdot,-T))
\nonumber
\end{equation}
where $E(\vp)$ is the functional (\ref{E}) we had encountered earlier
as the anti-derivative of the Monge-Amp\`ere measure $\o_\vp^n$.

\smallskip

We can now state the main theorems of this section and complete their proofs. For geodesic segments, we have \cite{PS06}:

\begin{theorem}
\label{segments-bergman}
Let
$h_1\in {\cal K}$
be another metric on
$L$ with
$\omega_1=-{i\over 2}\ddb\log\,h_1>0$, $M\equiv X\times \{e^{-1}< |w|< 1\}$.
Then the generalized solution $\Phi(z,w)$ of the Dirichlet problem (\ref{geo})
can be expressed as
\begin{equation}
\Phi(z,w)={\rm lim}_{k\to\infty}[sup_{\ell\geq k}\Phi_\ell(z,\log|w|)]^* \nonumber
\end{equation}
where the eigenvalues $\lambda_\alpha^{(k)}$ defining $\Phi_k(z,w)$ are the eigenvalues of the matrix $B_k$
of change of bases from an orthonormal basis of $H^0(X,L^k)$
with respect to $h_0,\o_0^n$ to an orthonormal basis of $H^0(X,L^k)$
with respect to $h_1,\o_1^n$.
\end{theorem}

\noindent
{\it Proof}: We compute explicitly the right hand side of (\ref{MAmass}).
Denote the approximating Fubini-Study metrics in (\ref{FS})
for the metrics $h_0$ and $h_1$ as $h_0(k)=e^{-\vp_0(k)}h_0$
and $h_1(k)=e^{-\vp_1(k)}h_0$ respectively. Let
$\o_a(k)$ be their curvatures. Let $\{s_{a,\alpha}^{(k)}\}$
be orthonormal bases with respect to $h_a,\o_a^n$.
Then we have, with $a=0,1$,
\bea
\int_X\dot\vp_a(k)\o_a(k)^n
&=&{2\over k^{n+1}}\int_X\sum_{\alpha=0}^{N_k}
\lambda_\alpha^{(k)}|s_{a,\alpha}^{(k)}(z)|^2h_a(k)^k\o_a(k)^n
\nonumber\\
&=&
{2\over k^{n+1}}\sum_{\alpha=0}^{N_k}\lambda_\alpha^{(k)}
+
O({1\over k^{n+2}})N_k\,{\rm max}_\alpha|\lambda_\alpha^{(k)}|,
\eea
where we have applied the Tian-Yau-Zelditch theorem. Thus the leading terms cancel
in the difference (\ref{MAmass}), giving the desired estimate. Q.E.D.

\medskip
By the uniqueness of the solution of the Dirichlet problem, this solution must coincide with the solution obtained by Chen \cite{C00}
from the method of a priori estimates (see Theorem \ref{Chen}), so it must be $C^{1,\alpha}$.
It has also been shown by Berndtsson \cite{Be2}
that the convergence described in \cite{PS06} can actually be strengthened to uniform convergence. 

\medskip
For geodesic rays defined by a test configuration, we have \cite{PS07, PS09b}:

\begin{theorem}
\label{rays-bergman}
Let $\rho$ be a test configuration for a positive line bundle
$L\to X$ over a compact complex manifold $X$. Let $h_0$ be a metric on $L$,
with $\omega_0=-{i\over 2}\ddb\,\log\,h_0>0$. Let 
$\Phi_k(z,w)$ be defined by (\ref{Phik}), where the eigenvalues $\lambda_\alpha^{(k)}$
are the eigenvalues of the endomorphisms $B_k$
on $H^0(X_0,L_0^k)$ induced by the
group action $\rho$. 
Then
\begin{equation}
\Phi(z,w)={\rm lim}_{k\to\infty}[{\rm sup}_{\ell\geq k}\Phi_k(z,w)]^*
\nonumber
\end{equation}
defines a generalized solution of the Dirichlet problem
\begin{equation}
(\pi^*\o_0+{i\over 2}\ddb\Phi)^{m+1}=0\
{\rm on}\ X\times D^\times,
\qquad\Phi(z,w)=0\ {\rm for}\ |w|=1.
\nonumber
\end{equation}
The solution is actually of class $C^{1,\alpha}(X\times D^\times)$ for any $0<\alpha<1$.
It is non-constant when the test configuration is non-trivial.

\end{theorem}

\noindent
{\it Proof}. As before, it remains only to prove the bound $O(k^{-1})$ on the mass of the Monge-Amp\`ere measure on $M=X\times D^\times$. 
Consider the functions $\Phi_k^\#(z,w)$ defined by the same formula as $\Phi_k(z,w)$, but with the eigenvalues $\lambda_\alpha^{(k)}$ replaced by their traceless
counterparts
\bea
\lambda_\alpha^{\#,(k)}
=\lambda_\alpha^{(k)}-{{\rm Tr}\,B_k\over N_k+1}.
\eea
Then $\Phi_k(z,w)=\Phi_k^\#(z,w)(z,w)+{{\rm Tr}\,B_k\over k(N_k+1)}\,\log\,|w|^2$.
Thus they have the same complex Hessian, and we can evaluate
(\ref{MAmass}) with $\Phi_k(z,w)$ replaced by $\Phi_k^\#(z,w)$.  
The formula (\ref{MAmass}) gives then, with obvious notations,
\bea
\int_{X\times D^\times}(\pi^*\o_0+{i\over 2}\ddb\Phi_k^\#)^{n+1}
=
{\rm lim}_{T\to\infty}\int_X\dot\vp(T)^\#\o_k(T)^n
-
\int_X\dot\vp(0)^\#\o_k(0)^n.
\eea
Since the eigenvalues $\lambda_\alpha^{\#(k)}$ sum to $0$,
the leading term in the second expression on the right hand side is $0$. As for the first expression, the lemma below shows that it is automatically $O(k^{-1})$. This establishes the fact that the limit $\Phi(z,w)$ satisfies the homogeneous complex Monge-Amp\`ere equation.

\begin{lemma}
\label{DF}
Let ${\cal T}$ be a test configuration. Then we have
\bea
{\rm lim}_{T\to\infty}\int_X\dot\vp(T)\o_k(T)^n={1\over k}F
\eea
where $F$ is the Donaldson-Futaki invariant of ${\cal T}$, defined as the second term $F$ in the following asymptotic expansion
\bea
{{\rm Tr}\,B_k\over k(N_k+1)}=F_0+F\,k^{-1}+O(k^{-2}).
\eea
\end{lemma}

We should say that this lemma was implicit in the paper of Donaldson \cite{D04}.
Its explicit statement and proof can be found in \cite{PS07}, Mabuchi \cite{M09},
and Donaldson \cite{D10}.

\medskip
Unlike in the case of the geodesic segments, the regularity of the geodesic rays obtained from the above theorem does not follow as yet, since the behavior of the ray near $w=0$ has not been addressed. The following lemma gives a complete description of this behavior
\cite{PS09b}:

\begin{lemma}
Let ${\cal T}$ be a test configuration, and let $p$ be an equivariant resolution of singularities as considered earlier in \S 13.2 . Then the function 
\bea
\Psi_k\equiv \Phi_k-\Phi_1
\eea
extends as a smooth function over the whole of $\tilde{\cal X}_D$. Furthermore, it satisfies the following uniform estimate
\bea
{\rm sup}_{k\geq 1}{\rm sup}_{\tilde{\cal X}_D}|\Psi_k|
\leq C <\infty.
\eea
\end{lemma}

With this lemma, we can show that the function
\bea
\Psi=\Phi-\Phi_1
\eea
is a bounded solution of a homogeneous complex Monge-Amp\`ere equation on $\tilde{\cal X}_D$ with a non-negative background form $\Omega_1$. This can be expressed in turn as a homogeneous complex Monge-Amp\`ere equation with a background form $\Omega$ which satisfies all the hypotheses of Theorem 14.
Thus the solution must be $C^{1,\alpha}$, and the proof of the theorem is complete.

\subsection{The toric case}

In the toric case, the previous constructions of solutions of homogeneous complex Monge-Amp\`ere equations as limits of Bergman metrics
can be analyzed more precisely.
We can obtain in this manner more detailed information on the approximating paths and their rates of convergence. A remarkable feature also emerges from
this study, which is an unexpected 
relation between the previous semiclassical constructions and the theory of large deviations \cite{SZ07, SZ10}.

\subsubsection{Bergman geodesics}

Let $X$ be an $n$-dimensional  toric manifold and $L\rightarrow X$ be a positive  toric line bundle over $X$.  Let $\cal K_{T}$ be the space of positively curved smooth toric hermitian metrics on $L$ which are invariant under the compact 
$(S^1)^n$ torus action. Let $h_0, h_1\in \cal K_T$ and let $h_t$ for  $0\leq t\leq 1$ be the Monge-Amp\`ere geodesic between them as defined in Section \S 13.4. We define 
\bea
\varphi_t (z)=\log \left(  h_t (h_0)^{-1}\right) .
\eea
Then $ Ricci(h_t) = Ricci(h_0) +{ i\over 2} \ddbar \varphi_t. $

\smallskip
The line bundle $L$ is associated to a convex polytope $P$ in $\R^n$ which coincides  with the image of the moment map by any toric K\"ahler metric in $c_1(L)$. Each integral point $\alpha \in k \overline P$ corresponds to a holomorphic section in $H^0(X, L^k)$. In particular, $\{z^\alpha\}_{\alpha \in k\overline{P}} $ form a basis for $H^0(X, L^k)$ for $z\in (\C^*)^n$. Let $\{ s_\alpha(z)\}_{\alpha\in k \overline P\cap \Z^n} $ be an orthornormal toric basis for $H^0(X, L^k)$ with respect to the $L^2$ norm defined by $(h_0^k, \omega_0^n)$, where $\omega_0= -{i\over 2}\ddb \log h_0$. Then we define the following 
paths of Bergman metrics,
\bea
h_{t,k} (z)= \sum_{\alpha\in k \overline P\cap \Z^n}   \frac{1}{ \left( Q_{h_0^k}(\alpha)\right)^{1-t} \left( Q_{h_1^k}(\alpha)\right)^t}  |s_\alpha (z)|^2,
\eea
where 
\bea
Q_{h_0^k} (\alpha) =\|s_\alpha\|^2_{h_0^k} = \int_X |s_\alpha(z)|_{h_0^k}^2 \omega_0^n, ~~~~Q_{h_1^k} (\alpha) =\|s_\alpha\|^2_{h_1^k} = \int_X |s_\alpha(z)|_{h_1^k}^2 \omega_1^n.
\eea
The corresponding potentials are given by
\bea
\varphi_k (t, z)  =\frac{1}{k}  \log \left( h_{t,k} h_0^{-k}\right).
\eea
They correspond to the potentials $\Phi_k(z,w)$ of Theorem
\ref{segments-bergman} with $t=\log |w|$.

\medskip

The following theorem \cite{SZ10} shows that
the potentials $\vp_k(t,z)$ actually converge in $C^2$:

\begin{theorem} \label{szel1}
We have
\bea
\lim_{k\rightarrow \infty} \|\varphi_k (t, z) - \varphi_t(z)\|_{C^2([0, 1]\times X)} = 0.
\eea

\end{theorem}

Theorem \ref{szel1} is a considerable strengthening of both Theorem \ref{segments-bergman} and Berndtsson's result
\cite{Be2} in the toric case. The advantage of studying Bergman metrics on toric manifolds is that toric holomorphic sections are naturally orthogonal to each other and one can analyze the norming constants $Q_{h^k}(\alpha)$. On the other hand, toric geodesics $\varphi_t$ are alway smooth \cite{Gd}.
Thus one may expect a higher order of convergence.  

\smallskip

In \cite{RZ08}, Theorem \ref{szel1} is generalized from geodesics of toric K\"ahler metrics to harmonic maps of a compact Riemannian manifold with boundary into the space of toric K\"ahler metrics.  More precisely, such a harmonic map equation can always be solved, and the solution approximated by harmonic maps into the space of toric Bergman metrics in the $C^2$-topology.  The case of geodesics corresponds to the case when the Riemannian manifold with smooth boundary is the interval $[0, 1]$.

\subsubsection{Geodesic rays and large deviations}

The approximation of the geodesic rays associated in Section \S 13.4 to a test configuration ${\cal T}$ for a polarization $L\to X$ can similarly be refined in the case of toric manifolds
\cite{SZ07}. An interesting observation is made in \cite{SZ07},
which  relates the geodesic rays constructed in \S 13.4 on toric manifolds to the large deviations principle \cite{V}.  

\medskip

Let $L \rightarrow X$ be a very ample toric line bundle over a toric manifold $X$. We use the same notations as in Section \S 13.5.1.  Let $h= e^{-\vp}$ a smooth toric hermitian metric on $L$ such that $\omega= {i\over 2} \ddbar \vp \in c_1(L)$ is a toric  K\"ahler metric on $X$. On $({\bf C}^n)^*$,the potential $\vp$ can be identified with a smooth convex function $\psi$ on ${\bf R}^n$,
$\vp(z)=\psi(\rho)$, with $\rho= \log |z|^2$. The Legendre transform of $\psi$
defined by
\bea
u(x)= \sup_{\rho \in \R^n} ( x\cdot \rho - \psi(\rho))
\eea 
is called the symplectic potential associated to $\psi$. 
The function $u(x)$ is a smooth convex function on $P$ with appropriate boundary singularities. 

\medskip

Let $\psi_t (\rho)$ be the geodesic ray constructed in
Theorem \ref{rays-bergman} in the toric setting. Then the symplectic potential associated to $\psi_t$ is given by $$ u_t (x) = u_0(x) - t(R- f(x))$$ for  some positive piecewise linear convex function $f(x)$ on $P$, and $R\in \R$ with $R-f(x)>0$ on $P$. The piecewise function $f(x)$ is an alternative way of describing a test configuration ${\cal T}$ in the toric setting \cite{D02}.

\smallskip

Define the pair $(d\mu_k^\rho, I^\rho(x))$  by
\bea
d \mu_k^\rho(x)&=& (\Pi_{h^k}(z,z))^{-1} \sum_{\alpha \in k \overline P \cap \Z^n} \frac{|s_\alpha|^2_{h^k}(z) }{Q_{h^k}(\alpha)} \delta_{\frac{\alpha}{k}}(x)  
\nonumber\\
I^\rho(x) &=& u_0(x) + \psi_0 - x\cdot \rho 
\eea
where $\Pi_{h^k}(z,z)$ is the Szeg\"o kernel for $(L^k, h^k)$ and $\delta_{\alpha /k} (x)$ is a delta function at $\alpha/k$.  
Then the pair $(d\mu_k^\rho, I^\rho(x))$ satisfies the large deviation principle.
The measure
$d\mu_k^\rho$ is a probability measure on $\overline P$ and the function $I^\rho(x)$ is called the rate function associated to $d\mu_k^\rho$.  Varadhan's lemma says that for each $t$ and $\rho$, 
\bea
\lim_{k\rightarrow \infty} \frac{1}{k} \log \int_P e^{k t(R- f(x))} d \mu_k^\rho(x) = \sup_{x\in P} \left(  t(R-f(x)) - I^\rho(x) \right).
\eea
However, it turns out that 
\bea
\psi_{t,k} (\rho) - \psi_0 (\rho) = \frac{1}{k} \log \int_P e^{kt(R- f(x))} d\mu_k^\rho(x)
\eea 
is exactly the Bergman geodesic ray constructed in Theorem \ref{rays-bergman}, while 
\bea
\sup_{x\in P} (  t(R- f(x)) - I^\rho(x)) &=& \sup_{x\in P} \left\{ x\cdot \rho - (u_0(x) - t (R- f(x)) \right\}  - \psi_0(\rho) \nonumber\\
&=&\psi_ t (\rho)- \psi_0(\rho) .
\eea 
Therefore Varadhan's lemma immediately gives the pointwise convergence of $\psi_{k, t} $ to $\psi_t$.
The uniform convergence in $C^1(X\times [0,1])$ topology is proved in \cite{SZ07}. The toric geodesic ray $\psi_t$ is also shown in \cite{SZ07} to be $C^{1,1}(X\times [0,1])$, but not  $C^2(X\times [0,1])$ in general.

\subsubsection{Counter-examples to regularity of higher order than $C^{1,1}$}

It is well-known that solutions of the homogeneous real Monge-Amp\`ere equation may be only of class $C^{1,1}$ and not higher. Such counterexamples have been extended to the complex case by Gamelin and Sibony \cite{GS}. A more general argument for why solutions cannot always be smooth has been given by Donaldson \cite{D}. In this section, we would like to give a simple example which illustrates the fact that, even when geodesic segments may be smooth, as in the case of toric varieties, geodesic rays associated to a test configuration may be again at most $C^{1,1}$
\cite{SZ07}.

\medskip
Consider the following simple example of a $C^{1,1}$ geodesic ray over ${\bf CP}^1$.  We consider the standard Fubini-Study metric $g_{FS} =\frac{i}{2}\ddbar \phi_0=  \frac{i}{2} \ddbar \log (1+|z|^2)=\frac{i}{2}\ddbar \log (1+ e^\rho)$ on ${\bf CP}^1$, where $\phi_0=\log(1+e^\rho)$ and $\rho = \log |z|^2$. Then the moment map can be constructed by 
\bea
x= \frac{ \partial \log (1+e^\rho)}{\partial \rho} = \frac{ e^\rho}{1+e^\rho} \in (0, 1).
\eea 
The symplectic potential $u_0$ corresponding to $\phi_0$ is given by 
\bea
u_0 (x)=  x \rho - \phi_0(\rho)=x \log x + (1-x) \log (1-x), ~~~~~x = \frac{e^\rho}{1+e^\rho} .
\eea 

Let $f(x) =  |x-1/2|$ be a piecewise linear convex function. Then 
\bea
u_t(x) = u_0(x) + t f(x)= x\log x+ (1-x) \log (1-x) + t |x-1/2|
\eea
induces a geodesic ray. Now we can calculate the K\"ahler potential $\phi_t(\rho)$ corresponding to $u_t$. 

By applying the Legendre transform to $u_t$, we can show by a straightforward calculations that
\bea
\phi_t (\rho) = \rho x - u_t( x), ~~~~ \rho = \frac{ \partial u_t}{\partial x},
\eea
and
\begin{equation} 
\phi_t (\rho) =  \left\{  \begin{array}{ll}
    - \frac{t}{2} + \log ( 1 + e ^{\rho+t}), &  \rho\in (-\infty, -t) \\ 
    \frac{\rho}{2} + \log 2,  & \rho \in (-t, t) \\ 
    \frac{t}{2} + \log (1+ e^{\rho- t} ), & \rho\in (t, \infty) \\ 
  \end{array} \right.
\end{equation}

Let $\varphi(t, z) = \phi_t(\rho) - \phi_0$. Then $\varphi$ is an  $\omega_{FS}$-psh function in $C^{1,1}({\bf R} \times X)$, but it is not of class $C^2 ({\bf R}\times X)$.

\subsection{The Cauchy problem for the homogeneous Monge-Amp\`ere euation}

The construction of geodesic rays associated to a test configuration can be viewed
as the solution of the Cauchy problem for the homogeneous complex Monge-Amp\`ere equation, with the initial velocity provided, implicitly, by the test configuration.
Even though the Cauchy problem is not well-posed in the sense of Hadamard, it is instructive to examine when there are solutions and when and how they break down. Such an analysis has been provided recently by Rubinstein and Zelditch \cite{RZ10a} for convex solutions of the real Monge-Amp\`ere equation.

Consider the Cauchy problem
\bea
&&
\det (\nabla^2 \phi) =0\ {\rm on}\ [0, T] \times\R^n,
\nonumber\\
&&
\phi(0, x)=\phi_0(x)\ {\rm on}\ \R^n,
\qquad 
 \frac{\partial \phi}{\partial t} (0, x)=\psi_0(x) \ {\rm on}\ \R^n. 
\eea

It is shown in \cite{RZ10b} that the above equation can be solved via Legendre transform until the Legendre transform of the solution stops being convex. More precisely, the Legendre transform of $\phi_t$ is given by 
\bea
u_t = u_0 + t  v.
\eea 
The function $u_t$ will stop being convex at a certain $T_{span}>0$ if $v$ is not convex, thus $\psi_t$ stops solving the HRMA after $T_{span}$. A candidate solution is constructed in \cite{RZ10b} for $t\geq T_{span}$, however, it does not solve the homogeneous real Monge-Amp\`ere equation even in a weak sense and it is not differentiable in general. Before $T_{span}$, the solution can also be approximated by Toeplitz quantization of the Hamiltonian flow defined by the Cauchy data \cite{RZ10a}.

\section{Envelopes and the Perron Method}
\setcounter{equation}{0}

The complex Monge-Amp\`ere equation satisfies the comparison principle. Thanks to this, generalized solutions can be obtained by the Perron method, as envelopes of families of plurisubharmonic functions. In this section, we describe some results obtained in this manner, focusing on the homogeneous case.

\subsection{Envelopes}

The Perron method for complex Monge-Amp\`ere equations was first developed by Bedford and Taylor \cite{BT76} for degenerate complex Monge-Amp\`ere equations on bounded domains in ${\bf C}^n$. A special case of their results of particular interest to our considerations is the following. Let $D\subset {\bf C}^n$ be a smooth,
bounded strictly pseudoconvex domain in ${\bf C}^n$, and let $f\in C(\p D)$. Define the following family of plurisubharmonic functions
\bea
\E_{D, f} = \{ u \in \PSH(D)~|~ u |_{\partial D} \leq f \}
\eea
and its upper envelope
\bea
\hat u (z)= \sup_{\E_{D, f}} u(z).
\eea
Then $\hat u\in \E_{D,f}\cap C(\overline D)$, $u$ is of class $C^{1,1}$ in the interior of $D$, and $u$ is the unique solution of the Dirichlet problem for the homogeneous complex Monge-Amp\`ere equation
\bea
(\frac{i}{2}\ddb u)^n=0 \ \ {\rm on}\ \ D,
\qquad
u |_{\partial D} = f \ \ {\rm on}\ \ \p D.
\eea
The local $C^{1,1}$ regularity follows from the $C^{1,1}$ regularity of the solution of the Dirichlet problem for the unit ball \cite{BT76}.

\medskip
The Perron method has been widely applied since for the complex Monge-Amp\`ere equations, and there has been considerable progress, thanks partly to the infusion of new techniques and improved approximation theorems for plurisubharmonic functions. In particular, the following results were obtained relatively recently by Berman and Demailly for certain homogeneous complex Monge-Amp\`ere equations, generalizing the geodesic equations considered in \S 13.

\medskip

Let $(X, \omega_X)$ be an $n$-dimensional compact K\"ahler manifold with a smooth K\"ahler form $\omega_X$. Let $\Sigma$ be a strictly pseudoconvex domain in ${\bf C}^m$ with $\rho$ being a smooth strictly plurisubharmonic defining function for $\Sigma$ and $\Sigma= \{ \rho <0\}$. We consider the product manifold $M = X\times \Sigma$ and let $\omega_M = \omega_X + \omega_\Sigma$ be a K\"ahler form on $M$, where $\omega_\Sigma = {i\over 2}\ddbar \rho.$ Let $\pi_X: M\rightarrow X$ and $\pi_\Sigma: M \rightarrow \Sigma$ be the natural projection maps.

Let $\a$ be closed real $(1,1)$ form on $M$ with bounded coefficients, with
$\alpha_s\equiv\alpha|_{\{s\}\times X} \geq \e\o_X$ for some $\e>0$ and all $s\in \Sigma$.
Let $f$ be a continuous function on $\overline{M}$ such that for each $s\in \p\Sigma$, $f|_{X_s} \in \PSH(X, \alpha_s)$, where $X_s=\pi_X^{-1}(s)$. 
We define
\bea
\E_{M, \alpha, f} = \{ u ~|~ u\in \PSH(M, \alpha)\cap C(M),~ u\leq f~\textnormal{on}~\partial M\}.
\eea
and the upper envelope of $\E_{M, \alpha, f}$ by
\bea
\vp = \sup_{\E_{M, \alpha, f}} u.
\eea
The following theorem was proved in \cite{BD};

\begin{theorem}
\label{BD}
Let $(X,\o_X)$, $\Sigma$, $\alpha$, and $f\in C(\p M)$ satisfy all the properties listed above, and define the family $\E_{M,\alpha,f}$ and its upper envelope $\vp$ as above. Then the function $\vp$ is the unique $\alpha$-plurisubharmonic solution of the Dirichlet problem
\bea
(\alpha+{i\over 2}\ddb\vp)^{{\rm dim}\,M}=0
\ \ {\rm on}\ \ M,
\qquad
\vp=f\ \ {\rm on}\ \ \p M.
\eea
Furthermore, if $f\in C^{1,1}(\p M)$, then for any $s\in\Sigma$, $i\ddb \vp|_{\{s\}\times X}$ is locally bounded,
uniformly in $s\in \Sigma$.
\end{theorem}

The proof in \cite{BD} depends on the Kiselman infimum principle and refined regularization techniques for plurisubharmonic functions. Here we discuss only the special case when $\alpha=\o_X$ (more precisely, the pull-back of $\o_X$ to $M$) and $f\in C^\infty(\p M)$. In this case, the fiber-wise regularity can be obtained by the method of elliptic regularization used earlier for the geodesic equation and the standard $C^2$ estimates of Yau, as described in \S 6. 

\smallskip
First we show that $\vp$ is the solution of the Dirichlet problem.
Let $\omega_{X, f} = \omega_X + \frac{i}{2} \ddbar f$ after extending $f$ to a smooth function on $\overline{M}$. The original problem is equivalent to the same problem formulated rather with
\bea
\E_{M, \omega_{X,f}, 0} = \{ u ~|~ u\in \PSH(M, \omega_{X,f}),~ u\leq 0~\textnormal{on}~\partial M\},
\qquad
\vp= \sup_{\E_{M, \omega_{X,f}, 0}} u.
\eea

We show that $\vp$ is continuous on $\partial M$. 
The form $\omega_{X, f} + A\frac{i}{2} \ddbar\rho$ is a K\"ahler form on $M$ for sufficiently large $A>0$. Hence $A\rho \in \E_{M, \omega_{X,f}, 0}$ and then  $A \rho \leq \vp.$ 
On the other hand, for any $u\in \E_{M, \omega_{X, f}, 0}$,  $u+ A\rho $ is plurisubharmonic on $\pi_X ^{-1}(z)$ for each $z\in X$, for sufficiently large $A>0$ independent of the choice $z\in X$. Thus $ u + A \rho \leq 0$ since $u+ A\rho \leq 0$ on $\partial (\pi_X^{-1}(z))$. It easily follows that 
\bea
B \rho \leq  \varphi \leq - B\rho
\eea
for some $B>0$. In particular, $\vp$ is continuous on $\partial M$.

\smallskip

Next we show that $\vp$ is continuous in $M$. First we fix $A>0$ with 
\bea
A\rho \leq \varphi \leq -A \rho.
\eea
For any compact subset $K$ in $M$ and any sufficiently small $\epsilon>0$, we choose $\delta = (4A)^{-1}\epsilon$ so that $K\subset M_{4\delta}$ and $\vp \leq \epsilon/4$ on $M\setminus M_{4 \delta}$, where 
\bea
M_\delta =\{ (z,s)\in M~|~ \rho(z,s)< -\delta\}.
\eea 
By Demailly's regularization techniques \cite{D}, there exists a decreasing sequence $\{ u_j\} \subset \PSH(M, \left(\frac{j+1}{j}\right)  \omega_{X,f})\cap C(M_{\delta/2})$ which converges to the $\vp^*$, the upper semi-continuous envelope of $\vp$. We define
\bea
\tilde u_j (z,s)= \{ \begin{array}{cc}
                 \max ( \frac{j}{j+1} u_j - \epsilon, 2A \rho ), &  (z,s) \in M_\delta\\
                 2A\rho & (z,s)\notin M_\delta.
               \end{array}
\eea
On $\partial M_\delta$, $\frac{j}{j+1} u_j - \epsilon \leq -\epsilon/2 \leq -2 A \delta.$ Hence  $\tilde u_j \in \E_{M, \omega_{X,f}, 0}$. Furthermore, on $K$, 
\bea
\frac{j}{j+1} u_j - \epsilon \geq \frac{j}{j+1} \vp - \epsilon  \geq \frac{j}{j+1}A\rho - 4A\delta \geq  A\rho + A\rho = 2A\rho  ,
\eea 
and so $\tilde u_j = \frac{j}{j+1} u_j - \epsilon$.
It follows immediately that
\bea
\vp^* \leq u_j = \frac{j+1}{j} ( \tilde u_j +\epsilon) \leq \frac{j+1}{j} \vp + \frac{j+1}{j}\epsilon
\eea 
or $0\leq u_j-\vp \leq \frac{1}{j} |\vp| + \epsilon$.
Therefore $ u_j$ converges to $\vp$ uniformly in $L^\infty(K)$ and so $\vp^*$ is continuous in $K$. In conclusion, $\vp= \vp^* \in \E_{M, \omega_{X,f}, 0}\cap C(\overline{M})$.

\smallskip
Finally, we show that $(\frac{i}{2}\ddbar \hat u)^{{\rm dim}\,M} =0$. Fix any Euclidean ball $B$ in $M$. There exists $\eta \in C^\infty(B)$ such that $\omega_{X, f} =\frac{i}{2}\ddbar \eta$ on $B$. Let 
\bea
\E_{B, (\vp+\eta)|_{\partial B}}=\{v~|~\varphi\in \PSH(B), ~v|_{\partial B} = (\vp+\eta)|_{\partial B} \}
\eea
and  $\psi = \sup_{v\in\E_{B, (\vp +\eta)|_{\partial B}}} v$.
Then by the above Bedford-Taylor theorem for ${\bf C}^{{\rm dim}\,M}$, 
we have $( \frac{i}{2} \ddbar \psi)^{{\rm dim}\,M}=0$ and $\vp = \psi - \eta$ on $\overline{B}$. Hence $(\omega_{X, f} + \frac{i}{2} \ddbar \vp)^{{\rm dim}\,M}=0$ on $B$.

\medskip

We turn to the proof of fiber wise regularity.
Consider the following elliptic regularization of the homogeneous Monge-Amp\`ere equation,
\begin{equation}
 (\omega_X +  \frac{i}{2} \ddbar \vp_\e)^{n+m}= \e\, \omega_M^{n+m} \ \textnormal{on}~ M, 
 \qquad
 \vp_\epsilon=f  \  \textnormal{on}~ \partial M.
\end{equation}
For any $\e>0$, there exists a unique smooth solution $\vp_\epsilon$ in $\overline M$. In fact, since
$(\omega_X +  \frac{i}{2} \ddbar \underline{u}_s)^{{\rm dim}\,M} > \epsilon \omega_M^{{\rm dim}\,M}$, the function
 $\underline{u}_s\equiv f+A\rho$ is a subsolution satisfying $\underline{u}_s|_{\p M}=f|_{\p M}$
 for sufficiently large $A>0$. Thus $\vp_\e\geq \underline{u}_\e$.
 Furthermore $\vp_\e\leq \vp\in \E_{M,\o_X,f}$, and we have
 \bea
 \|\vp_\e\|_{C^0(M)}\leq C
 \eea
 uniformly in $\e$.

It suffices now to show that
there exists $C>0$ such that for any $\epsilon\in (0,1)$ and $s\in \overline\Sigma$,
\bea
\label{2nd order}
(\omega_X +  \frac{i}{2} \ddbar \vp_\e)|_{X_s} \leq C \omega_X |_{X_s}.
\eea
It would follow then that, for any $s\in \overline\Sigma$, 
\bea
(\omega_X +  \frac{i}{2} \ddbar \vp )|_{X_s} \leq C \omega_X.
\eea
since $\vp_\epsilon$ is increasingly monotone as $\e \rightarrow 0$ and $\vp_\e$ converges to $\vp$ uniformly. 

To establish (\ref{2nd order}),
let $\omega_\e = \omega_X +  \frac{i}{2} \ddbar \vp_\e$. We denote $g_X$, $g_\Sigma$, $g_M$ and $g_\e$ be the K\"ahler metrics associated to $\omega_X$, $\omega_\Sigma$, $\omega_M$ and $\omega_\e$. We always use the product coordinates for $M$, where $(z, s) = (z_1, ..., z_n, s_1, ..., s_m)$ and $z\in X$ and write $\omega_M = {i\over 2} (g_X)_{\bar j i} dz^i \wedge d\bar z^j + {i\over 2} (g_\Sigma)_{\bar \beta\alpha} ds^\alpha \wedge d\bar s^\beta. $

We define 
\bea
H_\epsilon= \frac{ (\omega_X +  \frac{i}{2}\ddbar \vp_\e)\wedge \omega_X^{n-1} \wedge \omega_\Sigma^m}{\omega_X^n \wedge \omega_\Sigma^m}.
\eea
Notice that $H_\e$ is the trace of the relative endomorphism between
$\o_\e|_{X_s}$ and $\o_X$,
\bea
H_\e (z,s) = {\rm Tr}_{\omega_X}(\omega_\epsilon|_{X_s})(z) = \sum_{i, j=1}^n (g_X)^{i \bar j} (g_\e)_{\bar j i},
\eea 
where $g_X$ is the K\"ahler metric associated to $\omega_X$ and $g_\e$ is the K\"ahler metric associated to $\omega_X + {i\over 2}\ddbar \vp_\e$. The same calculations as in Yau's Schwarz lemma
\cite{Y78a} or in Yau's second order estimates, see \S 6, show that there exists $C_1>0$ such that for all $\e \in (0,1)$, we have on $M$,
\bea
\Delta_\e \log H_\e \geq - C_1 {\rm Tr}_{\omega_\e}(\omega_M)- C_1, 
\eea 
where $\Delta_\e$ is the Laplacian operator on $M$ with respect to $\omega_\e$.
We also have 
\bea
\Delta_\e (-\vp_\e+\rho) = -(n+m) + {\rm Tr}_{\omega_\e}(\omega_X + \omega_\Sigma)= -(n+m) + {\rm Tr}_{\omega_\e}(\omega_M).
\eea
Then there exist $C_2, C_3, C_4>0$ which are independent of $\epsilon\in (0,1)$ such that 
 \begin{eqnarray*}
\Delta_\e (\log H_\e -A \varphi_\e + A\rho) 
 &\geq&  (A-C_1) {\rm Tr}_{\omega_\e}( \omega_M)- C_1 - A(n+m)\nonumber\\
 &\geq& C_2 ({\rm Tr}_{\omega_M}(\omega_\e))^{1/(n+m-1)} \left(\frac{\omega_M^{n+m}}{\omega_\e^{n+m}}\right)^{1/(n+m-1)} - C_2\nonumber\\
 &\geq& C_3 H_\e - C_4.
 \end{eqnarray*}
On the other hand, 
\bea
H_\e|_{\partial M} = \frac{ (\omega_X +  \frac{i}{2}\ddbar f)\wedge \omega_X^{n-1}}{\omega_X^n}|_{\partial M}
\eea 
is uniformly bounded from above for all $\e\in (0,1)$.
Applying the maximum principle, we obtain a constant $C_5>0$ such that for all $\e\in (0,1)$, 
\bea
H_\epsilon \leq C_5
\eea 
since both $u_\e$ and $\rho$ are uniformly bounded in $C(\overline{M})$. The proof of the fiberwise regularity is complete.

\subsection{Envelopes with integral conditions}

In the previous section, we have seen how envelopes with pointwise Dirichlet conditions can produce solutions to the Dirichlet problem for complex Monge-Amp\`ere equation. It would be interesting to determine whether envelopes with integral conditions can be effectively used to produce other solutions of Monge-Amp\`ere equations, or other canonical metrics. We describe some examples of such envelopes below.

\medskip
One example is the following hermitian metric defined by Tsuji \cite{Ts07} on projective manifolds of general type, generalizing the metric introduced in \cite{NS}.
Let $X$ be a smooth projective variety of general type. Fix a smooth hermitian metric $h_0$ on $K_X$
and define
\begin{equation}
\varphi_{can} (z) = \sup\{ \varphi(z) ~|~ Ric(h_0) + \frac{i}{2} \ddbar \varphi \geq 0, ~~ \int_X e^\varphi h_0^{-1}=1\} \end{equation}
and
\begin{equation}
h_{can} = e^{-\varphi_{can}} h_0.
\end{equation}
It has been shown by Berman and Demailly \cite{BD} that this metric $h_{can}$ coincides with the
metric $\tilde h_{can}$ defined instead by
\begin{equation}
\tilde h_{can} (z)= \inf_{m\in {\bf Z}^+} \inf\{  (\sigma\wedge \bar \sigma(z) )^{-1/m}~|~ \int_X |\sigma\wedge\bar\sigma|^{1/m} =1,~ \sigma \in H^0(X, mK_X)\}.
\end{equation}
which is manifestly a birational invariant, 
$\tilde h_{can}$ is a birational invariant
since $H^0(X, mK_X)$ is invariant under birational transformations.

\medskip

Another example is the following. 
Let $\Omega$ be a bounded strictly pseudoconvex domain in ${\bf C}^n$. We define
\bea
\varphi_{can} (z) = \sup \{ \varphi(z) ~|~ \frac{i}{2} \ddbar \varphi\geq 0,~ \int_\Omega e^{\varphi} (\frac{i}{2}\ddbar |z|^2)^n=1\}.
\eea
We also define $(h_{can})^{-1} = e^{\varphi_{can}} ( \frac{i}{2}\ddbar |z|^2)^n$ to be the canonical measure on $\Omega$.

\begin{lemma}
$\varphi_{can}$ is a plurisubharmonic function on $\Omega$.
\end{lemma}

\noindent {\it Proof.}  First, we show that $\varphi_{can}$ is bounded from above in any compact subset of $\Omega$. Suppose not, then by taking a subsequence, there exist a sequence of points $z_j \rightarrow \hat z \in \Omega$ and a sequence of psh functions $\varphi_{j}$ with $\int_{\Omega} e^{\varphi_j} =1$ such that $$\varphi_j(z_j) \rightarrow \infty.$$
Without loss of generality, we can assume that $B(z_j, r)\subset \subset\Omega$ for all $j$ for some fixed $r>0$. Then by the mean value inequality and Jensen's inequality, there exist positive constants $C_1,~C_2,~C_3$ independent on $j$ such that
\bea
1= \int_{\Omega} e^{\varphi_j} \geq C_1 e^ { \int_{\Omega} \varphi_j } \geq C_2 e^ {\int_{B(z_j, r)} \varphi_j} \geq C_3 e^{\varphi_j(z_j)} \rightarrow \infty.
\eea
This is a contradiction.

\smallskip
Next we have to show that $\varphi_{can} = (\varphi_{can})^*$. By the definition of $\varphi_{can}$, for any $\hat z\in \Omega$, there exists a sequence $z_j\rightarrow \hat z$ and psh $\varphi_j$ with $\int_\Omega e^{\varphi_j}=1$ such that
\bea
\varphi_j(z_j) \rightarrow (\varphi_{can})^*(\hat z).
\eea
By taking a subsequence, we can assume that $\varphi_j$ converges to a psh function $\varphi$ in $L^1(\Omega)$ and thus almost everywhere. In particular, $\int_\Omega e^\varphi \leq 1$ by Fatou's lemma. On the other hand,
\begin{eqnarray*}\varphi_{can}(\hat z) &\geq& \varphi(\hat z) = \lim_{r\rightarrow 0} \frac{1}{vol(B(\hat z, r))} \int_{B(\hat z, r)} \varphi = \lim_{r\rightarrow 0} \lim_{j\rightarrow \infty} \frac{1}{vol(B(z_j, r))} \int_{B( z_j, r)} \varphi_j \\ &\geq& \lim_{j\rightarrow \infty} \varphi_j (z_j)= (\varphi_{can})^*(\hat z).\end{eqnarray*}
The lemma is proved.   Q.E.D.

\smallskip
The following theorem provides a more algebraic characterization of $\vp_{can}$:

\begin{theorem} 
$\varphi_{can} \in \PSH(\Omega) \cap C(\Omega)$ and 
\begin{equation}
\label{vpcan}
\varphi_{can}(z) = \sup_{m\in {\bf Z}^+} \sup \{ \frac{1}{m} \log |f|^2(z)~|~ f\in \cO(\Omega), ~~\int_\Omega |f|^{\frac{2}{m}} (\frac{i}{2} \ddbar |z|^2)^n=1\}.
\end{equation}

\end{theorem}

\noindent {\it Proof.} Denote by $\varphi_{can, alg}$ the right-hand side of
(\ref{vpcan}).
It is easy to see that $\varphi_{can,alg} \leq \varphi_{can}$,
so it suffices to show that $\varphi_{can, alg} \geq \varphi_{can}$.

\smallskip

By the Ohsawa-Takegoshi extension theorem, for any psh $\varphi$ on $\Omega$ with $\int_\Omega e^\varphi =1$ and any point $z \in \Omega$, there exists a holomorphic function $f$ on $\Omega$ such that for any $m\in {\bf Z}^+$, 
\bea
|f| ^2 e^{-(m-1)\varphi} (z) =1
\quad
{\rm and}
\quad
\int_\Omega |f|^2 e^{-(m-1)\varphi} \leq C,
\eea
where $C$ does not depend on $\varphi$ or $m$. By H\"older's inequality, we have
\begin{eqnarray*}  
\int_\Omega |f|^{2/m}\leq \left(   \int_\Omega |f|^2 e^{-(m-1)\varphi} \right)^{1/m} \left( \int_{\Omega} e^\varphi \right)^{(m-1)/m}
\leq (C)^{1/m}.
\end{eqnarray*}
Let $F= \frac{f}{ ( \int_\Omega |f|^{2/m})^m}$. Then $F\in \cO(\Omega)$ with $\int_\Omega |F|^{2/m}=1$, and
\bea
\frac{1}{m} \log |F|^2(z) \geq (m-1) \varphi(z) - \frac{1}{m} C = \varphi(z) - \frac{1}{m} ( \varphi(z) + C).
\eea
For fixed $z$, $\varphi(z)$ and $C$ are uniformly bounded from below. By letting $m\rightarrow \infty$, we get
\bea
\varphi_{can, alg}(z) \geq \varphi(z).
\eea
Since this is true for any $z$ and any psh $\varphi$ with $\int_\Omega e^\varphi =1$, we have
$\varphi_{can, alg}(z) \geq \varphi_{can} (z)$,
and hence $\varphi_{can, alg} = \varphi_{can}$.

\smallskip

Now we can show that $\varphi_{can}$ is continuous. It suffices to show that $\varphi_{can}$ is lower semi-continuous. Suppose not. Then there exists $\hat z\in \Omega$ and $\e>0$ and a sequence of points $z_j\in \Omega$ converging to $\hat z$ so that
\bea
\varphi_{can}(z_j)< \varphi_{can}(\hat z) - \epsilon.
\eea
Also there exist $f\in \cO(\Omega)$ and $m\in{\bf Z}^+$ such that 
\bea
\int_\Omega |f|^{2/m} = 1, ~~~~\frac{1}{m}\log |f|^2(\hat z) > \varphi_{can}(\hat z) - \frac{\epsilon}{2}.
\eea
Then there exist $r>0$ such that for all $z\in B(\hat z, r)$, $$\varphi_{can}(z) \geq \frac{1}{m}\log |f|^2(\hat z) > \varphi_{can}(\hat z) - \frac{\epsilon}{4}.$$ This is a contradiction. Q.E.D.

\medskip
A natural question to ask is whether $g_{can} = \frac{i}{2} \ddbar \varphi_{can}$ defines a complete metric on a bounded strictly pseudoconvex domain in ${\C}^n$ and how it is related to other invariant metrics such as the Bergman,  Carath\'eodory and Kobayashi metrics.

\section{Further Developments}
\setcounter{equation}{0}

As we had acknowledged in the introduction, we could not cover all the possible recent developments, and this survey has not touched on many important topics.
In this section, we would like to mention a few and provide some references, for readers who may be completely new to the subject.

\medskip
A first major omission is a discussion of the important
equation (\ref{CY}) with $F(z,\vp)=e^{f(z) - \vp}$, which corresponds to the open problem of K\"ahler-Einstein metrics on a compact K\"ahler manifold $(X,\o_0)$ with $\o_0\in c_1(K_X^{-1})$. From the discussion of a priori estimates in \S 6 and \S 7, we see that the equation would be solvable if we can obtain a $C^0$ estimate. The problem is to link such an estimate to stability in GIT, as required by the conjecture of Yau \cite{Y93}. Donaldson has recently laid out a program for achieving this \cite{D10, D11a, D11b, CDa, CDb}. Prior to this program,
K\"ahler-Einstein metrics with positive scalar curvature have 
been found in various geometric situations by Tian and Yau \cite{TY87}
and Tian \cite{T87} using the $\alpha$-invariant,
and by Siu \cite{Si} and Nadel \cite{N} using multiplier ideal sheaves.
Necessary conditions for K\"ahler-Einstein metrics have been
obtained by Tian \cite{T97}.  
It has been shown by Tian \cite{T90b} that, for surfaces, the existence of K\"ahler-Einstein metrics is equivalent to the vanishing of the Futaki invariant. A full account of the arguments in \cite{T90b} can be found in the paper of Tosatti \cite{T10}. The same characterization of the existence of K\"ahler-Einstein metrics by the vanishing of the Futaki invariant has been established by Wang and Zhu \cite{WZ2}. Their proof exploits the fact that the toric potentials on a toric variety satisfy a real Monge-Amp\`ere equation, and the image of their gradients is the polytope of the variety. For a survey of some of these developments and the related issue of stability, see \cite{PS03, PS08}.

\medskip
The related question of singularities for the Monge-Amp\`ere equation when the manifold is unstable is of similar considerable interest, and even less explored. The case of holomorphic vector bundles has seen remarkable progress, with the recent works of G. Daskalopoulos and R. Wentworth for complex surfaces \cite{DW1, DW2},
and A. Jacob \cite{J1,J2} for general K\"ahler manifolds, on the generalization 
to arbitrary dimensions of the Atiyah-Bott formula for complex curves. An analysis of a break-up of an unstable ruled surface by the Calabi flow has been given by G. Szekelyhidi \cite{Sz}.

\medskip
Another major omission is parabolic complex Monge-Amp\`ere equations, and particularly the K\"ahler-Ricci flow. As we had mentioned earlier, starting with the papers of Cao \cite{Cao} and Tsuji \cite{Ts}, there has been a constant feedback between developments for the elliptic and for the parabolic 
Monge-Amp\`ere equation. In fact, much of the material discussed in Sections 3, 6, 7 either arose from or are directly motivated by the study of the K\"ahler-Ricci flow on Fano manifolds
(e.g. \cite{ST09, PS06a, PSSW1, CZ, MS, Yu, Zh} and references therein) or on manifolds of general type (e.g. \cite{ST06,
TZ, ST09,SW1}). New powerful techniques have been introduced by Perelman
(see \cite{ST} for an account of Perelman's unpublished results on the K\"ahler-Ricci flow).
We refer to the papers we listed as well as to the recent survey
\cite{SW2} for a fuller list of references.
Related developments for the Sasaki-Ricci flow can be found in \cite{Co1, Co2, He1}.

\medskip
The Monge-Amp\`ere measure is uniquely defined by Bedford and Taylor \cite{BT76} for locally bounded potentials. The largest classes of possibly unbounded potentials for which a well-behaved measure can be defined have been identified by Cegrell \cite{Ceg} and Blocki \cite{B06}. The Monge-Amp\`ere measures can also be defined for unbounded potentials, if their singularity set is relatively compact within Stein neighborhoods \cite{D, Sib}. A prime example is the pluricomplex Green's function (see e.g.
\cite{L, BD, Gb, B00, Ze1} and references therein).
As we saw in Section \S 9, a non-pluripolar definition can be given, and the range of the corresponding Monge-Amp\`ere measures has been completely characterized by Guedj and Zeriahi \cite{GZ}. The investigation of Monge-Amp\`ere measures which charge pluripolar sets is still in its infancy, see \cite{CG} for examples on projective spaces and \cite{ACCH} for some general results. It is an important direction for research.

\medskip
In Section \S 13, we have seen how geodesics in the space of K\"ahler metrics lead to the homogeneous complex Monge-Amp\`ere equation. Similarly, Donaldson \cite{D07} has shown how geodesics in the space of volume forms on a Riemannian manifold lead to a non-linear equation now known as Donaldson's equation. He also showed how this equation can be interpreted as a PDE version of Nahm's equation in mathematical physics, and is closely related to well-known free boundary problems. The existence of $C^{1,\alpha}$ solutions of Donaldson's equation has been obtained by Chen and He \cite{CH} and He \cite{H}. The same questions of regularity and maximum rank arise for this equation as they do for the homogeneous complex Monge-Amp\`ere equation. Some early results in low dimensions can be found in \cite{GPa, GPb}. The existence and regularity of geodesics in the space of Sasaki metrics have also been investigated in \cite{GZ}.

\medskip
We have seen in Section 3.2 how the most basic Alexandrov-Bakelman-Pucci estimates can be applied to the complex Monge-Amp\`ere equations. It would be interesting to find out whether this method can be carried out further. In this context, we would like to mention the recent remarkable ABP estimates on Riemannian manifolds obtained by Wang and Zhang \cite{WZ1}, building on earlier works of Cabr\'e \cite{Ca1}.

\medskip

Finally, we would like to mention viscosity methods. They have been very successful in the investigation of non-linear equations where no complex structure plays a particular role \cite{CIL}. Even though the notion of plurisubharmonicity poses a number of difficulties, it may not be unreasonable to expect that viscosity methods can be developed and become of wider use for equations such as the complex Monge-Amp\`ere equation. Some major steps in incorporating plurisubharmonicity in viscosity methods have been undertaken by Harvey and Lawson \cite{HL1,HL2}, Eyssidieux, Guedj, and Zeriahi \cite{EGZ10}, and Wang \cite{W1}.
For example, a version of Theorem 12, establishing the existence and uniqueness of viscosity solutions to the Dirichlet problem for the equation (\ref{CY}) 
on domains in ${\bf C}^n$ for continuous data, has been established in \cite{W1}.

\bigskip

\begin{appendix}

\section{Plurisubharmonic functions}
\setcounter{equation}{0}

We gather here for the convenience of the reader some basic properties of plurisubharmonic functions and of their Monge-Amp\`ere measures.

\subsection{The exponential estimate}

Now the $L^2$ norm, in fact the $L^p$ norm for any $p<\infty$,
of any non-positive plurisubharmonic function $\vp$ is bounded by a constant depending
only on the K\"ahler class of $\o_0$. This is a consequence of
the following local estimate of H\"ormander, extended to K\"ahler manifolds by 
\cite{T87, TY87, Ze}

\begin{theorem}
\label{Hormander}
Let $\o_0$ be a K\"ahler form.
There exists a a strictly positive number $\a$ and
a constant $C$ depending only on $\o_0$  
so that
\bea
\label{Hormander}
{1\over [\o_0^n]}\int_X e^{-\a(u-{\rm sup}_Xu)}\o_0^n
\leq C
\eea
for all $u\in \PSH(X,\o_0)$.
\end{theorem}

Since $e^{\alpha t}\geq ({\a\over p})^p t^p$ for all $t\geq 0$ and all $p>0$,
it follows that $\|u-{\rm sup}_Xu\|_{L^p}$ is bounded by a constant depending
only on $p$ and $\o_0$.

For recent advances on exponential estimates for plurisubharmonic functions, see \cite{DNS}.

\subsection{Regularization of plurisubharmonic functions}

The existence of approximations of $\o_0$-plurisubharmonic functions by monotone sequences of smooth $\o_0$-plurisubharmonic functions is much more delicate
for K\"ahler manifolds than for domains in ${\bf C}^n$. Part of the difficulty
resides in the conflicting roles of the differential geometric and the complex structure. An early approximation theorem with loss of $\e$-positivity is due to Demailly \cite{D89}. Many others are now available, including the recent ones
of Demailly, Peternell, and Schneider
\cite{DPS} and of Demailly and Paun \cite{DP2}, which imply in particular the following statement:
let $(X,\o_0)$ be a K\"ahler manifold, and let $\gamma$ be a continuous non-negative $(1,1)$-form. Then for any $\vp\in PSH(X,\gamma)$ with Lelong numbers
$\nu_\vp(z)=0$ for all $z\in X$, and any 
subset $X'\subset X$ with compact closure, there exists a decreasing sequence $\e_j\downarrow 0$, and a sequence $\vp_j\in \PSH(X,\gamma+\e_j\o_0)\cap C^\infty(X)$
with $\vp_j\downarrow\vp$ in a neighborhood of $X'$.
This statement was given an independent proof by Blocki and Kolodziej \cite{BK}.
But particularly important for our purposes is the observation 
of Blocki and Kolodziej \cite{BK} that, for $X$ compact and $\gamma=\o_0$, no loss of positivity is necessary:

\begin{theorem}
\label{BK}
Let $(X,\o_0)$ be a compact K\"ahler manifold. Then for every $\vp\in PSH(X,\o_0)$, there exists a sequence $\vp_j\in PSH(X,\o_0)\cap C^\infty(X)$
with $\vp_j\downarrow \vp$.
\end{theorem}

For the convenience of the reader, we provide some details on how
to derive Theorem \ref{BK} from the result of \cite{DPS, DP2}.

\medskip

Let $\vp\in PSH(X,\o_0)$, $\vp\leq -1$ on $X$. For each $j$, the function ${\rm max}(\vp,-j)$ is bounded, and hance has vanishing Lelong numbers.
Thus the result of \cite{DPS, DP2} implies the existence of a sequence 
of smooth functions $\psi_{jk}\in PSH(X,(1+\e_{j,k})\o_0)$ with
$\psi_{j,k}\downarrow {\rm max}(\vp,-j)$ and $\e_{j,k}\downarrow 0$
as $k\to\infty$. By passing to a subsequence,
we may assume that $\e_{j,k}\geq \e_{j+1,k}$ for all $j,k$. 
It suffices to show that
there exists a sequence $k_1<k_2<\cdots$ so that
$\tilde \vp_j\equiv \psi_{j,k_j}+{1\over 2^j}$ is a decreasing sequence converging to $\vp$. The sequence
\bea
\vp_j\equiv {\tilde\vp_j\over 1+\e_{j,k_j}}
\eea
is then a sequence of smooth functions in $PSH(X,\o_0)$ with $\vp_j\downarrow \vp$.

\medskip
We choose $k_j$ inductively as follows. Fix $j$. Let
\bea
C_j=\{\psi_{j+1,k}\geq \psi_{l,k_l+j}+{1\over 2^{j+1}} \ {\rm for\ some}\ l\leq j\}.
\eea
Since $\cap_k C_k=\emptyset$, we can choose $k_{j+1}$ so that $C_{k_{j+1}}=\emptyset$. Note that 
\bea
\psi_{j+1,k_j}\leq \psi_{j,k_j}+{1\over 2^{j+1}}
\ \ {\rm and\ so}
\ \ \tilde\vp_{j+1}\leq\tilde\vp_j.
\eea
Now fix $x\in X$, $\e>0$, and assume that $\vp(x)\geq -j_0$. Choose $j$ so that
\bea
0\leq \psi_{j_0,k_{j_0}+j}(x)-\vp(x)<\e
\eea
Then
$\psi_{j,k_j}-\vp(x)<\e+{1\over 2^j}$, which implies
$\tilde\vp_j(x)-\vp(x)<\e+{1\over 2^j}+{1\over 2^j}$. The argument is complete.

\subsection{The comparison principle}

The following is a useful version of the comparison principle. It follows from
the standard arguments of Bedford and Taylor \cite{BT82}, using the above approximation theorem for plurisubharmonic functions on K\"ahler manifolds.

\begin{theorem}
\label{comparison}
Let $(X,\o_0)$ be a compact K\"ahler manifold with smooth boundary $\p X$
and dimension $n$, and let $\o$ be a smooth, non-negative, closed $(1,1)$-form.
Then we have
\bea
\int_{\{\vp<\psi\}}(\o+{i\over 2}\ddb \psi)^n
\leq
\int_{\{\vp<\psi\}}(\o+{i\over 2}\ddb\vp)^n
\eea
for all $\vp,\psi\in \PSH(X,\o)\cap L^\infty(X)$ satisfying
${\rm liminf}_{z\to\p X}(\vp(z)-\psi(z))\geq 0$.
\end{theorem}

\end{appendix}

\bigskip

\bigskip
{\bf Acknowledgements}

\smallskip

The first-named author would like to thank Professors Huai-Dong Cao
and Xiaofeng Sun of Lehigh University for their very warm hospitality during
his visit there during the 2010 conference in honor of Professor C.C. Hsiung. The authors would like to thank Professor Pengfei Guan for providing them with some invaluable references. They have also benefitted greatly from interactions with Professors O. Munteanu, G. Szekelyhidi, V. Tosatti, Ben Weinkove, and the participants of the informal Complex Geometry and PDE seminar at Columbia University.

\end{document}